\documentstyle[12pt,diagram]{article}
\newtheorem{thm}{Theorem}
\newtheorem{prop}[thm]{Proposition}

\newtheorem{thm-defi}[thm]{Theorem/Definition}
\newtheorem{example}[thm]{Example}
\newtheorem{cor}[thm]{Corollary}

\newtheorem{new-lemma}[thm]{Lemma}
\newtheorem{defi}[thm]{Definition}
\newtheorem{rem}[thm]{Remark}

\newtheorem{condition}{Condition}

\newcommand{\A}{{\cal A}}
\newcommand{\C}{{\cal C}}

\newcommand{\E}{{\cal E}}
\newcommand{\F}{{\cal F}}

\newcommand{\Hy}{{\cal H}}

\newcommand{\LB}{{\cal L}}

\newcommand{\M}{{\cal M}}
\newcommand{\U}{{\cal U}}
\newcommand{\V}{{\cal V}}
\newcommand{\W}{{\cal W}}

\newcommand{\PP}{{\Bbb P}}

\newcommand{\Integers}{{\Bbb Z}}
\newcommand{\ComplexNumbers}{{\Bbb C}}

\newcommand{\LieAlg}[1]{{\frak #1}}
\newcommand{\linsys}[1]{{\mid}#1{\mid}}

\newcommand{\IsomRightArrowOf}[1]{
\stackrel
{\stackrel{#1}{\cong}}
{\rightarrow}
}
\newcommand{\IsomRightArrow}{\stackrel{\cong}{\rightarrow}}
\newcommand{\LongIsomRightArrow}{\stackrel{\cong}{\longrightarrow}}
\newcommand{\RightArrowOf}[1]{\stackrel{#1}{\rightarrow}}

\newcommand{\LongLeftArrowOf}[1]{\stackrel{#1}{\longleftarrow}}

\newcommand{\HookRightArrowOf}[1]{\stackrel{#1}{\hookrightarrow}}
\newcommand{\LongRightArrowOf}[1]{\stackrel{#1}{\longrightarrow}}

\newcommand{\StructureSheaf}[1]{{\cal O}_{#1}}
\newcommand{\EndProof}{\hfill  $\Box$}
\newcommand{\restricted}[2]{#1_{\mid_{#2}}}

\newcommand{\Hilb}{\rm Hilb}

\newcommand{\rank}{\rm rank}
\newcommand{\coker}{\rm coker}
\newcommand{\Pic}{\rm Pic}
\newcommand{\Sym}{\rm Sym}
\newcommand{\Ext}{\rm Ext}

\newcommand{\Hom}{\rm Hom}
\newcommand{\Aut}{\rm Aut}
\newcommand{\End}{\rm End}
\newcommand{\Abs}[1]{\mid\!#1\!\mid}

\newcommand{\SheafHom}{{\cal H}om}
\newcommand{\SheafEnd}{{\cal E}nd}
\newcommand{\SheafExt}{{\cal E}xt}
\newcommand{\RelExt}{{\cal E}xt}
\newcommand{\Ideal}[1]{{\cal I}_{#1}}

\newcommand{\Wedge}[1]{\stackrel{#1}{\wedge}}
\newcommand{\Contract}{\rfloor}

\newcommand{\DoubleTilde}[1]{\stackrel{\approx}{#1}}
\input amssymb.sty
\input psfig.sty

\begin{document}
\begin{center}
\begin{Large}
{\bf 
\noindent
 Brill-Noether duality for moduli spaces of sheaves on K3 surfaces
}
\end{Large}
\\
Eyal Markman
\footnote{Partially supported by NSF grant number DMS-9802532}
\end{center}

{\scriptsize
\tableofcontents 
}

\section{Introduction}
\label{sec-introduction}

Hilbert schemes of points are among the simplest moduli spaces of sheaves 
on an algebraic surface $S$. Compactified relative Picards 
over a linear system of curves in $S$ may be considered as moduli space 
of sheaves with pure one-dimensional support. The latter 
is complicated in comparison. 
For example, the punctual Hilbert schemes are always smooth
while the compactified relative Picard is rarely smooth. 
This work was motivated by a desire to understand (resolve) the 
birational Abel-Jacobi isomorphism 
\[
S^{[g]} \ \ \leftarrow\cdots\rightarrow \ \ 
\Pic^g_{\linsys{\StructureSheaf{S}(1)}}
\]
between the punctual Hilbert scheme 
$S^{[g]}$ of a symplectic surface with a linear system 
$\linsys{\StructureSheaf{S}(1)}$ of curves of genus $g$ and 
the relative compactified Picard $\Pic^g_{\linsys{\StructureSheaf{S}(1)}}$. 
The latter is a completely integrable hamiltonian system 
\cite{hurtubise-local-geometry,mukai-symplectic-structure}. 
Of particular interest is the case where $S$ is the cotangent bundle 
of a smooth algebraic curve and 
$\Pic_{\linsys{\StructureSheaf{S}(1)}}$ ``is'' a Hitchin system.  
The recursive nature of the geometry involved forced us to consider 
a more general duality among the moduli spaces of sheaves of arbitrary rank
on the symplectic surface. However, due to the complexity of the birational 
isomorphisms, we consider here only the simplest setup; that of 
a K3 surface $S$ and the collection of moduli spaces of sheaves whose 
first Chern class is minimal (Condition \ref{cond-linear-system} 
Section \ref{sec-admissible-collections}). 
The structure, both local and global, of the birational Abel-Jacobi 
isomorphism turns out to be surprisingly beautiful. 
Local analytically, the birational isomorphism is modeled after 
two dual Springer resolutions 
\[
T^*G(t,H) \ \rightarrow  \overline{{\cal N}^t} \leftarrow T^*G(t,H^*)
\]
of the closure of the nilpotent orbit in $\End(H)$ of square-zero matrices 
of rank $t$. 
We proceed to describe the global structure (Theorems 
\ref{thm-lift-of-symmetry-to-elementary-transformations} and 
\ref{thm-correspondence-is-a-linear-combination} below).

\bigskip
A projective polarized K3 surface 
is a simply connected surface $S$ with a trivial 
canonical bundle  and a choice of an ample line bundle
$\StructureSheaf{S}(1)$. There is a sequence
of 19-dimensional irreducible moduli spaces of polarized K3 surface
parametrized by the genus $g\geq 2$ of a curve in the linear system
$\linsys{\StructureSheaf{S}(1)}$. 
The linear system $\linsys{\StructureSheaf{S}(1)}$ is $g$ dimensional and 
the line bundle $\StructureSheaf{S}(1)$
gives rise to a morphism $\varphi : S \rightarrow \PP^g$ which is an
embedding as a surface of degree $2g-2$ if $g > 2$ and a
double cover of $\PP^2$ branched along a sextic if $g=2$. 

Consider first the generic case of a K3 surface with a cyclic Picard group
generated by $\StructureSheaf{S}(1)$. 
A theorem of Mukai implies
that all moduli spaces of Gieseker-Simpson stable 
sheaves on $S$, which happen to be compact, 
are smooth projective hyperkahler 
varieties (they admit an algebraic symplectic structure)
\cite{mukai-symplectic-structure}. 
Compactness of the moduli space is automatic 
for the collection of components 
parametrizing sheaves whose rank $r$, determinant $\StructureSheaf{S}(d)$,
and Euler characteristic $\chi$ satisfy $\mbox{gcd}(r,d,\chi)=1$
(these correspond to primitive vectors in the Mukai lattice introduce below). 
In particular, 
compactness of the stable locus holds for the relative 
compactified Picards over the linear system $\linsys{\StructureSheaf{S}(1)}$
(union of compactified Picards of all curves in the linear system). 
Here, the fact that all curves in the linear system  
$\linsys{\StructureSheaf{S}(1)}$ are reduced and irreducible
is crucial. 
The polarized weight 2 Hodge structure 
of a moduli space of stable sheaves 
is identified as a sub-quotient  of the Mukai-lattice. 
This is the cohomology lattice 
$
\widetilde{H}(S,\Integers):= 
H^0(S,\Integers)\oplus
H^2(S,\Integers)\oplus
H^4(S,\Integers) 
$ 
endowed with the non-degenerate symmetric pairing 

\smallskip
\noindent
\hspace{20ex}
$
\langle(r',c_1',s'),(r'',c_1'',s'')\rangle \ \ := \ \ 
c_1'\cdot c_1'' - r's''-r''s'.
$

\smallskip
\noindent
The Euler characteristic of a coherent sheaf $F$ on $S$ of rank $r$ and
Chern classes $c_1$, $c_2$ is $\chi(F)=2r+\frac{(c_1)^2}{2}-c_2$. 
Following Mukai, 
we associate to $F$ its Mukai vector 
$v(F):=ch(F)\sqrt{Td(S)} =(r,c_1,s)$ where $s=\chi(F)-r$. 
$\Pic(S)$ acts on the Mukai lattice by tensorization 

\smallskip
\noindent
$
\hspace{16ex}
\StructureSheaf{S}(d)\otimes (r,c_1,s) \ = \ 
\left(r,c_1+rd,s+
\frac{1}{2}\StructureSheaf{S}(d)\cdot(2c_1+r\StructureSheaf{S}(d))\right). 
$

\smallskip
\noindent
The dimension of the moduli space  $\M(v)$
of stable sheaves with Mukai vector $v$ is 
$\langle{v},v \rangle+2=(c_1)^2-2rs+2$. Conjecturally, 
the weight 2 Hodge structure of a smooth projective  moduli space  $\M(v)$
is isomorphic to $v^\perp$  if $\dim(\M(v)) > 2$ and to 
$v^\perp/\Integers\cdot v$  if $\dim(\M(v)) = 2$ 
(see \cite{mukai-hodge,ogrady-hodge-str,yoshiyoka} for many known cases). 
Above, $v^\perp$ is the orthogonal hyperplane  with respect to the Mukai 
pairing.

The collection of non-empty moduli spaces consists of the points in the 
``hyperboloid'' 

\smallskip
\hspace{16ex}
$
\V \ := \ \ 
\{v=(r,\StructureSheaf{S}(d),s) \ \mid \ 
\frac{1}{2}\dim \M(v) = 1 + d^2(g-1) -rs \geq 0 
\} 
$

\begin{figure}[h]
\psfig{figure=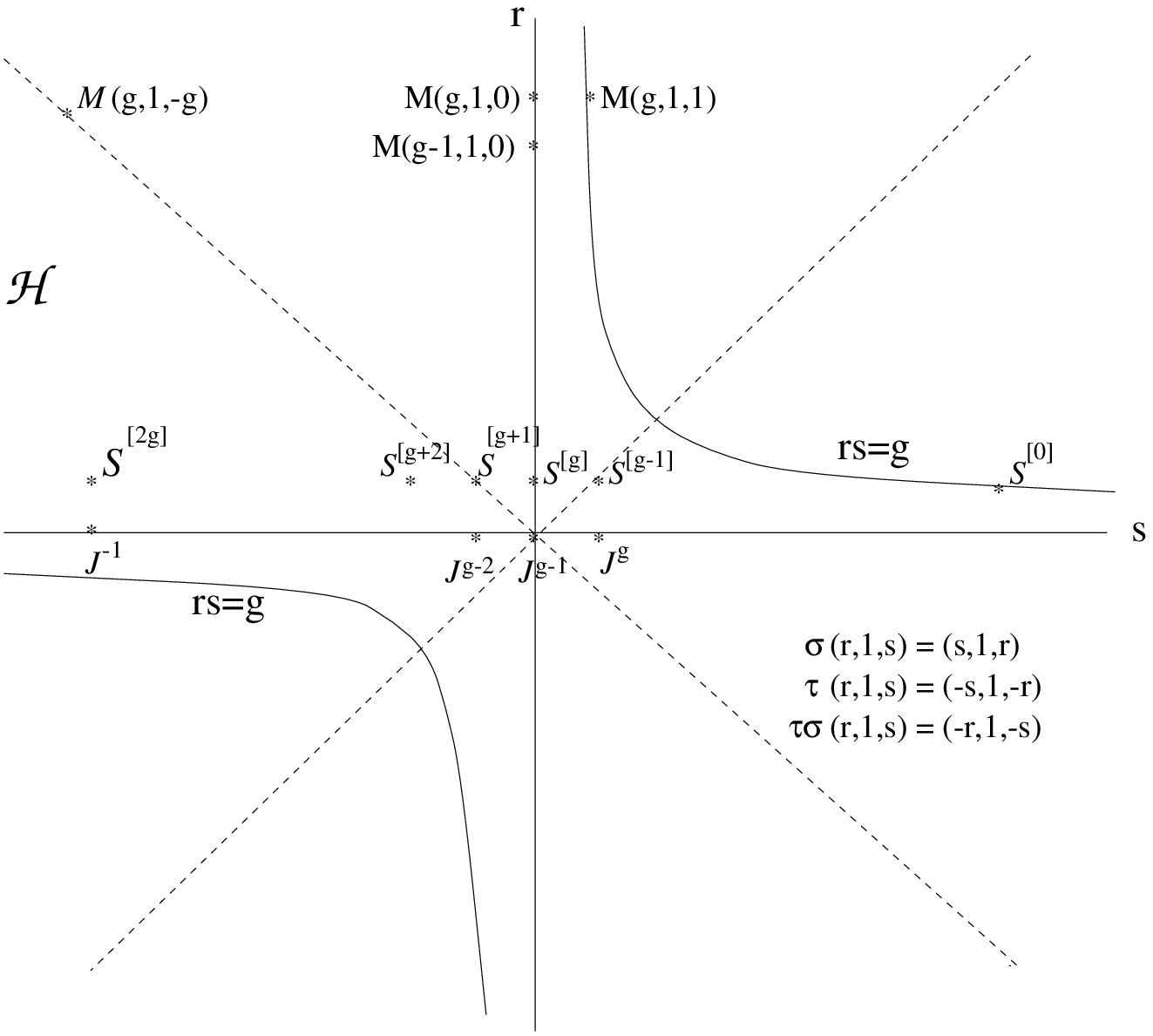}
\caption{\label{eq-graph-of-hyperbola} The region in the Mukai lattice of 
non-empty moduli spaces with  $c_1=\StructureSheaf{S}(1)$.}
\end{figure}

\smallskip
\noindent
in $\Integers^{3}$ satisfying the additional condition that the rank 
$r$ is $\geq 0$. 
$\V$ is symmetric with respect to $r$ and $s$ but the last 
constraint is not. 
Let us illustrate this lack of symmetry in 
the plane $c_1=\StructureSheaf{S}(1)$. 
We get the planar region 
$\Hy \ := \ \V\cap (c_1=1) \ := \ \{v=(r,\StructureSheaf{S}(1),s) \ \mid \ 
\frac{1}{2}\dim \M(v) = g - rs \geq 0
\}$ bounded by a hyperbola 
(see Figure \ref{eq-graph-of-hyperbola}). 
{\em The symmetry is restored via  Brill-Noether theory.}
%
There is a natural $\Integers/2 \times \Integers/2$
action by a group of isometries of the Mukai lattice. 
The hyperbola $\Hy$ is invariant and we get three involutions 
$\sigma$, $\tau$ and $\tau\circ\sigma$ where 
$\sigma$ and $\tau$ are reflections with respect to the 
hyperplanes $r\!-\!s\!=\!0$ and $r\!+\!s\!=\!0$. 
%
%

The Hilbert scheme $S^{[n]}$ of length $n$ zero dimensional 
subschemes is represented by $\M(1,0,1\!-\!n)$ as well as by any 
$\Pic(S)$-translate, hence  
also by $\M(1,1,g\!-\!n)$. The relative compactified Picard of degree 
$i$ over the linear system $\linsys{\StructureSheaf{S}(1)}$ 
is the moduli space 
$J^i := \M(0,1,i\!+\!1\!-\!g)$. Observe that the birational isomorphism 
between $J^g$ and $S^{[g]}$, realized by the Abel-Jacobi map, 
is a lift of the involution $\sigma(1,1,0) = (0,1,1)$. 
The local Torelli theorem holds for the weight 2 Hodge structure 
of projective hyperkakler varieties \cite{beauville-varieties-with-zero-c-1}.
Hence, 
given a universal isometry of the Mukai lattice, it is natural to ask if
it ``lifts'' to a birational transformation on the level of moduli spaces. 
Tyurin observed  that this is indeed the case 
for $\sigma$ and moduli spaces parametrized by vectors in the first
quadrant of the region $\Hy$ 
(see \cite{tyurin-cycles-curves-surfaces} (4.11) and Theorem 4.1). 
In \cite{gottsche-huybrechts} the reflection
$\tau(1,1,-2)=(2,1,-1)$ was
lifted to a birational isomorphism $S^{[g+2]} \leftrightarrow \M(2,1,-1)$.  

\bigskip
\noindent
{\bf Stratified Mukai elementary transformations:}
An {\em elementary Mukai transformation} is a birational transformation
$(M,P)\leftrightarrow (W,P^*)$ between a symplectic variety $M$ 
containing a {\em smooth} codimension $n$ subvariety $P$ which is a 
$\PP^n$-bundle $P\rightarrow Y$. The blow-up of $M$ along $P$ admits a 
second ruling, the contraction of which results in a symplectic variety $W$
containing the dual bundle \cite{mukai-symplectic-structure}. 
In Section \ref{sec-stratified-transformations} Theorem
\ref{thm-stratified-elementary-transformation} we 
construct a {\em stratified analogue of a Mukai elementary transformation} 
- a birational transformation of a symplectic variety $M$ 
admitting a stratification with a highly recursive structure, 
which we call a {\em dualizable stratified collection}
(see (\ref{eq-diagram-of-X-v}) and Definition 
\ref{def-admisible-stratified-collection}). 
The birational transformation produces another  symplectic variety 
$W$ and has the affect of 
``replacing'' each stratum in $M$ - a grassmannian bundle - 
by the dual grassmannian bundle. The base of each grassmannian bundle
is itself the dense open stratum in a smaller dualizable collection. 
The  elementary transformation 
is a duality; when applied twice it recovers the original variety. 
In Section \ref{sec-dual-springer-resolutions} the  simplest 
example is introduced: A Springer resolution of the 
closure of a square-zero nilpotent (co)adjoint orbit in $\LieAlg{gl}_{n}$ 
is related to the dual Springer resolution by a 
stratified elementary transformation. 

We use the Brill-Noether stratification 
(\ref{eq-brill-noether-stratification})
of the moduli spaces in $\Hy$, and prove:

\begin{thm}
\label{thm-lift-of-symmetry-to-elementary-transformations}
Given a Mukai vector $v$ in $\Hy$ with negative rank, define the moduli space 
$\M(v)$ to be identical to $\M(\sigma\circ\tau(v))$. Then, 
the $\Integers/2 \times \Integers/2$ symmetry of $\Hy$, as a set 
of Mukai vectors, 
lifts to an action by  stratified elementary transformations 
on the level of moduli spaces.
Consequently, we obtain a {\bf resolution} of the
birational isomorphisms $\M(v)\leftrightarrow \M(\sigma(v))$ 
and $\M(v)\leftrightarrow \M(\tau(v))$  as a sequence of blow-ups 
along smooth subvarieties, 
followed by a dual sequence of blow-downs. 
\end{thm}

For a more detailed statement, see Theorem 
\ref{thm-mukai-reflection-extends-to-a-stratified-transformation}. 
The formal identification $\M(v)=\M(\sigma\circ\tau(v))$
is well defined because the only vectors in $\Hy$, for which 
both $v$ and $\sigma\circ\tau(v)$ have non-negative rank, are 
Mukai vectors with $r=0$ and $c_1=\StructureSheaf{S}(1)$. 
They correspond to  the compactified Jacobians over 
$\linsys{\StructureSheaf{S}(1)}$ and 
$\sigma\circ\tau$ takes the Mukai vector 
of $J^n$ to that of $J^{2g-2-n}$. These two moduli spaces are naturally
isomorphic (see \cite{le-potier-coherent} Theorem 5.7). 
More conceptually, since $\sigma$ and $\tau$ lift to  stratified 
elementary transformations with respect to the same
(Brill-Noether) stratification, they coincide 
as birational transformations of $\M(\tau(v))$. 
Hence, $\sigma\circ\tau[\M(v)] = \tau\circ\tau [\M(v)] = \M(v)$. 
When the rank of $v$ is negative, $\M(v)$ can be also described as 
a moduli space of (equivalence classes) of complexes of sheaves. 

The general definition of a dualizable stratified collection 
is illustrated in the context of moduli spaces of sheaves. 
Given an integer $t$, denote by $\vec{t}$ the Mukai vector 
$(t,0,t)$ of the trivial rank $t$ bundle. 
Given a Mukai vector $v$ in $\Hy$,
let $\mu(v)$, the {\em distance of $v$ from the boundary of $\Hy$}, be 

\begin{equation}
\label{eq-length-of-bn-stratification}
\mu(v) \ := \ 
\left\{
\begin{array}{c}
\max \{t \ : \ v+\vec{t} \in \Hy, \ t\in \Integers_{\geq 0} \}, 
\ \ \mbox{if} \ \chi(v) \geq 0
\\
\max \{t \ : \ v-\vec{t} \in \Hy, \ t\in \Integers_{\geq 0} \}, 
\ \ \mbox{if} \ \chi(v) \leq 0.
\end{array}
\right.
\end{equation}

\noindent
Then $\mu(v)\!+\!1$ is the length of the Brill-Noether stratification. 
For example, on the line $r\!=\!0$ of compactified relative Jacobians in 
$\Hy$, $\mu(v)\!+\!1$ is equal to the usual 
length of the Brill-Noether stratification of a Petri-generic curve, 
$\mu(J^{g-1+n})=\mu(0,1,n)=\max\{0,\frac{-n+\lceil\sqrt{n^2+4g}\rceil}{2}\}$. 
Choose, for example, $v\in \Hy$ with non-negative 
Euler characteristic 
$\chi(v)=r\!+\!s\geq 0$.
If $\mu(v)>0$ then $\mu(v+\vec{1})=\mu(v)-1$. The square 
$(\mu(v)\!+\!1)\!\times\!(\mu(v)\!+\!1)$ upper triangular matrix
\begin{equation}
\label{eq-diagram-of-M-v}
\begin{array}{cccccc}
\M(v)\supset      & \M(v)^1 \supset       & \cdots             & \supset \M(v)^t  
&   \cdots       & \supset \M(v)^{\mu}
\\
                 & \downarrow           & 
\\
                 & \M(v\!+\!\vec{1}) \supset & \M(v\!+\!\vec{1})^1 \supset & 
&   \cdots       & \supset \M(v\!+\!\vec{1})^{\mu\!-\!1}
\\
                 &                      & \downarrow
\\
                 &                      & \M(v\!+\!\vec{2}) \supset   & 
&   \cdots       & \supset \M(v\!+\!\vec{2})^{\mu\!-\!2}
\\
                 &                      &                        &
&                & \vdots
\\
                 &                      &                        &
&                & \downarrow
\\
                 &                      &                        &
&                & \M(v\!+\!\vec{\mu})
\end{array}
\end{equation}
is a stratified dualizable collection. If $\chi(v) \leq 0$, replace 
$\M(v\!+\!\vec{i})^t$ by $\M(v\!-\!\vec{i})^t$ in the  matrix
(\ref{eq-diagram-of-M-v}) to obtain 
the analogous stratified dualizable collection (here, even if we start with
$v$ satisfying $\rank(v)\geq 0$, the Mukai vector $v\!-\!\vec{i}$ may have 
negative rank and we use the convention of Theorem 
\ref{thm-lift-of-symmetry-to-elementary-transformations}).
The diagonal entries are
symplectic projective moduli spaces. Each row is the
Brill-Noether stratification of  the diagonal entry. 
When $\chi\geq 0$, we set 

\begin{equation}
\label{eq-brill-noether-stratification} 
\M(v)^t \ \ := \ \ \{F \in \M(v) \ \mid \ h^1(F) \geq t  \}
\end{equation}

\smallskip
\noindent
and when $\chi\leq 0$ we use $h^0$ instead. 
Every space $\M(v+\vec{i})^{t-i}$ in the $t$-th column admits a dominant 
rational morphism to the diagonal symplectic entry $\M(v+\vec{t})$ which is 
regular away from the smaller stratum and realizes
\[ 
\M(v\!+\vec{i})^{t\!-\!i}\setminus \M(v\!+\!\vec{i})^{t\!-\!i\!+\!1} 
\ \  \longrightarrow \ \ 
\M(v\!+\!\vec{t})\setminus \M(v\!+\!\vec{t})^1
\] 
as a Grassmannian bundle. The
fiber over $E\in \M(v+\vec{t})\setminus \M(v+\vec{t})^1$ is
$G(t\!-\!i,H^0(E))$. As a subvariety of $\M(v)$, 
the codimension of each Grassmannian bundle is equal to the dimension
of the Grassmannian fiber. Hence, the  projectivized normal bundle of each 
Grassmannian bundle is isomorphic to its relative 
projectivized cotangent bundle. The latter is a
bundle of homomorphisms and admits a canonical determinantal stratification.
The stratified elementary transformation is performed by blowing up the matrix
(\ref{eq-diagram-of-M-v}) column by column from right to left. 
The key to the success of the recursion is the fact that the intersection
of the strict transform 
of the Brill-Noether stratification with a  new exceptional 
divisor (which is fibered by projectivized cotangent bundles of grassmannians)
is precisely the determinantal stratification of the latter. 
The existence of two dual sequences of blow-downs on the top iterated 
blow-up can be seen already on the 
level of the projectivised cotangent bundle of a single grassmannian: 
Given a vector space $H$, the top iterated blow-up
of $\PP{T}^*G(t,H)$ is naturally isomorphic to that of $\PP{T}^*G(t,H^*)$.

The top iterated blow-up $B^{[1]}\PP{T}^*G(t,H)$ is a particularly nice
compactification of the open dense $GL(H)$-orbit in $\PP{T}^*G(t,H)$. 
It admits an interpretation as a moduli space of complete collineations 
and its cohomology ring is well understood \cite{bifet-deconcini-processi}.
Complete collineations and complete quadrics play an important role in 
enumerative geometry \cite{laksov-review}. They received a modern treatment in 
\cite{vainsencher,laksov,kleiman-thorup,thadeus-colineations}. 

The top iterated blow-up $B^{[1]}\M_S(v)$ is a coarse moduli space of 
{\em complete stable sheaves}. A complete sheaf 
$(F,\rho_2,\rho_3,\dots,\rho_k)$ consists of a sheaf $F$ on $S$ and a 
sequence of non-zero homomorphisms, up to a scalar factor, 
defined recursively by
\begin{eqnarray*}
\rho_2 & : &  H^0(F) \rightarrow H^1(F) \ \ \ \ \ \mbox{and} 
\\
\rho_{i+1} & : & \ker(\rho_i) \rightarrow 
\coker(\rho_i), \ \ 2 \leq i \leq k-1,
\end{eqnarray*}
such that the last homomorphism $\rho_k$ is surjective if $\chi(F)\geq 0$
and injective if $\chi(F)\leq 0$. 
The sequence is empty if either $H^0(F)$ or  $H^1(F)$ vanishes. The sequence 
$\{\rho_2,\rho_3,\dots, \rho_k\}$ should be viewed as a truncation of
a complete collineation $\{\rho_1,\rho_2,\dots, \rho_k\}$. The latter 
behaves well when the sheaf $F$ varies in a flat family. 
A choice of a section  $\gamma$ of a sufficiently ample line bundle on $S$ 
gives rise to a homomorphism 
$\rho_1: V_0(F) \rightarrow V_1(F)$ between vector spaces
which ranks depend only on the Mukai vector $v(F)$ 
(see (\ref{eq-exact-seq-is-locally-free-presentation-of-R1})). 
The kernel of $\rho_1$ is $H^0(F)$ and cokernel is $H^1(F)$. 
The notion of families of complete collineations is subtle but
well understood \cite{laksov,kleiman-thorup}. 
It leads to the notion of families of complete sheaves
once we fix a section $\gamma$ as above. Thaddeus' work 
on complete collineations \cite{thadeus-colineations} combined  
with Theorem \ref{thm-lift-of-symmetry-to-elementary-transformations} 
suggests that the stratified elementary transformations 
in Theorem \ref{thm-lift-of-symmetry-to-elementary-transformations}
come from a variation of Geometric Invariant Theory quotients in the sense of 
\cite{dolgachev-hu,thadeus-git}. 

The moduli space of complete sheaves $B^{[1]}\M_S(v)$ is instrumental 
in studying the intersection theory on moduli spaces. Indeed, it plays 
a central role in the proof of Theorem 
\ref{thm-correspondence-is-a-linear-combination}. This is not surprising, 
considering the important role played by complete collineations in 
classical enumerative geometry. $B^{[1]}\M_S(v)$ is different in general 
from the closure of the graph of the birational isomorphism
in $\M(v)\times \M(\sigma(v))$. The latter is the moduli space 
of coherent systems $G^0(\chi(v),\M(\sigma(v)))$ discussed in section 
\ref{coherent-systems}. The two are equal only in the case of a
Mukai elementary transformation, i.e., when $\mu(v) \leq 1$. While both 
moduli spaces are useful in the proof of Theorem 
\ref{thm-lift-of-symmetry-to-elementary-transformations}, it is the moduli 
space of complete sheaves which admits the full recursive structure 
which makes transparent the analogy with dual Springer resolutions of 
nilpotent orbits. 


\bigskip
\noindent
{\bf A correspondence inducing isomorphism of cohomology rings:}

Birational hyperkahler varieties $M'$, $M''$  are especially similar
(under mild conditions which are satisfied by our stratified
elementary transformations):
{\em There exists a family $\M  \rightarrow D^\times$
of hyperkahler varieties over the punctured disk
with two extensions $\M'$, $\M''$ to smooth families over the disk with 
special fibers $M'$, $M''$ respectively}
\cite{huybrects-iseparable-points-in-moduli}. 
Consequently, the Hodge structure,
the cohomology {\em ring} structure and all continuous invariants
of $M'$ and $M''$ coincide. 
The Hodge conjecture suggests the existence of a correspondence 
(an algebraic cycle in $M' \times M''$) which induces 
the isomorphism of Hodge structures. 
This correspondence contains the closure of the graph of the
birational isomorphism as a component. It is easy to see that there are other
components. For example, we saw that $J^g$ and $S^{[g]}$ are birational
$\sigma(J^g)=S^{[g]}$ (see Figure \ref{eq-graph-of-hyperbola}). 
The $g$-th symmetric product $C^{[g]}$ of a smooth member (curve) $C$
of $\linsys{\StructureSheaf{S}(1)}$, is a subvariety 
in $S^{[g]}$. The birational image of $C^{[g]}$ in $J^g$ is the
Jacobian $J^g_C$ which has self intersection $0$. However, 
the self intersection of $C^{[g]}$ in $S^{[g]}$ is 
{\scriptsize
${\displaystyle
\left(
\begin{array}{c}
2g\!-\!2
\\
g
\end{array}
\right).
}
$
}
Recently, progress has been made in understanding the ring structure
of the cohomology of the Hilbert schemes 
\cite{ellingsrud-stromme-ring-str,lehn}. 
In order to translate this to an understanding of the cup-product 
in the cohomology of birational components such as
$\sigma(S^{[n]})$, one needs to compute explicitly the isomorphism 
$H^*(S^{[n]},\Integers)\IsomRightArrow H^*(\sigma(S^{[n]}),\Integers)$. 

\medskip
\begin{thm}
\label{thm-correspondence-is-a-linear-combination}
Let 
$
\{M(i)^j\}  \longleftrightarrow \{W(i)^j\}
$ 
be a stratified elementary transformation 
between two stratified dualizable collections associated with 
irreducible projective symplectic varieties $M=M(0)^0$ and $W=W(0)^0$. Then,  
the natural isomorphism $H^*(M,\Integers) \cong H^*(W,\Integers)$
is induced by a correspondence
which, as a cycle $\Delta$ in $M\times W$, is a sum
\begin{equation}
\label{eq-linear-combination}
\Delta \ = \Delta_0 + \sum_{t=1}^\mu \Delta_t
\end{equation}
where $\Delta_0$ is the closure of the graph of the birational transformation 
$M \leftrightarrow W$ and $\Delta_t$ is the closure of the fiber product 
$M^t\times_{M(t)}W^t$ of the two grassmannian bundles. 
\end{thm}

\noindent
Note that the dimension of $\Delta_t$ is equal to $\dim(M)$. 
The Theorem is proven in Section 
\ref{proof-that-correspondence-induces-isomorphism-of-coho-rings}. 
A few applications of Theorem 
\ref{thm-correspondence-is-a-linear-combination}
are discussed in Section \ref{subsec-self-dual-moduli-spaces}.

\bigskip
\noindent
{\bf Auto-equivalences of the derived category:}

Let $S$ be a K3 surface and 
consider the group $G$ of Hodge isometries of the Mukai lattice 
$\widetilde{H}(S,\Integers)$ of $S$. 
These are integral isometries which send $H^{2,0}(S)$ to itself when 
the lattice is complexified. Any Hodge-isometry $\phi$ 
leaves the algebraic sublattice $\widetilde{H}_{Alg}(S,\Integers)$ invariant.
If $\phi$  restricts to the identity automorphism of 
$\widetilde{H}_{Alg}(S,\Integers)$ then it is induced by an automorphism of $S$
(the Torelli theorem).
It follows that $G$ fits in the exact sequence:
\begin{equation}
\label{eq-short-exact-seq-Hodge-isometries}
0 \rightarrow 
\Aut(S,\Pic(S)) \rightarrow 
G  \rightarrow 
\Aut(\widetilde{H}_{Alg}(S,\Integers)) \rightarrow 0.
\end{equation}
Let $G_{Der}$ be the group of (covariant and contravariant) 
auto-equivalences of the bounded derived category of $S$. 
We have a map  
\begin{eqnarray*}
\eta : D^b(S) & \rightarrow & \widetilde{H}(S,\Integers)
\\
(F_i,\partial_i) & \mapsto & \sum_i (-1)^i v(F_i)
\end{eqnarray*}
sending an object represented by a complex $(F_i,\partial_i)$ to 
the Mukai vector of the associated class in K-theory, i.e., to 
the alternating sum of the Mukai vectors of its coherent sheaves. 
Notice that $\eta$ is equivariant with respect to 
the cyclic subgroup $\langle T \rangle \subset G_{Der}$ generated by the 
translation auto-equivalence  $T : D^b(S) \rightarrow D^b(S)$ once we send
$T$ to $-Id \in G$. 
Orlov proved that there is a surjective homomorphism
\begin{equation}
\label{eq-surjective-homo-G-Der-onto-G}
G_{Der}/\langle T^2\rangle \ \ \rightarrow \ \ G.
\end{equation}
Any auto-equivalence $\Phi$ of the derived category $D^b(S)$
induces a Hodge isometry $\phi$ of $\widetilde{H}(S,\Integers)$
and $\eta$ is $(\Phi,\phi)$-equivariant (\cite{orlov}, Theorem 2.2). 
Moreover, any Hodge isometry can be lifted to an 
auto-equivalence of the derived category $D^b(S)$
(\cite{orlov}, proof of Theorem 3.13). 
For example, 
\begin{enumerate}
\item
Tensorization by a line bundle $\gamma$ on $S$ is induced by the Fourier-Mukai 
functor 
$R^\cdot\pi_{2_*}\E\stackrel{L}{\otimes}\pi^*_1(\cdot) \ : \ D^b(S)
\rightarrow D^b(S)$ associated to the sheaf $\E$ on $S\times S$ where $\E$
is supported on the diagonal as the line bundle $\gamma$. 
\item
The Hodge-isometry $-\tau$ ($\tau$ as in Theorem 
\ref{thm-lift-of-symmetry-to-elementary-transformations}) 
is induced by the Fourier-Mukai functor
associated to the ideal sheaf of the diagonal.
\item
The Hodge isometry $-\sigma\circ\tau$ is induced by the {\em contravariant} 
involutive functor 
$R^{\cdot}\SheafHom_S(\cdot,\StructureSheaf{S}) : 
D^b(S)\rightarrow D^b(S)^{op}$. 
\end{enumerate}

Assume that the K3 surface $S$ has a cyclic Picard group and 
let us introduce yet a third group $G_{bir}$. Let 
\[
\M_S(\bullet) := \cup \{ \M_S(v) : \ v \in \widetilde{H}_{Alg}(S,\Integers)
\ \ \mbox{is a primitive Mukai vector}, \ \ \langle v,v\rangle \geq -2 \}
\]
be the disjoint union of moduli spaces of stable sheaves on $S$
with primitive algebraic Mukai vectors of non-negative dimension. 
Note that these moduli spaces  are all smooth and compact. 
We use the convention that, if $\rank(v) \leq 0$, then
$\M_S(v)$ is equal to $\M_S(-v)$. The group $G_{bir}$ is defined to be 
the group of all possible lifts of elements of $G$ to birational 
automorphisms of $\M_S(\bullet)$. 
A natural question arises:

\smallskip
\noindent
{\bf Question:} 
{\em 
Does the homomorphism 
(\ref{eq-surjective-homo-G-Der-onto-G}) factor through $G_{bir}$? 
}

\smallskip
The Fourier-Mukai functor lifting an isometry to an auto-equivalence
sends, in general, a stable sheaf to the class of a complex supported at
more than one degree. 
Nevertheless, 
there seems to exist, at least, a surjective 
homomorphism $G_{bir} \rightarrow G$.
Clearly, $\Pic_S$ and $\Aut(S)$ 
lift to  $G_{bir}$. The element $-Id$ lifts, by definition.
Theorem \ref{thm-lift-of-symmetry-to-elementary-transformations}
suggests that the Hodge isometries $\sigma$ and $\tau$ 
lift to $G_{bir}$. In the Theorem the lift is carried out only for the 
collection with $c_1=\StructureSheaf{S}(1)$, but our definitions of
$\sigma$ and $\tau$ on the Zariski open Brill-Neother stratum 
seem to extend to birational isomorphisms for other values of $c_1$.  
If indeed $\sigma$ and $\tau$ lift, then the image of 
$G_{bir} \rightarrow G$ has at most a finite index (and is 
surjective if $g=2$). This follows from the exactness of 
(\ref{eq-short-exact-seq-Hodge-isometries}) and the fact that 
$\{-Id,\sigma,\tau\}$ and $\Pic_S$ generate a finite index 
subgroup $\Gamma$ of the group of isometries of the rank $3$ lattice 
$\widetilde{H}_{Alg}(S,\Integers)$. If $g=2$, $\Gamma$ is the whole group,
but in general it is a proper subgroup. For example, if $g=7$ the isometry
{\footnotesize
\[
\left(\begin{array}{c}
r\\d\\s
\end{array}\right)
\ \ \ \mapsto \ \ \ 
\left(\begin{array}{ccc}
2 & 12 & 3\\
1 & 5 & 1\\
3 & 12 & 2
\end{array}\right)
\left(\begin{array}{c}
r\\d\\s
\end{array}\right)
\]
}
is not in the subgroup  $\Gamma$
because any isometry in $\Gamma$ takes $(0,0,1)$ to a vector $(r,d,s)$
such that $(g-1)$ divides exactly one of the two integers $\{r,s\}$
and is relatively prime to the other.

\smallskip
A much harder problem would be to 
resolve the birational isomorphisms in $G_{bir}$. It is natural to
try to generalize our results to other reflections of the Mukai lattice. 
The group of isometries of the rank $3$ lattice 
$\widetilde{H}_{Alg}(S,\Integers)$ 
contains a finite index subgroup $W \subset G$ generated by reflections
with respect to Mukai vectors with $\langle v, v\rangle = \pm 2$. 
In fact, setting 
$\sigma' := \StructureSheaf{}(1)\circ \sigma\circ \StructureSheaf{}(-1)$ 
and $\tau' := \StructureSheaf{}(1)\circ \tau\circ \StructureSheaf{}(-1)$
we have the equality 
\[
\StructureSheaf{}(2) \ = \ \sigma'\tau'\sigma\tau \in W.
\]
Our results suggest a relationship between dual pairs of hyperkahler 
resolutions of singularities and reflections in $G_{bir}$. 
While at the level $c_1 = \StructureSheaf{}(1)$ the reflections 
$\sigma$ and $\tau$ correspond to Springer resolutions of a nilpotent
orbit with a simple (well ordered) stratification, it seems that
for other values of $c_1$ more complicated singularities will arise.
It would be interesting, for example, to interpret the resolution
$\M_S(1,0,1-n) \ \cong \ S^{[n]}\rightarrow \Sym^n S$ of the symmetric 
product as a special case of a lift of some reflection ($\tau$ ?) to 
$G_{bir}$. 

\medskip
The rest of the paper is organized as follows.
Sections \ref{sec-stratified-transformations}, 
\ref{sec-determinantal-varieties},
and \ref{sec-construction-of-the-dual-collection}
are devoted to the general study of stratified elementary transformations.
In section \ref{sec-stratified-transformations}
we construct the stratified elementary 
transformations. Section \ref{sec-determinantal-varieties} contains the 
background information on determinantal
varieties needed for the proof of the construction. In section 
\ref{sec-construction-of-the-dual-collection} we complete the proof of the 
construction. We also prove Theorem 
\ref{thm-correspondence-is-a-linear-combination} identifying the
correspondence inducing the cohomology ring isomorphism.
Section \ref{sec-brill-noether-for-k3} is devoted to our main example, 
the moduli spaces of sheaves on a K3 surface. The organization of section
\ref{sec-brill-noether-for-k3} is described at the beginning of that section.

\medskip
{\em Acknowledgments:} It is a pleasure to acknowledge fruitful
conversations with D. Cox, R. Donagi, 
B. Fantechi, L. G\"{o}ttsche, V. Ginzburg, 
D. Huybrechts, S. Kleiman, I. Mirkovic, T. Pantev and M. Thaddeus.

\section{Stratified Elementary Transformations} 
\label{sec-stratified-transformations}

In section \ref{sec-dual-springer-resolutions} we introduce the prototypical 
example of the symplectic birational isomorphisms considered in this paper. 
In section \ref{subsec-constraction-of-stratified-transformations} we 
introduce a global analogue for projective symplectic varieties. 

\subsection{Dual resolutions of the closure of a nilpotent orbit}
\label{sec-dual-springer-resolutions}

Let $H$ be a vector space of dimension $h$, $t\leq h/2$ an integer and 
consider the natural morphism $\pi_1 : T^*G(t,H) \rightarrow \End(H)$ onto
the closure $\overline{{\cal N}^t}$ in $\End(H)$ of the nilpotent orbit 
${\cal N}^t$ of square-zero nilpotent elements of rank $t$. The natural 
isomorphism $\End(H) \cong \End(H^*)$ provides another resolution 
\begin{equation}
\label{eq-flop}
T^*G(t,H) \ \LongRightArrowOf{\pi_1} \ \overline{{\cal N}^t} \ 
\LongLeftArrowOf{\pi_2} \ T^*G(t,H^*). 
\end{equation}

There is a simultaneous deformation of (\ref{eq-flop}) over $\ComplexNumbers$ 
which smoothes $\overline{{\cal N}^t}$ away from the special fiber and 
deforms $\pi_i$ to isomorphisms. 
The cotangent bundle of  $G(t,H)$ 
is isomorphic, as a bundle of Lie subalgebras  
of $\SheafEnd(H_{G(t,H)})$, to the nilpotent radical of the corresponding parabolic
subalgebra ${\cal P}_{G(t,H)} \subset \SheafEnd(H_{G(t,H)})$. 
We have  a unique non-trivial extension

\begin{equation}
\label{eq-idempotent-extension-of-cotangent-bundle}
0 \rightarrow T^*_{G(t,H)} \rightarrow E(H) 
\rightarrow \StructureSheaf{G(t,H)} \rightarrow 0. 
\end{equation}

\noindent
A non-zero section $\gamma$ of $\StructureSheaf{G(t,H)}$ determines a 
symplectic $T^*_{G(t,H)}$ torsor $X_\gamma$ embedded in $E(H)$. 
$X_\gamma$ does not have any section. Hence the
lagrangian zero-section of $T^*_{G(t,H)}$ does not deform. 
The vector bundle $E(H)$ admits an embedding as a bundle of subalgebras
in $\End(H_{G(t,H)})$. Each point in $G(t,H)$ determines a decomposition of 
the Levi factor of the parabolic subalgebra (with ranks $t$ and $h\!-\!t$). 
$E(H)$ is the subalgebra of ${\cal P}_{G(t,H)}$ which projects to the 
center in the first factor and to zero in the second factor. 
In other words, $E(H)$ is the subalgebra of endomorphisms leaving the subbundle
$\tau(H)$ invariant, whose image in $\SheafEnd(\tau(H))$ is 
a scalar multiple of the identity and whose image in $\SheafEnd(q(H))$ 
is zero. The morphism $E(H)\rightarrow \End(H)$ restricts to an embedding 
of $X_\gamma$, $\gamma\neq 0$, as a smooth closed orbit. 

\bigskip
The birational isomorphism (\ref{eq-flop}) is a special case of the setup 
of Section \ref{subsec-constraction-of-stratified-transformations}.
Consider the square 
$(t\!+\!1)\times(t\!+\!1)$ upper triangular matrix whose rows are the 
determinantal stratification of $T^*G(k,H)$, $0\leq k \leq t$, 
indexed as follows. 
$T^*G(k,H)$ is isomorphic to the homomorphism bundle 
$\SheafHom(q_k,\tau_k)$ from the tautological quotient bundle to the 
tautological sub-bundle. The generic point corresponds to a surjective 
homomorphism. $T^*G(k,H)^i$ denotes the locus with $i$-dimensional cokernel.

\begin{equation}
\label{eq-dualizable-diagram-for-a-cotangent-bundle-of-grass}
\begin{array}{c}
\left(
\begin{array}{cccc}
T^*G(t,H) \supset & T^*G(t,H)^1           &     \cdots              
& \supset T^*G(t,H)^{t}
\\
                   &  \downarrow          &                   
\\
                   & T^*G(t\!-\!1,H) \supset & T^*G(t\!-\!1,H)^1   
& \cdots \ T^*G(t\!-\!1,H)^{t\!-\!1}
\\
                   &                        & \downarrow 
\\
                   &                        &  T^*G(t\!-\!2,H) \supset  
& \cdots \ T^*G(t\!-\!2,H)^{t\!-\!2}
\\
                   &                        &     
& \vdots
\\
                   &                        &     
& \downarrow 
\\
                   &                        &     
& T^*G(0,H)
\end{array}
\right)
\\
\begin{array}{cccc}
\hspace{7ex}\downarrow \hspace{7ex} & \hspace{5ex}\downarrow \hspace{11ex}&  
\hspace{5ex}\dots \hspace{5ex}& \hspace{8ex} \downarrow \hspace{7ex}
\end{array}
\\
\left(
\begin{array}{cccc}
\hspace{8ex} \overline{{\cal N}^t} \ \ \supset  \hspace{4ex} & 
\hspace{2ex}\overline{{\cal N}^{t\!-\!1}} \ \ \supset \hspace{6ex}
& \hspace{5ex} \dots \hspace{5ex} & \hspace{5ex}  \{0\} \hspace{7ex}
\end{array}
\right)
\end{array}
\end{equation}

\noindent
Above, $T^*G(0,H)$ is a point. 
For the sake of compatibility with the  notation in the 
following sections, set $n:=h\!- \!2t\!+1$ and $\mu=t$. We have a natural 
$G(r,n\!+\!2r\!-\!1)$-bundle: 
\begin{equation}
\label{eq-stratum-of-cotangent-bundle-is-a-grassmannian-bundle}
f_{t,r} \ : \ [T^*G(t,H)^{r} \setminus T^*G(t,H)^{r\!+\!1}]
\ \rightarrow \ 
[T^*G(t-r,H)\setminus T^*G(t-r,H)^1], \ \ \    0\leq r \leq \mu. 
\end{equation}

\noindent
$[T^*G(t,H)^{r} \setminus T^*G(t,H)^{r\!+\!1}]$ is isomorphic to the
Zariski open stratum in the homomorphism bundle 
$\SheafHom(q_{h\!-\!t\!+\!r},\tau_{t\!-\!r})$ over 
$Flag(t\!-\!r,t,h\!-\!t\!+\!r,H)$. 
$[T^*G(t-r,H)\setminus T^*G(t-r,H)^1]$ 
is isomorphic to the
Zariski open stratum in the homomorphism bundle 
$\SheafHom(q_{h\!-\!t\!+\!r},\tau_{t\!-\!r})$ over 
$Flag(t\!-\!r,h\!-\!t\!+\!r,H)$. 
$Flag(t\!-\!r,t,h\!-\!t\!+\!r,H)$
is isomorphic to the bundle $G(r,\tau_{h\!-\!t\!+\!r}/\tau_{t\!-\!r})$
over $Flag(t\!-\!r,h\!-\!t\!+\!r,H)$. 
The morphism (\ref{eq-stratum-of-cotangent-bundle-is-a-grassmannian-bundle})
is the pullback of the morphism 
\[
Flag(t\!-\!r,t,h\!-\!t\!+\!r,H) \ \rightarrow 
Flag(t\!-\!r,h\!-\!t\!+\!r,H)
\]
and is hence a Grassmannian fibration. In a similar fashion, every 
entry in each column admits a rational morphism to the symplectic 
diagonal entry in that column. 
Observe that the morphism to $\overline{{\cal N}^t}$ 
contracts these grassmannian fibrations. 

\bigskip
In this paper we show that the birational isomorphisms 
(\ref{eq-flop}) describe the local structure of the birational isomorphisms
of moduli spaces of sheaves in Theorem 
\ref{thm-mukai-reflection-extends-to-a-stratified-transformation}. 
This fact depends heavily on the condition \ref{cond-linear-system} 
imposed on the first Chern class. 
It is plausible that the structure of birational 
isomorphisms of moduli spaces of sheaves with a more general first Chern class
would also admit a group theoretic model. There is a generalization of the 
above construction valid for a nilpotent orbit ${\cal N}_\eta$ in $\End(H)$
associated to any partition $\eta:=(p_1\geq p_2 \geq \cdots\geq p_k)$ of $h$. 
In the general case there are many Springer resolutions 
of the closure $\overline{{\cal N}}_\eta$ but they come in dual pairs. 
Given a partition $\eta$ we get a Young diagram whose $i$-th row consists of 
$p_i$ boxes. The dual partition 
$\hat{\eta}=(\hat{p}_1\geq \hat{p}_2 \geq \cdots\geq \hat{p}_m)$
is defined setting $\hat{p}_i$ as the length of the $i$-th column of 
the diagram of $\eta$. Choose a permutation $\theta\in \mbox{Sym}_m$
and set 
$n(\theta) := (n_1\leq n_2 \leq \cdots\leq n_{m\!-\!1})$ where 
\[
n_j :=  \sum_{i=1}^j \hat{p}_{\theta(i)}. 
\]
Then the cotangent bundle $T^*Flag(n_1,\dots,n_{m\!-\!1},H)$ is a resolution of
the closure $\overline{{\cal N}}_\eta$ where ${{\cal N}}_\eta$ is the 
nilpotent orbit of 
matrices whose Jordan canonical form has $k$ blocks and the $i$-th block
is a regular nilpotent $p_i\times p_i$ matrix. The analogue of 
(\ref{eq-flop}) is then
\begin{equation}
\label{eq-flop-for-flags}
T^*Flag(n_1,\dots,n_{m\!-\!1},H) \ \ \rightarrow \ 
\overline{{\cal N}}_\eta \ \leftarrow \ \ 
T^*Flag(h\!-\!n_{m\!-\!1},\dots,h\!-\!n_1,H). 
\end{equation}
The recursive structure of a Flag variety as a Grassmannian bundle over a 
smaller flag variety suggests an inductive way to reduce the study of 
(\ref{eq-flop-for-flags}) to that of (\ref{eq-flop}). 
It seems likely that the work of Thaddeus \cite{thadeus-colineations}
could be generalized to exhibit the collection of birational 
varieties $\PP{T}^*Flag(n(\theta),H)$, $\theta\in \mbox{Sym}_m$, as different 
geometric invariant theory quotients associated to polarizations in 
different faces of a convex polytope. The torus acting is likely to be the 
center of the Levi factor of the parabolic subgroups in $SL(H)$. 

\subsection{Construction  of stratified transformations}
\label{subsec-constraction-of-stratified-transformations}

Determinantal stratifications of symplectic varieties tend to have a 
recursive structure. Roughly, the recursive structure arises by taking the 
quotient of each stratum by the null-foliation of the symplectic form. 
In nice situations, these foliations are Grassmannian bundles. 
Dualizing them gives rise to a birational symplectic variety. 
We formalize the recursive setup is this section. 

Let $\mu$ be a non-negative integer, 
$X(r)$, $0 \leq r \leq \mu$ a collection of smooth projective 
symplectic varieties and 
\[
X(r) = X(r)^0 \supset X(r)^1 \supset \dots \supset X(r)^{\mu-r} \supset
X(r)^{\mu+1-r}= \emptyset
\]
a flag of closed subschemes. We denote $X(0)$ also by $X$ and assume that 
$X$ is connected.
It is convenient to arrange the data in an upper triangular 
$(\mu\!+\!1)\times(\mu\!+\!1)$-matrix with symplectic diagonal entries:

\begin{equation}
\label{eq-diagram-of-X-v}
\begin{array}{cccccc}
X(0)\supset      & X(0)^1 \supset       & \cdots             & \supset X(0)^t  
&   \cdots       & \supset X(0)^{\mu}
\\
                 & \downarrow           & 
\\
                 & X({1}) \supset & X({1})^1 \supset & 
&   \cdots       & \supset X({1})^{\mu-1}
\\
                 &                      & \downarrow
\\
                 &                      & X({2}) \supset   & 
&   \cdots       & \supset X({2})^{\mu-2}
\\
                 &                      &                        &
&                & \vdots
\\
                 &                      &                        &
&                & \downarrow
\\
                 &                      &                        &
&                & X({\mu})
\end{array}
\end{equation}

\noindent
We will assume below (in Condition \ref{cond-abstract-tyurin-morphism}) that
every entry $X(i)^t$ admits a rational morphism 
to the symplectic diagonal entry $X(i\!+\!t)$ in the same column. This morphism
is regular away from $X(i)^{t\!+\!1}$ and realizes
\[
X(i)^t \setminus X(i)^{t\!+\!1} \ \ \longrightarrow 
X(i\!+\!t) \setminus X(i\!+\!t)^1
\]
as a Grassmannian bundle. 

We will sometimes denote the general $X(r)$ by 
$M$ and we set $\mu(M):=\mu-r$, $M^t:=X(r)^t$ and $M(r'):=X(r+r')$.
Denote by $B^{k}M^t$, $t < k \leq \mu(M)$ the blow-up of $M^t$
along the subscheme $M^k$. We introduce the following notation for iterated
blow-ups:
\begin{enumerate}
\item[$k=\mu$]
Set $B^{[\mu(M)]}M^t := B^{\mu(M)}M^t$, \ \ $t<\mu(M)$.
\item[$t < k$]
Define recursively 
$B^{[k]}M^t$, $t < k \leq \mu(M)$ to be the blow-up of $B^{[k+1]}M^t$ 
along $B^{[k+1]}M^k$. Note that we used here recursively 
the base change property of blowing-up to conclude that 
$B^{[k+1]}M^k$ is isomorphic to the strict transform
of $M^k$ in $B^{[k+1]}M^t$ for all $t\leq k$.
\item[$t=k$]
We further denote by $B^{[k]}M^k$ the exceptional divisor in 
$B^{[k]}M$ corresponding to $M^k$.
\item[$t > k$]
Let $B^{[k]}M^t$, $t > k$ be the strict 
transform in $B^{[k]}M$ of the Cartier divisor $B^{[k+1]}M^{t}$.
\end{enumerate}
If $\mu \geq 1$, denote by $n(M)$ the codimension of $M^1$ in $M$.

\begin{condition} \label{cond-codimension}
$\dim(M(r)) \ = \ \dim(M) - 2r[n(M)+r-1]$. 
\end{condition}

\begin{condition} \label{cond-abstract-tyurin-morphism}
(Quotient by the null foliation of the symplectic form) 
There exists a morphism 
\begin{equation} \label{eq-abstract-tyurin-morphism}
f_{r,t}: B^{[t+1]}X(r)^t \rightarrow X(r+t).
\end{equation}

\noindent
When using the notation  $M=X(r)$ we denote $f_{r,t}$ by $f_t$.
The morphism $f_t$ lifts to a smooth projective morphism

\begin{equation} \label{eq-lifted-abstract-tyurin-morphism}
\tilde{f}_{t}: B^{[t+1]}M^t \rightarrow B^{[1]}M(t)
\end{equation}
realizing $B^{[t+1]}M^t$ as a 
$G(\PP^{t-1},\PP^{n(M)+2t-2})$-bundle. 
\end{condition}

\begin{new-lemma}
\label{lemma-smoothness-of-iterated-blow-up}
$B^{[t+1]}X(r)^t$ is smooth. 
In particular, $X(r)^t \setminus X(r)^{t+1}$ is smooth.
\end{new-lemma}

\noindent
Note that the Lemma implies that each of the iterated blow-ups is
a sequence of blow-ups  of smooth varieties along smooth subvarieties.

\medskip
\noindent
{\bf Proof:} (of lemma \ref{lemma-smoothness-of-iterated-blow-up})
The proof is by descending induction on $r+t$.
If $r+t=\mu$ then $B^{\mu-r+1}X(r)^{\mu-r}=X(r)^{\mu-r}$ and Condition
\ref{cond-abstract-tyurin-morphism}
implies that $B^{[\mu-r+1]}X(r)^{\mu-r}$ is a 
Grassmannian-bundle 
over $B^{[1]}X(\mu)=X(\mu)$. Since $X(\mu)$ is smooth,
so is $B^{\mu-r+1}X(r)$.

\noindent
Induction step:
Condition \ref{cond-abstract-tyurin-morphism}
implies that 
$B^{[t+1]}X(r)^t$ is a Grassmannian bundle over 
$B^{[1]}X(r+t)$. It suffices to prove that $B^{[1]}X(r+t)$
is smooth. $B^{[1]}X(r+t)$ is the blow-up of $B^{[2]}X(r+t)$ along
$B^{[2]}X(r+t)^1$. By the induction hypothesis, $B^{[2]}X(r+t)^1$
is smooth. Hence,
so is $B^{[1]}X(r+t)$. 
\EndProof

\medskip
Conditions \ref{cond-codimension} 
and \ref{cond-abstract-tyurin-morphism} imply the equality 
\begin{equation}
\label{eq-n}
n(X(r)) = n + 2r. 
\end{equation}

\noindent
A short calculation shows that we have an equality 
\begin{equation}
\label{eq-codimension-is-sum-of-codimensions}
codim(M^k,M) = codim(M^t,M) +
codim(M(t)^{k-t},M(t)), \ \ \mbox{for} \ \
t \leq k \leq \mu(M). 
\end{equation}

Observe that for $t=1$ we get that 
$\tilde{f}_1 : B^{[2]}M^1 \rightarrow B^{[1]}M(1)$
is a $\PP^{n(M)}$-bundle. Similarly, 
\[
\tilde{f}_{r,1}:B^{[2]}X(r)^1 \rightarrow B^{[1]}X(r+1) 
\]
is a $\PP^{n+2r}$-bundle. For $1 \leq r \leq \mu$, we denote
the $\PP^{n+2r-2}$-bundle over $B^{[1]}X(r)$ by $\PP{W}_{B^{[1]}X(r)}$

\begin{equation}
\label{eq-projective-bundle}
\PP{W}_{B^{[1]}X(r)} := B^{[2]}X(r-1)^1.
\end{equation}

\noindent
We do not assume that there exists a global vector bundle $W_{B^{[1]}X(r)}$ on 
$ B^{[1]}X(r)$ whose projectivization is $\PP{W}_{B^{[1]}X(r)}$. We will
continue this abuse of notation and denote 
$G(\PP^{k-1},\PP{W}_{B^{[1]}X(r)})$
by $G(k,{W}_{B^{[1]}X(r)})$. In section
\ref{sec-brill-noether-for-k3} $X(r)$ will be a moduli space of sheaves 
on a K3 surface $S$ and over $B^{[1]}X(r)$ we may or may not 
have a vector bundle $W_{B^{[1]}X(r)}$ whose fiber over a sheaf $F$ is a 
subspace of $H^0(S,F)$. The existence of $W_{B^{[1]}X(r)}$ depends on 
the existence of a universal sheaf. Nevertheless,  
$\PP{W}_{B^{[1]}X(r)}$ exists and is canonical. 

\begin{condition} \label{cond-compatibility-of-grassmanian-bundles}
The $G(\PP^{t-1},\PP^{n(M)+2t-2})$-bundle 
(\ref{eq-lifted-abstract-tyurin-morphism})
over $B^{[1]}M(t)$ is identified with the Grassmannian bundle 
$G(t,W_{B^{[1]}M(t)})$. 
\end{condition}

The identification in Condition
\ref{cond-compatibility-of-grassmanian-bundles}
introduces on $B^{[t+1]}M^t$ natural projectivised sub and quotient bundles 
of $\tilde{f}_{t}^*(\PP{W}_{B^{[1]}M(t)})$
which we denote by $\PP\tau_{B^{[t+1]}M^t}$
and $\PP{q}_{B^{[t+1]}M^t}$. 
Denote the relative cotangent bundle of the Grassmannian bundle
$G(t,W_{B^{[1]}M(t)})$ by 
$(\Omega^1 G)(t,W_{B^{[1]}M(t)})$ 
and its projectivization by 
$(\PP\Omega^1 G)(t,W_{B^{[1]}M(t)})$. 

\begin{condition} 
\label{cond-compatibility-of-stratifications-wrt-abstract-tyurin}
The morphisms $f_t$ and $\tilde{f}_t$ are compatible with the 
stratifications. In other words we have equality of Cartier divisors
on $B^{[t+1]}M^t$:
\[
\left[B^{[t+1]}M^t \cap B^{[t+1]}M^k\right] = 
\tilde{f}_t^{-1}\left(B^{[1]}M(t)^{k-t} \right), \ \ \mbox{for} \ \ 
t+1\leq k \leq \mu(M).
\]
\end{condition}

\noindent
(Compare Condition 
\ref{cond-compatibility-of-stratifications-wrt-abstract-tyurin} 
with Theorem \ref{thm-vainsencher} part \ref{thm-item-vain4-5}). 

The conditions above 
imply the following (compare with Theorem \ref{thm-vainsencher}
part \ref{thm-item-vain2}).

\begin{new-lemma}
\label{lemma-normal-bundle-is-twisted-relative-cotangent-bundle}
The symplectic structure of $M$ induces  a canonical isomorphism 
\begin{equation}
\label{eq-normal-bundle-is-twisted-relative-cotangent-bundle}
\phi : N_{B^{[t]}M^{t-1}/B^{[t]}M} \IsomRightArrow 
\Omega^1_{\tilde{f}_{t-1}}\left(
-\sum_{k=t}^{\mu(M)}B^{[t]}M^k
\right)
\end{equation}
between the normal bundle of $B^{[t]}M^{t-1}$ in $B^{[t]}M$ and 
the twist of the relative cotangent bundle of the Grassmannian bundle
$\tilde{f}_{t-1}: B^{[t]}M^{t-1} \rightarrow B^{[1]}M(t-1)$. 
\end{new-lemma}

\noindent
{\bf Proof:}
The case $t=1$ is trivial as $\tilde{f}_0$ is the identity. 
Assume from now on that $t\geq 2$. 
The symplectic structure $\sigma_{M}$ of $M$ pulls back to a degenerate 
$2$-form on $B^{[t]}M$ and restricts to a $2$-form 
$\sigma_{B^{[t]}M^{t-1}}$ on $B^{[t]}M^{t-1}$. 
As $\tilde{f}_{t-1}: B^{[t]}M^{t-1} \rightarrow B^{[1]}M(t-1)$
is a Grassmannian bundle, the $2$-form $\sigma_{B^{[t]}M^{t-1}}$
vanishes on each fiber and hence 
it is a pullback of a $2$-form on $B^{[1]}M(t-1)$. Contraction with 
$\sigma_M$ induces a sheaf homomorphism 
\[
\phi: N_{B^{[t]}M^{t-1}/B^{[t]}M} \rightarrow \Omega^1_{\tilde{f}_{t-1}}.
\]
It is easy to check that $\phi$ vanishes along 
$(B^{[t]}M^{t-1})\cap (B^{[t]}M^{k})$ for $k\geq t$. 
Hence, $\phi$ in
(\ref{eq-normal-bundle-is-twisted-relative-cotangent-bundle}) 
is an isomorphism if and only if 
\[
\det(N_{B^{[t]}M^{t-1}/B^{[t]}M})^*\otimes 
\omega_{\tilde{f}_{t-1}}\otimes 
\StructureSheaf{B^{[1]}M^{t-1}}\left(
-codim(M^{t-1},M)\cdot \sum_{k=t}^{\mu(M)}B^{[t]}M^k 
\right)
\]
is the trivial line-bundle.
The two short exact sequences
\begin{eqnarray*}
0 \rightarrow T_{B^{[t]}M^{t-1}} \rightarrow
& (T_{B^{[t]}M}\restricted{)}{B^{[t]}M^{t-1}} &
\rightarrow
N_{B^{[t]}M^{t-1}/B^{[t]}M} \rightarrow 0, \ \ \mbox{and} 
\\
0 \rightarrow T_{\tilde{f}_{t-1}} \rightarrow 
& T_{B^{[t]}M^{t-1}} & \rightarrow
\tilde{f}_{t-1}^*\left(T_{B^{[1]}M(t-1)}\right) \rightarrow 0
\end{eqnarray*}
imply that we have an isomorphism  
\begin{equation}
\label{eq-two-trivial-line-bundles}
\det(N_{B^{[t]}M^{t-1}/B^{[t]}M})^*\otimes \omega_{\tilde{f}_{t-1}}
\cong 
(\omega_{B^{[t]}M}\restricted{)}{B^{[t]}M^{t-1}}
\otimes \tilde{f}_{t-1}^*\left(\omega_{B^{[1]}M(t-1)}^*\right).
\end{equation}
Since $M$ and $M(t-1)$ are symplectic, we have the isomorphisms 
\[
\omega_{B^{[t]}M} \cong
\StructureSheaf{B^{[t]}M}\left(
\sum_{k=t}^{\mu(M)}
[codim(M^k,M)-1]\cdot
B^{[t]}M^k\right), 
\ \ \mbox{and}
\]
\[
\omega_{B^{[1]}M(t-1)} \cong
\StructureSheaf{B^{[1]}M(t-1)}\left(
\sum_{k=1}^{\mu(M(t-1))}
[codim(M(t-1)^k,M(t-1))-1]\cdot
B^{[1]}M(t-1)^k\right).
\]

\noindent
Condition \ref{cond-compatibility-of-stratifications-wrt-abstract-tyurin}
and equation (\ref{eq-codimension-is-sum-of-codimensions})
imply that the
right hand side of (\ref{eq-two-trivial-line-bundles})
is the line-bundle
\[
\StructureSheaf{B^{[t]}M^{t-1}}\left(
codim(M^{t-1},M)\cdot
\sum_{k=t}^{\mu(M)}B^{[t]}M^k
\right).
\] 
This completes the proof of Lemma
\ref{lemma-normal-bundle-is-twisted-relative-cotangent-bundle}.
\EndProof

\bigskip
As a consequence of the lemma we see that 
$\pi_t:B^{[t]}M^t\rightarrow B^{[t+1]}M^t$ is identified
with the projectivised homomorphism bundle
$\PP\Hom(q_{B^{[t+1]}M^t};\tau_{B^{[t+1]}M^t})$.

\begin{equation} \label{eq-normal-bundle-is-relative-cotangent-bundle}
B^{[t]}M^t = \PP\Hom(q_{B^{[t+1]}M^t},\tau_{B^{[t+1]}M^t}).
\end{equation}

\noindent
In fact, the following all denote the same space: 
\[
B^{[t]}M^t = 
\PP{N}_{B^{[t]}M^{t-1}/B^{[t]}M}= 
\PP\Omega^1_{\tilde{f}_t} = 
(\PP\Omega^1 G)(t,W_{B^{[1]}M(t)}) =
\PP\Hom(q_{B^{[t+1]}M^t};\tau_{B^{[t+1]}M^t}).
\]

\noindent
Consequently, $B^{[t]}M^t$ admits two canonical stratifications.
One stratification is induced by that of $B^{[t]}M$, namely,
\[
B^{[t]}M^t \cap 
\{
B^{[t]}M \supset B^{[t]}M^1 \supset \dots \supset B^{[t]}M^{t-1}
\}.
\]
The other stratification is the determinantal stratification indexed
by the nullity of the homomorphisms

\begin{equation} \label{eq-nullity-stratification}
\PP\Hom(q,\tau)^k := 
\{ 
\varphi  \ \ \mid \ \ 
\min[nullity(\varphi),nullity(\varphi^*)] \geq k
\}.
\end{equation}

\noindent
The inequality 
$\rank(q_{B^{[t+1]}M^t}) = 
n(M) + t - 1 \geq t = \rank(\tau_{B^{[t+1]}M^t}) 
$ 
translates (\ref{eq-nullity-stratification}) to: 

\begin{equation} \label{eq-nullity-stratification-consice}
\PP\Hom(q_{B^{[t+1]}M^t};\tau_{B^{[t+1]}M^t})^k := 
\{ 
\varphi  \ \ \mid \ \ 
nullity(\varphi) \geq 
k+ n(M) - 1    \ \
\}, \ \ 0 \leq k \leq t-1.
\end{equation}

\begin{condition} 
\label{cond-compatibility-of-stratifications-on-proper-transforms}
The two stratifications coincide scheme theoretically:
\[
\PP\Hom(q_{B^{[t+1]}M^t};\tau_{B^{[t+1]}M^t})^k = 
\left(B^{[t]}M^t \cap B^{[t]}M^k \right), \ \ 0 \leq k \leq t-1. 
\]
\end{condition}

\noindent
(Compare Condition 
\ref{cond-compatibility-of-stratifications-on-proper-transforms}
with Theorem \ref{thm-vainsencher} part 
\ref{thm-item-vain6}).

\begin{defi}
{\rm
\label{def-admisible-stratified-collection}
A collection of stratified quasi-projective schemes 
\begin{equation}
\label{eq-abstract-admissible-collection}
\{ X(r)^t, f_{r,t}\mid \ \ 0 \leq r \leq \mu, \ \ 0 \leq t \leq \mu-r \}
\end{equation}
satisfying all the conditions above will be 
called a
{\em symplectically dualizable stratified collection
}
or a 
{\em dualizable collection
} 
for short.
}
\end{defi}

Note that if the collection (\ref{eq-diagram-of-X-v}) is  dualizable, then 
the sub-collection obtained by deleting the top row
is also dualizable.

\begin{thm}
\label{thm-stratified-elementary-transformation}
Given a symplectically dualizable stratified collection
$\{X(r)^t, f_{r,t}\}$ with $X=X(0)^0$, 
there exists another dualizable collection
$\{ Y(r)^t , f'_{r,t}\}$ with $Y=Y(0)^0$, $n(Y)=n(X)$ and $\mu(Y)=\mu(X)$
satisfying
\begin{enumerate}
\item 
The full iterated blow-ups are isomorphic
\[
\tilde{q}_r :B^{[1]}X(r) \IsomRightArrow B^{[1]}Y(r), \ \ 0 \leq r \leq \mu.
\]
These isomorphisms restrict to  isomorphisms on the exceptional divisors
\[
\tilde{q}_{r,t} : B^{[1]}X(r)^t  \IsomRightArrow 
B^{[1]}Y(r)^t, \ \ 0 \leq t \leq \mu-r.
\]
\item
The two 
$\PP^{n+2r-2}$-bundles $\PP{W}_{B^{[1]}X(r)}$ and $\PP{W}_{B^{[1]}Y(r)}$ 
\begin{eqnarray*}
\tilde{f}_{r-1,1} & : B^{[2]}X(r-1)^1 \rightarrow B^{[1]}X(r) \\
\tilde{f}'_{r-1,1}& : B^{[2]}Y(r-1)^1 \rightarrow B^{[1]}Y(r) 
\end{eqnarray*}
are dual. In other words, $\tilde{q}_r$ lifts to an isomorphism
\begin{equation}
\label{eq-projective-bundles-are-dual}
\DoubleTilde{q}_r \ : \  \PP{W}_{B^{[1]}X(r)} \ \LongIsomRightArrow \ 
\PP{W}^*_{B^{[1]}Y(r)}.
\end{equation}
\item
We have a commutative diagram of isomorphisms (for $1\leq t \leq \mu-r$) 
\[
{
\divide\dgARROWLENGTH by 2
\begin{diagram}
\node{B^{[1]}X(r)^t}
\arrow{e,t}{\cong}
\arrow{s,lr}{\cong}{\tilde{q}_{r,t}}
\node{(B^{[1]}\PP\Omega^1 G)(t,W_{B^{[1]}X(r+t)})}
\arrow{s,r}{\cong}
\\
\node{B^{[1]}Y(r)^t}
\arrow{e,t}{\cong}
\node{(B^{[1]}\PP\Omega^1 G)(t,W_{B^{[1]}Y(r+t)})}
\end{diagram}
}
\]
where the right vertical isomorphism is the 
relative transposition isomorphism (see section \ref{sec-transposition}). 
\end{enumerate}
\end{thm}

We will refer to the collection $\{ Y(r)^t , f'_{r,t}\}$ as
the {\em  dual collection}. 
The proof of Theorem \ref{thm-stratified-elementary-transformation} 
is carried out in 
section \ref{sec-construction-of-the-dual-collection}.

\section{Determinantal Varieties - Background} 
\label{sec-determinantal-varieties}

\subsection{Blowing up determinantal ideals}
\label{subsec-blowing-up-det-ideals}

We recall Vainsencher's results about blowing up determinantal ideals. 
Let $V_0$ and $V_1$ be 
vector bundles over a scheme $M$ of finite type over 
$\ComplexNumbers$. Given a vector 
bundle $E$ over $M$ 
we denote the tautological exact sequence over the 
Grassmanian bundle $G(i,E)$ by 
\[
0 \rightarrow 
\tau_i(E) \rightarrow 
E_{G(i,E)} \rightarrow 
q_i(E) \rightarrow 0.
\]
Given another vector bundle $F$ we denote by $\SheafHom(q_i(E),\tau_j(F))$ 
the corresponding vector bundle over $G(i,E)\times_{M} G(j,F)$. 
We index the determinantal loci in $\PP\SheafHom(V_0,V_1)$
as in (\ref{eq-nullity-stratification}) in terms of co-rank. 
(Inxeding by rank would simplify the notation in Theorem \ref{thm-vainsencher}
but would be cumbersome in the context of 
Brill-Noether stratifications of moduli spaces of sheaves.) 
Denote the intersection $B^{[i]}\PP\SheafHom(V_0,V_1)^{j} \cap 
B^{[i]}\PP\SheafHom(V_0,V_1)^{k}$ by $B^{[i]}\PP\SheafHom(V_0,V_1)^{j\cap k}$.

\begin{thm} 
\label{thm-vainsencher}
(\cite{vainsencher} Theorem 2.4)
\begin{enumerate}
\item
\label{thm-item-vain1}
$
B^{[i+1]}\PP\SheafHom(V_0,V_1)^i \ = \ 
B^{[i+1]}\PP\SheafHom(q_{r_0-r_1+i}(V_0),\tau_{r_1-i}(V_1))
$
for $0\leq i < r_1$.
\item
\label{thm-item-vain2}
Let $\Ideal{i}$ be the product of the invertible ideals 
$\Ideal{}B^{[i\!+\!1]}\PP\SheafHom(V_0,V_1)^k$, 
$k=i\!+\!1, \dots ,r_1\!-\!1$.
Then the conormal sheaf $N^*_i$ of 
$B^{[i+1]}\PP\SheafHom(V_0,V_1)^i$ in $B^{[i+1]}\PP\SheafHom(V_0,V_1)$
is isomorphic to 
$\Ideal{i}\cdot\StructureSheaf{}(1)\otimes
\SheafHom(\tau_{r_0-r_1+i}(V_0),q_{r_1-i}(V_1))$.
\item
\label{thm-item-vain3}
Let $Z(i)$ be $G(r_0\!-\!i,V_0) \times_M G(i,V_1)$. 
For $i\leq k$, 
$
B^{[i]}\PP\SheafHom(V_0,V_1)^k 
\ = \ 
$
\\
$
B^{[i]}\PP\SheafHom(\tau_{r_0-r_1+k}(V_0),q_{r_1-k}(V_1))
\ \times_{Z(r_1-k)} \ 
B^{[1]}\PP\SheafHom(q_{r_0-r_1+k}(V_0),\tau_{r_1-k}(V_1)).
$
\\
Moreover, under this identification, we have 
$
N^*_{B^{[k]}\PP\SheafHom(V_0,V_1)^k} \ = \ 
\StructureSheaf{}(1)\otimes \Ideal{k}^{-1}. 
$
\item
\label{thm-item-vain4-5}
For $k > i \geq 0$, 
$
B^{[i\!+\!1]}\PP\SheafHom(V_0,V_1)^{k \cap i}
\ = \ 
B^{[i\!+\!2]}\PP\SheafHom
(q_{r_0\!-\!r_1\!+\!i}(V_0),\tau_{r_1\!-\!i}(V_1))^{k\!-\!i} \ =
$
\\
$
B^{[i\!+\!1]}\PP\SheafHom(\tau_{r_0\!-\!r_1\!+\!k}(V_0),q_{r_1\!-\!k}(V_1))^i
\ \times_{Z(r_1-k)} \
B^{[2]}\PP\SheafHom(q_{r_0\!-\!r_1\!+\!k}(V_0),\tau_{r_1\!-\!k}(V_1)).
$
\item
\label{thm-item-vain6}
$
B^{[i]}\PP\SheafHom(V_0,V_1)^{k \cap i}
\ = 
$
\\
$
\PP\SheafHom(\tau_{r_0\!-\!r_1\!+\!i}(V_0),q_{r_1\!-\!i}(V_1))^k \ 
\times_{Z(r_1-i)} \
B^{[1]}\PP\SheafHom
(q_{r_0\!-\!r_1\!+\!i}(V_0),\tau_{r_1\!-\!i}(V_1)) 
$
for $i>k$.
\item
\label{thm-item-vain8}
$
\left(\Ideal{}B^{[i\!+\!1]}\PP\SheafHom(V_0,V_1)^k\right)\cdot
\StructureSheaf{B^{[i]}\PP\SheafHom(V_0,V_1)} 
\ = \ 
\left(\Ideal{}B^{[i]}\PP\SheafHom(V_0,V_1)^k\right)\cdot
\left(\Ideal{}B^{[i]}\PP\SheafHom(V_0,V_1)^i\right)^{i\!-\!k\!+\!1}
$
\\
for $i>k$.
\item
\label{thm-item-vain9-10}
For $k>j$, 
$
B^{[j\!-\!i]}\PP\SheafHom(V_0,V_1)^{k \cap j} \ = 
$
\\
$ 
B^{[j\!-\!i]}\PP\SheafHom(\tau_{r_0\!-\!r_1\!+\!j}(V_0),q_{r_1\!-\!j}(V_1))
\ \times_{Z(j)} \
B^{[2]}\PP\SheafHom(q_{r_0\!-\!r_1\!+\!j}(V_0),\tau_{r_1\!-\!j}(V_1))^{k\!-\!j}
\ = $
\\
$
B^{[j\!-\!i]}\PP\SheafHom(\tau_{r_0\!-\!r_1\!+\!k}(V_0),q_{r_1\!-\!k}(V_1))^j
\ \times_{Z(r_1\!-\!k)} \
B^{[2]}\PP\SheafHom(q_{r_0\!-\!r_1\!+\!k}(V_0),\tau_{r_1\!-\!k}(V_1)).
$
\end{enumerate}
\end{thm}

The iterated blow-up $B^{[k\!+\!1]}\PP\Hom(V_0,V_1)^k$ has several beautiful 
modular interpretations: as a space of complete collineations,
as a closure in the fiber product of Grassmannian bundles
\cite{vainsencher,laksov,kleiman-thorup}, and three new modular interpretations
discovered in \cite{thadeus-colineations}. 

\subsection{Transpositions}
\label{sec-transposition}

Let $H$ be an $h$-dimensional  vector space. 
The cotangent bundle $T^*G(t,H)$ of the
Grassmannian of $t$-dimensional subspaces of $H$ is isomorphic to the
homomorphism bundle
$\Hom(q,\tau)$. We obtain a natural determinantal stratification of 
$\PP{T}^*G(t,H)$. Similarly, we have a natural determinantal 
stratification of $\PP{T}^*G(t,H^*)$. Assume that $h \geq 2t$. 
We construct in this section a canonical isomorphism
\[
B^{[1]}\PP{T}^*G(t,H) \IsomRightArrow
B^{[1]}\PP{T}^*G(t,H^*)
\]
between the two iterated blow-ups. 
Notice that $Flag(t,h\!-\!t,H)$ and $Flag(t,h\!-\!t,H^*)$ are isomorphic and
under this identification we have the isomorphisms 
$\tau_t(H) \cong q_{h\!-\!t}^*(H^*)$ and
$\tau_{h\!-\!t}(H) \cong q_t^*(H^*)$. The birational transposition isomorphism 
``factors'' through 
\begin{equation}
\label{eq-factorization-of-birational-iso-through-flag}
\PP{T}^*G(t,H) 
\ \leftarrow \
\PP\left([\tau_t(H)\otimes\tau_t(H^*)\restricted{]}{Flag(t,h\!-\!t,H)}\right)
\ \rightarrow \ 
\PP{T}^*G(t,H^*). 
\end{equation}

%
%

The whole construction works in the relative setting in which $H$ is a 
vector bundle over a base $Z$. 
For each 
pair of integers $(t,j)$ satisfying
$0 \leq t \leq \frac{h}{2}$ and $0 \leq j \leq t+1$, we get 
\[
(B^{[j]}\PP\Omega^1 G)(t,H) \rightarrow  G(t,H)  \rightarrow Z.
\]
For $j = t$ and $j=t+1$ we set
\begin{eqnarray*}
(B^{[t]}\PP\Omega^1 G)(t,H) &=& 
(\PP\Omega^1 G)(t,H), \ \ \mbox{and} 
\\
(B^{[t+1]}\PP\Omega^1 G)(t,H) &=& G(t,H).
\end{eqnarray*}
For $0 \leq j \leq t-1$, 
$(B^{[j]}\PP\Omega^1 G)(t,H)$ is the iterated blow-up of the determinantal
stratification of $(\PP\Omega^1 G)(t,H)$. 

\begin{prop}
\label{prop-transposition-for-cotangent-of-grassmannian-bundle}
$
(B^{[1]}\PP\Omega^1 G)(t,H) \ \cong \
B^{[1]}\PP\left([\tau_t(H^*)\otimes \tau_t(H)
\restricted{]}{Flag(t,h\!-\!t,H)} \right) 
\ \cong \
(B^{[1]}\PP\Omega^1 G)(t,H^*). 
$
\end{prop}

\noindent
{\bf Proof:}
Theorem \ref{thm-vainsencher} part \ref{thm-item-vain1} (with $i=0$,
$V_0=q_t(H)$, and $V_1=\tau_t(H)$) 
implies the identity 
\[
B^{[1]}\PP\SheafHom(q_t(H),\tau_t(H)) \ = \
B^{[1]}\PP\SheafHom(q_{h\!-\!2t}(q_t(H)),\tau_t(H)). 
\]
But $Flag(t,h\!-\!t,H)$ is isomorphic to $G(h\!-\!2t,q_t(H))$. 
\EndProof

\medskip
We denote by 
\[
W_{(B^{[1]}\PP\Omega^1G)(t,H)}
\]
the rank $h-2t$ bundle over the full iterated blow-up 
$(B^{[1]}\PP\Omega^1 G)(t,H)$ obtained by pulling back the bundle
$\tau_{h\!-\!t}/\tau_t$ from $Flag(t,h\!-\!t,H)$. 

\medskip
Given a third integer $i$ satisfying $0 \leq i \leq t-1$ we get the variety
$
(B^{[j]}\PP\Omega^1 G)(t,H)^i 
$
which is the strict transform of the strata $(\PP\Omega^1 G)(t,H)^i$ under the
iterated blow-up. 
The following identity will be useful in section 
\ref{sec-construction-of-the-dual-collection}:

\begin{new-lemma}
\label{lemma-inner-square-commutes}
There is a natural isomorphism for $j \leq i+1$
\[
(B^{[j]}\PP\Omega^1G)(t,H)^i \cong
(B^{[j]}\PP\Omega^1G)(i,W_{(B^{[1]}\PP\Omega^1G)(t-i,H)}).
\]
In particular, for $j=i$ and $j=i+1$ we have a commutative diagram 
in which the horizontal morphisms are isomorphisms

\begin{equation}
\label{diagram-inner-square-commutes}
{
\divide\dgARROWLENGTH by 2
\begin{diagram}
\node{(B^{[i]}\PP\Omega^1G)(t,H)^i}
\arrow{e,t}{\cong}
\arrow{s,l}{\beta_{i+1}}
\node{(\PP\Omega^1G)(i,W_{(B^{[1]}\PP\Omega^1G)(t-i,H)})}
\arrow{s,r}{\beta_{i+1}}
\\
\node{(B^{[i+1]}\PP\Omega^1G)(t,H)^i}
\arrow{e,t}{\cong}
\node{G(i,W_{(B^{[1]}\PP\Omega^1G)(t-i,H)}).}
\end{diagram}
}
\end{equation}
\end{new-lemma}

\noindent
{\bf Proof:} 
Theorem \ref{thm-vainsencher} part \ref{thm-item-vain3}, 
and the identity 
$G(t\!-\!i,\tau_t(H))\ \times_{Z}G(h\!-\!2t\!+\!i,q_t(H)) \ = \ 
Flag(t\!-\!i,t,h\!-\!t\!+\!i,H)$,
imply that
$(B^{[j]}\PP\Omega^1G)(t,H)^i$ is isomorphic to \\
$
B^{[j]}\PP\SheafHom(\tau_{h\!-\!2t\!+\!i}(q_t(H)),q_{t\!-\!i}(\tau_t(H))) 
\ \times_{Flag(t\!-\!i,t,h\!-\!t\!+\!i,H)} \
B^{[1]}\PP\SheafHom(q_{h\!-\!t\!+\!i}(H),\tau_{t\!-\!i}(H)).
$
\\
The latter is $B^{[j]}\PP\Omega^1
G\left(i,W_{B^{[1]}\PP\SheafHom(q_{h\!-\!t\!+\!i}(H),\tau_{t\!-\!i}(H))}
\right). $
Theorem \ref{thm-vainsencher} part \ref{thm-item-vain1} (with $i=0$)
identifies $B^{[1]}\PP\SheafHom(q_{h\!-\!t\!+\!i}(H),\tau_{t\!-\!i}(H))$
with $(B^{[1]}\PP\Omega^1G)(t\!-\!i,H)$. 
\EndProof

\medskip
Theorem \ref{thm-vainsencher} and 
Lemma \ref{lemma-inner-square-commutes} imply that the stratified collection 

\begin{equation}
\label{eq-dualizable-diagram-for-a-proj-cotangent-bundle-of-grass}
\begin{array}{cccccccc}
P\Omega^1G(t,H) \supset & P\Omega^1G(t,H)^1           &     \cdots              
& \supset P\Omega^1G(t,H)^{t-1}
\\
                   &  \downarrow \pi_{t,1}  &                   
\\
                   & P\Omega^1G(t-1,H) \supset & P\Omega^1G(t-1,H)^1   
& \cdots P\Omega^1G(t-1,H)^{t-2}
\\
                   &                        & \downarrow \pi_{t-1,1} 
\\
                   &                        &  P\Omega^1G(t-2,H) \supset  
& \cdots P\Omega^1G(t-2,H)^{t-3}
\\
                   &                        &     
& \vdots
\\
                   &                        &     
& \downarrow \pi_{2,1}
\\
                   &                        &     
& P\Omega^1G(1,H)
\end{array}
\end{equation}

\noindent
is  dualizable 
in a sense analogous to that of Definition 
\ref{def-admisible-stratified-collection}. 
Proposition \ref{prop-transposition-for-cotangent-of-grassmannian-bundle} 
identifies the dual collection as the one associated to $H^*$.
Note, however, that
the diagonal entries of 
(\ref{eq-dualizable-diagram-for-a-proj-cotangent-bundle-of-grass}) 
are (fibrations over the base scheme $Z$ of) contact 
(rather than symplectic) varieties. 

\bigskip
Denote by $\tau'_{G(t,H^*)}$ and 
$q'_{G(t,H^*)}$ the tautological sub- and quotient bundles. 
Let $h:=h_{(\PP\Omega^1G)(t,H)}$ 
(resp. $h':=h_{(\PP\Omega^1G)(t,H^*)}$) be the tautological line sub-bundle of
$\beta_t^*\Hom(q,\tau)$ (resp. $\beta'^*_t\Hom(q',\tau')$). 


\begin{new-lemma}
\label{lemma-two-tautological-h-subbundles-are-isomorphic}
Over $(B^{[1]}\PP\Omega^1G)(t,H)$ we have a natural isomorphism 
\[
\beta^*h_{(\PP\Omega^1G)(t,H)} \cong \beta'^* h_{(\PP\Omega^1G)(t,H^*)}.
\]
\end{new-lemma}

\noindent
{\bf Proof:}
The bundle $\tau'_t\otimes\tau_t$ over $Flag(t,h\!-\!t,H)$ is a subbundle of 
the pullback of both $\Hom(q,\tau)$ and $\Hom(q',\tau')$. 
Working over the three spaces in
(\ref{eq-factorization-of-birational-iso-through-flag}) 
we see that both $h$ and $h'$ pullback to the tautological line subbundle of 
the pullback of $\tau'_t\otimes\tau_t$ to $\PP[\tau'_t\otimes\tau_t]$. 
\EndProof


\medskip
The proof of Theorem \ref{thm-correspondence-is-a-linear-combination}
will require an extension of 
Proposition \ref{prop-transposition-for-cotangent-of-grassmannian-bundle}. 
If the base scheme $Z$ is a point, the cotangent bundle of  $G(t,H)$ 
admits a non-trivial extension $E(H)$ 
(see (\ref{eq-idempotent-extension-of-cotangent-bundle})). 
The complement $[\PP{E}(H)\setminus (\PP\Omega^1G)(t,H)]$ 
parametrizes $\ComplexNumbers^{\times}$-orbits of idempotents and maps 
isomorphically onto the Zariski open subset in $G(t,H)\times G(h\!-\!t,H)$ of 
decompositions of $H$. We have a relative analogue of the  extension 
(\ref{eq-idempotent-extension-of-cotangent-bundle}) over any base $Z$. 
Over $G(t,H)\times_{Z}G(h\!-\!t,H)$ we have a canonical homomorphism
\begin{equation}
\label{eq-homomorphism-over-product-of-grassmannians}
\alpha : \tau_t(H) \ \rightarrow \ q_{h\!-\!t}(H).
\end{equation}

\noindent
It is the composition of the injection $\tau_t(H)\hookrightarrow H$
with the projection $j:H \rightarrow q_{h\!-\!t}(H)$. If the pair 
$\tau_t(H)$ and $\tau_{h\!-\!t}(H)$ provides a decomposition of $H$, then
$\alpha$ is an isomorphism and $(j\circ \alpha^{-1})$ is the projection 
to $\tau_t(H)$. 
The restriction of $\alpha$ to the flag variety vanishes. 
The section $\alpha$ is transversal to the (affine) determinantal 
stratification of the affine bundle $\SheafHom(\tau_t(H),q_{h\!-\!t}(H))$. 
We get a determinantal stratification of $G(t,H)\times_{Z}G(h\!-\!t,H)$. 

\begin{prop}
\label{prop-extended-transposition}
There is a canonical isomorphism
\[
B^{[1]}\PP{E}(H) \ \cong \ 
B^{[1]}[G(t,H)\times_{Z}G(t,H^*)] \ \cong \ 
B^{[1]}\PP{E}(H^*).
\]
\end{prop}

\noindent
{\bf Proof:} 
$B^{[1]}\PP{E}(H)$ is smooth because the stratification of $\PP{E}(H)$ 
is supported on its divisor ${\PP}T^*G(t,H)$. 
$B^{[1]}[G(t,H)\times_{Z}G(t,H^*)]$
is smooth because each determinantal stratum in $G(t,H)\times_{Z}G(t,H^*)$
is smooth away from the next stratum and has the expected dimension. 
If we regard $\PP{E}(H)$ as a subbundle of $\PP[H^*\otimes\tau_t]$, then 
Theorem \ref{thm-vainsencher} part \ref{thm-item-vain1}
implies that there is a regular morphism from $B^{[1]}\PP{E}(H)$
to $G(t,H)\times_{Z}G(t,H^*)$. An inductive application of the 
universal property of blowing-up implies that the morphism 
lifts to $B^{[1]}[G(t,H)\times_{Z}G(t,H^*)]$. 
It is easy to check that the lift is an isomorphism away from 
codimension two. By purity of the ramification locus, the lift must be
an isomorphism.
\EndProof

\subsection{The Petri map}
\label{subsec-petri-map}

In this section we assume that $M$ is a smooth algebraic variety. 
Let 
\begin{equation}
\label{eq-homomorphism}
e : V_0 \rightarrow V_1
\end{equation}
be a homomorphism of vector bundles over $M$. 
Denote by $r_i$ the rank of $V_i$ and 
Assume that $e$ is generically surjective (in particular 
the inequality $r_0\geq r_1$). 
Let
\begin{equation}
\label{eq-determinantal-stratification-of-e}
M = M^0 \supset M^1 \supset \cdots \supset M^{\mu} \supset \emptyset
\end{equation}
be the determinantal stratification. 
Over $M^t\setminus M^{t+1}$ we have two vector bundles
$\ker(\restricted{e}{M^t\setminus M^{t+1}})$ and
$\coker(\restricted{e}{M^t\setminus M^{t+1}})$. 

\begin{new-lemma}
\label{lemma-existence-of-petri-map}
\cite{fulton,a-c-g-h}
\begin{enumerate}
\item
\label{lemma-item-existence-of-petri-map}
There exists a canonical  homomorphism of 
$\StructureSheaf{M^t\setminus M^{t+1}}$-modules 
\begin{equation}
\phi: TM\otimes \ker(\restricted{e}{M^t\setminus M^{t+1}}) \rightarrow 
\coker(\restricted{e}{M^t\setminus M^{t+1}}). 
\end{equation}
Given a point $f$ in the fiber of $\ker(\restricted{e}{M^t\setminus M^{t+1}})$ 
at $x\in M^t\setminus M^{t+1}$ and a tangent vector 
$\xi\in T_{x}\restricted{M}{M^t\setminus M^{t+1}}$, 
their image  $\phi(\xi\otimes f)$ 
is the first order infinitesimal obstruction to extending $f$ as a section of 
$\ker(e)$ in the direction of $\xi$. 
\item
\label{lemma-item-petri-map-identifies-conormal-bundle}
If $M^t$ has the expected dimension and
\[
G(r_0-r_1+t,M^t) := \{(x,W)\mid \ W \in G(r_0-r_1+t,\ker(e_x))\}
\]
is smooth, then the tangent cone of $M^t$ at
$x\in M^k\setminus M^{k+1}$, 
$k\geq t$, is the subscheme of $T_xM$ of zeroes of the section
$\Wedge{k+1-t}\Phi_x$. Here we interpret $\phi_x$ as a section
$\Phi_x$ of the vector bundle 
$Hom(\ker(e_x),\coker(e_x))\otimes\StructureSheaf{T_xM}$
over $T_xM$.  
\item 
\label{lemma-item-petri-map-of-dual-homo-is-the-same}
As a section of 
\[
\ker(\restricted{e}{M^t\setminus M^{t+1}})\otimes
\coker(\restricted{e}{M^t\setminus M^{t+1}})^* \rightarrow
T^*\restricted{M}{M^t\setminus M^{t+1}},
\]
$\phi$ is equal to the section obtained from the dual homomorphism
\[
e^*: V^*_1 \rightarrow V^*_0.
\]
\end{enumerate}
\end{new-lemma}

%
%
%
%
%

Part \ref{lemma-item-petri-map-identifies-conormal-bundle} of the Lemma
indicates that the expected dimension of $M^t$ is
\begin{equation}
\label{eq-expected-codimension}
\rho(t):= \dim(M)- t(r_0-r_1+t). 
\end{equation}

\subsection{Blowing up the smallest stratum}
\label{subsec-blowing-up-smallest-strata}

Let 
\begin{equation}
\label{eq-pulled-bach-homomorphism}
\beta_{\mu}: B^{\mu}M^t \rightarrow M^t, \ 0\leq t \leq \mu-1,
\end{equation}
be the blow-up of $M^t$ along $M^{\mu}$. 
Pulling back (\ref{eq-homomorphism}) we get the homomorphism
\begin{equation}
\label{eq-pulled-back-homomorphism-via-blow-up-map}
\beta_{\mu}^*(e): \beta_{\mu}^*V_0 
\rightarrow
\beta_{\mu}^*V_1.
\end{equation}
It defines a determinantal stratification
\[
B^{\mu}M  = \beta_{\mu}^{-1}M^0 \supset \beta_{\mu}^{-1}M^1 
\supset \cdots \supset 
\beta_{\mu}^{-1}M^{\mu} \supset \emptyset
\]
which is the {\em total}-transform of the stratification
(\ref{eq-determinantal-stratification-of-e}).

Next we perform an elementary transformation of 
(\ref{eq-pulled-back-homomorphism-via-blow-up-map}) 
whose determinantal stratification
is the {\em strict}-transform of 
(\ref{eq-determinantal-stratification-of-e}). 
Let 
\begin{eqnarray*}
B^{\mu}V_0 & := & \beta_{\mu}^*V_0 
\\
B^{\mu}V_1 & := & \ker\left(
\beta_{\mu}^*V_1 \rightarrow \beta_{\mu}^*\coker(\restricted{e}{M^{\mu}})
\right).
\end{eqnarray*}

\noindent
$E:=\beta_{\mu}^{-1}M^{\mu}$ is a Cartier divisor on $B^{\mu}M$ and
$\beta_{\mu}^*(e)$ has constant rank over $E$. 
Thus, $B^{\mu}V_1$ is locally free. 
We get a natural homomorphism

\begin{equation}
\label{eq-blown-up-homomorphism}
B^{\mu}(e):
B^{\mu}V_0
\rightarrow
B^{\mu}V_1.
\end{equation}

\noindent
It defines a determinantal stratification

\begin{equation}
\label{eq-blown-up-stratification}
B^{\mu}M  = B^{\mu}M^0 \supset B^{\mu}M^1 \supset \cdots \supset 
B^{\mu}M^{\mu-1} \supset (B^{\mu}M)^{\mu}(B^{\mu}e) \supset \emptyset.
\end{equation}

\noindent
If $\mu=1$ and $n(M)=1$ then $M^1$ is already a divisor on $M$.
In that case $\beta_{\mu}:B^{\mu}M\rightarrow M$ is the identity,
however $B^{\mu}V_1$ is different from $V_1$ and
$B^{\mu}(e)$ is different from $e$. 

\begin{new-lemma}
\label{lemma-stratification-of-blown-up-homomorphism-is-the-proper-trensform}
The subschemes $B^{\mu}M^t$, $0\leq t \leq \mu-1$ in the 
stratification
(\ref{eq-blown-up-stratification})
are the strict transforms of the subschemes 
$M^t$ in the stratification 
(\ref{eq-determinantal-stratification-of-e}).
\end{new-lemma}

\noindent
In general, the subscheme $(B^{\mu}M)^{\mu}(B^{\mu}e)$
may or may not be empty. 

\medskip
\noindent
{\bf Proof:}
The Lemma holds 
for a general determinantal stratification
over a scheme of finite type over $\ComplexNumbers$ 
(without assuming smoothness). 
The Lemma is a reformulation of  Theorem \ref{thm-vainsencher}
part \ref{thm-item-vain8}. 
The homomorphism $e$ fits in a commutative diagram:

\begin{equation}
{
\divide\dgARROWLENGTH by 2
\begin{diagram}
\node{
B^{\mu}M}
\arrow{e,t}{\widetilde{\PP(e)}}
\arrow{s,l}{\beta}
\node{B^{\mu}\PP\Hom(V_0,V_1)}
\arrow{s}
\\
\node{M} 
\arrow{e,b}{\PP{e}}
\node{\PP\Hom(V_0,V_1)}
\end{diagram}
}
\end{equation}

\noindent
Above, $\PP{e}$ is a closed immersion and its image is disjoint from 
the loci $\PP\Hom(V_0,V_1)^i$ for all $i > \mu$. 
The section $e$ is the pullback 
of the tautological section $\tilde{e}$. 
The morphism $\widetilde{\PP(e)}$ is the canonical closed immersion 
(see \cite{hartshorne} Corollary 7.15 page 165). 
The determinantal ideal-sheaves  on $M$ and and their strict and total 
transforms on $B^{\mu}M$ are inverse images of the 
corresponding ideals on $\PP\Hom(V_0,V_1)$ and 
$B^{\mu}\PP\Hom(V_0,V_1)$ via the closed immersions $\PP{e}$ and 
$\widetilde{\PP(e)}$. Hence, it suffices to prove the Lemma 
with $M$ equal to the complement of $\PP\Hom(V_0,V_1)^{\mu+1}$ in 
the bundle $\PP\Hom(V_0,V_1)$ 
(renaming $V_1\otimes\StructureSheaf{}(1)$ by $V_1$ and $\tilde{e}$ by $e$). 

The determinantal locus $M^t$ is the zero subscheme of $\Wedge{r_1-t+1}e$.
Its ideal sheaf $\Ideal{M^t}$ is the image of
\[
(\Wedge{r_1-t+1}V_0) \ \otimes \ (\Wedge{r_1-t+1}V_1)^* \ \ 
\LongRightArrowOf{\Wedge{r_1-t+1}e} \ \
\StructureSheaf{M}.
\]
The ideal of the total transform $\beta^{-1}\Ideal{M^t}$ is the 
image of $\Wedge{r_1-t+1}\beta^*(e)$. Let $\Ideal{t}(B^{\mu}e)$ be the 
ideal sheaf which is the image of $\Wedge{r_1-t+1}B^{\mu}e$ in 
$\StructureSheaf{B^{\mu}M}$ and let $\Ideal{B^{\mu}M^t}$ be the ideal of 
the strict transform $B^{\mu}M^t$ of $M^t$. 
We need to prove the equality 
$\Ideal{B^{\mu}M^t} = \Ideal{t}(B^{\mu}e)$. The equality follows from
the two equalities
\begin{eqnarray*}
\beta^{-1}\Ideal{M^t} & = & \Ideal{t}(B^{\mu}e)\cdot (\Ideal{E})^{\mu-t+1}
\ \ \ \mbox{and} 
\\
\beta^{-1}\Ideal{M^t} & = & \Ideal{B^{\mu}M^t} \cdot (\Ideal{E})^{\mu-t+1}.
\end{eqnarray*}
The first equality follows from a short local consideration. 
The second equality follows from 
Theorem \ref{thm-vainsencher} part \ref{thm-item-vain8}.
\EndProof

\bigskip
Denote by
\[
E^t := E \cap B^{\mu}M^t, \ \ 0\leq t \leq \mu, 
\]
the induced stratification of $E$. 
We describe below the determinantal stratification 
(\ref{eq-intersection-with-proper-transform-stratification}) 
of the exceptional divisor
under the  assumption that the smallest determinantal 
locus $M^{\mu}$ is smooth and of the expected dimension. 
We will see that
the determinantal stratification 
(\ref{eq-blown-up-stratification})
of the blown-up homomorphism $B^{\mu}(e)$ is shorter
than the determinantal stratification of $e$. 
In other words, $(B^{\mu}M)^\mu(B^{\mu}e)$ is empty.

Assume that $M$ is smooth  and the lowest strata
$M^{\mu}$ of (\ref{eq-determinantal-stratification-of-e})
is smooth of the expected codimension $\mu(r_0-r_1+\mu)$. 
By Lemma
\ref{lemma-existence-of-petri-map},
The Petri map of $e$ is an isomorphism
\begin{equation}
\label{eq-homomorphism-induced-by-petri-map}
\phi:\Hom(\coker(\restricted{e}{M^{\mu}}),\ker(\restricted{e}{M^{\mu}})) 
\IsomRightArrow N^*_{M^{\mu}/M}.
\end{equation}

\noindent
In particular, 
the exceptional divisor 
\[
E:=\beta_{\mu}^{-1}M^{\mu}
\] 
in $B^{\mu}M$ is isomorphic to 
the projectivized homomorphism bundle
\begin{equation}
\label{eq-exceptional-divisor-is-a-projectivized-hom-bundle}
E \cong 
\PP\Hom(\ker(\restricted{e}{M^{\mu}}),\coker(\restricted{e}{M^{\mu}})).
\end{equation}

\noindent
The exceptional divisor $E$ admits two determinantal
stratifications:

\begin{equation}
\label{eq-intersection-with-proper-transform-stratification}
E \cap \{ B^{\mu}M = B^{\mu}M^0 \supset B^{\mu}M^1 \supset \cdots
\supset B^{\mu}M^{\mu-1} \supset (B^{\mu}M)^{\mu}(B^{\mu}e)
\} \ \ \mbox{and}
\end{equation}

\begin{equation}
\label{eq-hom-stratification}
E^k :=  
\{\eta\mid \ nullity(\eta) \geq k+r_0-r_1\}.
\end{equation}

\noindent
The equality of the two stratifications 
(\ref{eq-intersection-with-proper-transform-stratification}) and 
(\ref{eq-hom-stratification}) is proven in Theorem 
\ref{thm-vainsencher}  part \ref{thm-item-vain6} (see also Lemma 
\ref{lemma-compatibility-of-two-stratifications-on-exceptional-divisor}
part \ref{lemma-item-two-stratifications-coincide} below). 
Denote by 
$\StructureSheaf{\PP\Hom(
\ker(\restricted{e}{M^{\mu}}),
\coker(\restricted{e}{M^{\mu}}))}(1)$ 
the relative ample line-bundle and let 

\begin{equation}
\eta: 
\beta_{\mu}^* \ker(\restricted{e}{M^{\mu}})
\rightarrow 
\left(\beta_{\mu}^*\coker(\restricted{e}{M^{\mu}})\right)(1)
\end{equation}

\noindent
be the tautological homomorphism. 

\begin{new-lemma}
\label{lemma-compatibility-of-two-stratifications-on-exceptional-divisor}
\begin{enumerate}
\item
\label{lemma-item-two-tautological-line-bundles-isomorphic}
The two line bundles 
$\StructureSheaf{\PP\Hom(
\ker(\restricted{e}{M^{\mu}}),
\coker(\restricted{e}{M^{\mu}}))}(1)$ 
and
$\StructureSheaf{E}(-E)$ are canonically isomorphic.
\item
\label{lemma-item-composition-is-projection}
The composition
\[
\beta^*V_{\restricted{0}{E}} 
\RightArrowOf{B^{\mu}e}
B^{\mu}V_{\restricted{1}{E}} 
\rightarrow
Im(\beta^*{\restricted{e}{E}}) 
\cong
\beta^*V_{\restricted{0}{E}}/ \ker(\beta^*{\restricted{e}{E}})
\]
is the natural projection.
\item
\label{lemma-item-blown-up-homomorphism-restricts-to-eta}
The restriction of $B^{\mu}(e)$ to the subbundle
$\beta^*\ker(\restricted{e}{M^{\mu}})$
\begin{equation}
\label{eq-blown-up-homomorphism-restricts-to-eta}
B^{\mu}(e\restricted{)}{\ker(e)} : 
\beta^*\ker(\restricted{e}{M^{\mu}}) \rightarrow 
\beta^*\coker(\restricted{e}{M^{\mu}})(-E)
\end{equation}
is equal to the composition of the Petri map with 
the codifferential of $\beta$
\[
\beta^*\ker(\restricted{e}{M^{\mu}}) \RightArrowOf{\beta^*\phi} 
\beta^*\coker(\restricted{e}{M^{\mu}})\otimes 
N^*_{M^{\mu}/M} 
\RightArrowOf{d^*\beta}
\beta^*\coker(\restricted{e}{M^{\mu}})(-E).
\] 

\noindent
We used above the identification
\[
\StructureSheaf{E}(-E) \cong N^*_{E/B^{\mu}M}.
\]
\item
\label{lemma-item-petri-map-vs-eta}
The pullback $\beta^*\phi$ of the Petri map $\phi$ of $e$ 
is related to $\eta$ by the commutative diagram:
\begin{equation}
{
\divide\dgARROWLENGTH by 2
\begin{diagram}
\node{
\beta^*\ker(\restricted{e}{M^{\mu}})\otimes
\beta^*\coker(\restricted{e}{M^{\mu}})^*
}
\arrow{e,t}{\beta^*\phi}
\arrow{se,t}{\eta}
\node{\beta^*\left[ N^*_{M^{\mu}/M}\right]}
\arrow{s,l}{d^*\beta}
\\
\node[2]{
\StructureSheaf{E}(-E).
} 
\end{diagram}
}
\end{equation}
\item
\label{lemma-item-ker-of-blown-up-homomorphism-is-ker-of-eta}
Over $E$ we have a canonical identifications
\begin{eqnarray*}
\coker(B^{\mu}\restricted{e}{E}) & = &
\coker\left[
\eta: \ker(\beta^*\restricted{e}{E}) \rightarrow
\coker(\beta^*\restricted{e}{E})\otimes(-E)
\right],
\\
\ker(B^{\mu}\restricted{e}{E}) & = &
\ker\left[
\eta: \ker(\beta^*\restricted{e}{E}) \rightarrow
\coker(\beta^*\restricted{e}{E})\otimes(-E)
\right].
\end{eqnarray*}
\item
\label{lemma-item-two-stratifications-coincide}
The two stratifications 
(\ref{eq-hom-stratification}) and 
(\ref{eq-intersection-with-proper-transform-stratification})
of $E$ coincide 
\[
\PP\Hom(
\ker(\restricted{e}{M^{\mu}}),
\coker(\restricted{e}{M^{\mu}}))^t = E^t. 
\]
In particular, $(B^{\mu}M)^{\mu}(B^{\mu}e)$ is empty and 
$B^{\mu}(e)$ has rank $\geq \rank(V_1)+1-\mu$ throughout $B^{\mu}M$. 
\end{enumerate}
\end{new-lemma}

\noindent
{\bf Proof:}
\ref{lemma-item-two-tautological-line-bundles-isomorphic}) 
The line-bundle
$N_{E/B^{\mu}M}$ is the tautological line-sub-bundle
of $\beta^*N_{M^{\mu}/M}$, while the line-bundle
$\StructureSheaf{\PP\Hom(
\ker(\restricted{e}{M^{\mu}}),
\coker(\restricted{e}{M^{\mu}})
)}(-1)$ 
is the tautological line-sub-bundle of 
$\beta^*\Hom(\ker(\restricted{e}{M^{\mu}}),
\coker(\restricted{e}{M^{\mu}}))$. 
We assumed that the two vector bundles $N_{M^{\mu}/M}$ and 
$\Hom(
\ker(\restricted{e}{M^{\mu}}),
\coker(\restricted{e}{M^{\mu}})
)$
are canonically isomorphic. Hence, so do the
two tautological line-sub-bundles of their pull-backs. 

\smallskip
\noindent
\ref{lemma-item-composition-is-projection}) is clear.

\smallskip
\noindent
\ref{lemma-item-blown-up-homomorphism-restricts-to-eta})
Notice first that $B^{\mu}V_{\restricted{1}{E}}$ is an extension
\[
0 \rightarrow
\beta^*\coker(\restricted{e}{M^{\mu}})(-E)
\rightarrow 
B^{\mu}V_{\restricted{1}{E}}
\rightarrow 
\beta^*Im(\restricted{e}{M^{\mu}}) 
\rightarrow 0.
\]

\noindent
Clearly, $B^{\mu}(e)$ maps $\beta^*\ker(\restricted{e}{M^{\mu}})$
to $\beta^*\coker(\restricted{e}{M^{\mu}})(-E)$. 
Tracing through the definition of $\phi$ shows the desired equality.

\smallskip
\noindent
\ref{lemma-item-petri-map-vs-eta})
$\eta^*$ is the embedding of $N^*_{E/B^{\mu}M}$ as
a line-sub-bundle of $\beta^*N_{M^{\mu}/M}$ via $d^*\beta$. 

\smallskip
\noindent
\ref{lemma-item-ker-of-blown-up-homomorphism-is-ker-of-eta})
Follows immediately from parts \ref{lemma-item-composition-is-projection},
\ref{lemma-item-blown-up-homomorphism-restricts-to-eta},
and \ref{lemma-item-petri-map-vs-eta}. 

\smallskip
\noindent
\ref{lemma-item-two-stratifications-coincide})
Part \ref{lemma-item-ker-of-blown-up-homomorphism-is-ker-of-eta}
 implies that, 
over $E$, the ranks of $\beta^*e$ and $B^{\mu}(e)$ are related by
\[
\rank(B^{\mu}e) = 
\rank(\beta^*e)+ \rank(\eta). 
\]
Hence the two stratifications coincide set theoretically. 
Moreover, the equality of the cokernels of $\eta$ and
$B^{\mu}(e\restricted{)}{E}$ implies that both stratifications 
are the determinantal stratifications of a locally free presentation
of the same sheaf on $E$. 
\EndProof

\bigskip
Our construction of the elementary transform 
$B^{\mu}e : B^{\mu}V_0 \rightarrow B^{\mu}V_1$ was not symmetric in 
$V_0$ and $V_1$. If we start instead with the dual homomorphism
\[
e^*:V_1^* \rightarrow V_0^*,
\]
its blow-up
\[
B^{\mu}(e^*) : B^{\mu}V_1^* \rightarrow B^{\mu}V_0^*
\]
maps the pullback $B^{\mu}V_1^*:=\beta^*(V_1^*)$
to the elementary transform $B^{\mu}V_0^*$ of $\beta^*(V_0^*)$. 
Combining Lemma 
\ref{lemma-compatibility-of-two-stratifications-on-exceptional-divisor}
part \ref{lemma-item-ker-of-blown-up-homomorphism-is-ker-of-eta}
with its analogue  for $B^{\mu}(e^*)$ restores the symmetry:

\begin{cor}
\label{cor-ker-of-blown-up-e-vs-coker-of-blown-up-dual-e}
\begin{enumerate}
\item
\label{cor-item-domain-of-petri-map-of-Be-equal-that-of-Be-dual}
We have canonical isomorphisms:
\begin{eqnarray*}
\ker\left(
B^{\mu}(e^*\restricted{)}{B^{\mu}M^{t}\setminus B^{\mu}M^{t+1}}
\right)
& \cong &
\coker\left(
B^{\mu}(e\restricted{)}{B^{\mu}M^{t}\setminus B^{\mu}M^{t+1}}
\right)^*(-E),
\\
\coker\left(
B^{\mu}(e^*\restricted{)}{B^{\mu}M^{t}\setminus B^{\mu}M^{t+1}}
\right)
& \cong &
\ker\left(
B^{\mu}(e\restricted{)}{B^{\mu}M^{t}\setminus B^{\mu}M^{t+1}}
\right)^*(-E).
\end{eqnarray*}

\noindent
Moreover, the sheaves $\coker\left((B^{\mu}e^*)^*\right)(-E)$
and $\coker(B^{\mu}(e))$ are isomorphic globally over $B^{\mu}M$.
Similarly, the sheaves $\coker\left((B^{\mu}e^*)\right)(E)$ and
$\coker\left((B^{\mu}e)^*\right)$ are isomorphic globally over $B^{\mu}M$.

\item
\label{cor-item-petri-maps-of-Be-equal-that-of-Be-dual}
The Petri maps of $B^{\mu}(e)$ and $B^{\mu}(e^*)$ are equal

\[
\phi_{e}=\phi_{e^*}:
\ker\left((B^{\mu}e\restricted{)}{B^{\mu}M^t\setminus B^{\mu}M^{t+1}}\right)
\otimes
\coker\left((B^{\mu}e\restricted{)}{B^{\mu}M^t\setminus 
B^{\mu}M^{t+1}}\right)^*
\rightarrow
T^*B^{\mu}\restricted{M}{B^{\mu}M^t\setminus B^{\mu}M^{t+1}}.
\]
\end{enumerate}
\end{cor}

\noindent
{\bf Proof:}
\ref{cor-item-domain-of-petri-map-of-Be-equal-that-of-Be-dual})
We have an injective sheaf homomorphism of short exact sequences:
\[
{
\divide\dgARROWLENGTH by 2
\begin{diagram}
\node{0}
\arrow{e}
\node{(\beta^*V_1)^*}
\arrow{e}
\arrow{s,l}{B(e^*)}
\node{(B^{\mu}V_1)^*}
\arrow{e}
\arrow{s,l}{B(e)^*}
\node{\left[
\beta^*\coker(\restricted{e}{M^{\mu}})\right]^*(E)}
\arrow{e}
\arrow{s,l}{\eta^*}
\node{0}
\\
\node{0}
\arrow{e}
\node{B^{\mu}(V_0^*)}
\arrow{e}
\node{\beta^*V_0^*}
\arrow{e}
\node{\beta^*\ker(\restricted{e}{M^{\mu}})^*}
\arrow{e}
\node{0.}
\end{diagram}
}
\]
We get the quotient short exact sequence
\[
0
\rightarrow
\coker(B(e^*))
\rightarrow
\coker((Be)^*) 
\rightarrow
\coker(\eta^*)
\rightarrow 0.
\]
Lemma 
\ref{lemma-compatibility-of-two-stratifications-on-exceptional-divisor}
part \ref{lemma-item-ker-of-blown-up-homomorphism-is-ker-of-eta}
(applied to $e^*$) 
implies that the surjective homomorphism factors through
an isomorphism
\[
\coker((Be)^*\restricted{)}{E} \IsomRightArrow 
\coker(\eta^*).
\]
Hence, $\coker(B(e^*))$ is isomorphic to
$\coker((Be)^*)\otimes \StructureSheaf{B^{\mu}M}(-E)$ as sheaves
over $B^{\mu}M$. 

\medskip
\noindent
\ref{cor-item-petri-maps-of-Be-equal-that-of-Be-dual})
By Lemma \ref{lemma-existence-of-petri-map}
part \ref{lemma-item-petri-map-of-dual-homo-is-the-same}) it
suffices to prove that the Petri maps of 
$B^{\mu}(e)$ and $B^{\mu}(e^*)^*$ are equal. 
Part \ref{cor-item-domain-of-petri-map-of-Be-equal-that-of-Be-dual} 
implies that the two Petri maps arise from two locally free
presentations of the same sheaf (up to a twist by the line bundle
$\StructureSheaf{B^{\mu}M}(E)$). 
\EndProof

\section{Construction of the dual collection}
\label{sec-construction-of-the-dual-collection}

We have collected the background material on determinantal varieties 
necessary for the proofs of Theorems
\ref{thm-stratified-elementary-transformation} and 
\ref{thm-correspondence-is-a-linear-combination}. 
We prove 
Theorem \ref{thm-correspondence-is-a-linear-combination} is section
\ref{proof-that-correspondence-induces-isomorphism-of-coho-rings}. 

\subsection{Proof of Theorem \ref{thm-stratified-elementary-transformation}}
\label{subsec-proof-of-thm-construction-of-dual-collection}
The main facts needed for the proof are 
the existence of the relative transposition isomorphism 
(Proposition \ref{prop-transposition-for-cotangent-of-grassmannian-bundle}) 
and Lemma's \ref{lemma-inner-square-commutes} and
\ref{lemma-two-tautological-h-subbundles-are-isomorphic}. 
Lemma \ref{lemma-two-tautological-h-subbundles-are-isomorphic} implies
that the normal bundle of the exceptional divisor 
$B^{[k]}Y(r)^k$ (which is to be blown-down
next) restricts as $\StructureSheaf{}(-1)$ to the fibers of the
ruling of $B^{[k]}Y(r)^k$. It follows that we can blow-down 
$B^{[k]}Y(r)$ along $B^{[k]}Y(r)^k$. Lemma 
\ref{lemma-inner-square-commutes} is needed in order to identify the 
affect of this blow-down on the remaining exceptional divisors.

\begin{prop}
\label{prop-recursive-construction-of-blow-down}
We construct, recursively with respect to $k$:
\begin{enumerate}
\item[i.]
Schemes $B^{[k]}Y(r)^t$, parametrized by integers $(r,k,t)$ in the ranges 
$0 \leq r \leq \mu(X)$, 
$1 \leq k \leq \mu(X(r))+1$ and 
$0 \leq t \leq \mu(X(r))$. We denote $B^{[\mu(X(r))+1]}Y(r)^t$ also by 
$Y(r)^t$ and $Y(r)^0$ by $Y(r)$.
\item[ii.]
Morphisms
\[
\tilde{f}'_{r,k-1}:B^{[k]}Y(r)^{k-1} \rightarrow B^{[1]}Y(r+k-1),
\]
and
\item[iii.]
blow-down morphisms
\[
\beta'_{k-1}:B^{[k-1]}Y(r)^t \rightarrow B^{[k]}Y(r)^t, \ \ 
1 \leq k \leq \mu(X(r))+1,
\]
\end{enumerate}

\noindent
satisfying 

\begin{enumerate}
\item \label{prop-item-dual-top-blow-up-is-original-top-blow-up}
$B^{[1]}Y(r)^t := B^{[1]}X(r)^t$, for $0 \leq r \leq \mu(X(r))$ 
and $0 \leq t \leq \mu(X(r))$. We define $\PP{W}_{B^{[1]}Y(j)}$ 
to be the dual $\PP{W}^*_{B^{[1]}X(j)}$ of (\ref{eq-projective-bundle}).
\item 
\label{prop-item-smoothness-of-dual-varieties}
$B^{[k]}Y(r)$ is a smooth variety.
\item
\label{prop-item-def-of-dual-stratification}
$B^{[k]}Y(r)^t$ is the scheme theoretic image of $B^{[k-1]}Y(r)^t$ via
\[
\beta'_{k-1}: B^{[k-1]}Y(r) \rightarrow B^{[k]}Y(r).
\]
\item
\label{prop-item-smooth-exceptional-divisors-in-dual-collection}
If $k \leq t$, then $B^{[k]}Y(r)^t$ is a smooth divisor in $B^{[k]}Y(r)$. 
\item
\label{prop-item-canonical-line-bundle-of-dual-collection}
The canonical line-bundle of $B^{[k]}Y(r)$ is given by
\begin{equation}
\label{eq-canonical-bundle-of-iterated-blow-down}
\omega_{B^{[k]}Y(r)} \cong
\StructureSheaf{B^{[k]}Y(r)}\left(
\sum_{t=k}^{\mu(X(r))}[codim(X(r)^t,X(r))-1]\cdot B^{[k]}Y(r)^t
\right).
\end{equation}
\item
\label{prop-item-blow-down-of-dual-strata}
If $t\geq k$ then we have an isomorphism 
\begin{equation}
\label{eq-stratification-of-iterated-blow-down-t-geq-k}
B^{[k]}Y(r)^t \cong (B^{[k]}\PP\Omega^1G)(t,W_{B^{[1]}Y(r+t)})
\end{equation}
inducing isomorphisms of the exceptional divisors in $B^{[k]}Y(r)^t$
\begin{equation}
\label{eq-exceptional-divisors-on-strata-of-iterated-blow-down-t-geq-k}
B^{[k]}Y(r)^{t\cap j} \cong (B^{[k]}\PP\Omega^1G)(t,W_{B^{[1]}Y(r+t)})^j
\ \ \mbox{for} \ \ j \geq k.
\end{equation}
Moreover, we have a commutative diagram of blow-down morphisms
\begin{equation}
\label{diagram-of-blow-downs}
\begin{array}{ccc}
B^{[k-1]}Y(r) & \LongRightArrowOf{\beta'_{k-1}} & B^{[k]}Y(r)
\\
\cup \ \uparrow \ \hspace{1ex} & & \hspace{1ex}  \ \uparrow \ \cup
\\
B^{[k-1]}Y(r)^t & \LongRightArrowOf{\beta'_{k-1}} & B^{[k]}Y(r)^t
\\
\cong \ \downarrow \ \hspace{1ex} & & \hspace{1ex}  \ \downarrow \ \cong
\\
(B^{[k-1]}\PP\Omega^1G)(t,W_{B^{[1]}Y(r+t)}) & 
\LongRightArrowOf{\beta'_{k-1}} &
(B^{[k]}\PP\Omega^1G)(t,W_{B^{[1]}Y(r+t)})
\\
\downarrow & & \downarrow
\\
B^{[1]}Y(r\!+\!t) & \LongRightArrowOf{=} & B^{[1]}Y(r\!+\!t).
\end{array}
\end{equation}


\noindent
If $t=k-1$, then we have an isomorphism
\begin{equation}
\label{eq-stratification-of-iterated-blow-down-t-eq-k-minus-1}
B^{[k]}Y(r)^{k-1} \cong G(k\!-\!1,W_{B^{[1]}Y(r+k-1)}).
\end{equation}
\item
\label{prop-item-normal-bundle-to-exceptional-divisors-in-dual-collection}
For $t\geq k$, the normal bundle 
$N_{B^{[k]}Y(r)^t/B^{[k]}Y(r)}$ 
is isomorphic to 
\begin{equation}
\label{eq-normal-bundle-to-exceptional-divisors-t-geq-k-in-dual-collection}
\beta'^*\left\{
h_{B^{[t]}Y(r)^t}
\otimes
(\tilde{f}'_{r,t}\circ\beta'_t)^*\StructureSheaf{B^{[1]}X(r+t)}\left(
-\sum_{i=1}^{\mu(X(r+t))}B^{[1]}X(r+t)^i
\right)
\right\}.
\end{equation}
Above, the line bundle 

\begin{equation}
\label{eq-tautological-line-sub-bundle-h}
h_{B^{[t]}Y(r)^t}
\end{equation} 

\noindent
on $B^{[t]}Y(r)^t$ is the 
tautological sub-bundle of 
\[
(\beta'_t)^*\left(
\SheafHom(q_{B^{[t+1]}Y(r)^t},\tau_{B^{[t+1]}Y(r)^t})
\right)
\]
and $\PP\tau_{B^{[t+1]}Y(r)^t}$, $\PP{q}_{B^{[t+1]}Y(r)^t}$ are the 
sub- and quotient bundles of
$(\tilde{f}'_{r,t})^*(\PP{W}_{B^{[1]}Y(r+t)})$. (Here
$B^{[t+1]}Y(r)^t$ involves the abuse of notation introduced in 
Remark \ref{rem-abuse-notation-of-strata-before-defined}).
Note that although only the projectivization of 
${W}_{B^{[1]}Y(r+t)}$ is defined, the vector bundle 
$\SheafHom(q_{B^{[t+1]}Y(r)^t},\tau_{B^{[t+1]}Y(r)^t})$ is well defined.
\end{enumerate}
\end{prop}

\begin{rem}
\label{rem-abuse-notation-of-strata-before-defined}
{\rm
The varieties $B^{[k]}Y(r)^t$, being subvarieties of 
$B^{[k]}Y(r)$, are not defined until $B^{[k]}Y(r)$ is constructed.
Nevertheless, we will  sometimes (ab)use the notation 
$B^{[k]}Y(r)^t$, with $t\geq k-1$, even before $B^{[k]}Y(r)$ is constructed. 
In that case, $B^{[k]}Y(r)^t$ will only refer to the variety
$(B^{[k]}\PP\Omega^1G)(t,W_{B^{[1]}Y(r+t)})$ without claiming yet that
the latter is a subvariety of $B^{[k]}Y(r)$. The isomorphism between 
$(B^{[k]}\PP\Omega^1G)(t,W_{B^{[1]}Y(r+t)})$ and $B^{[k]}Y(r)^t$
is proven once $B^{[k]}Y(r)$ is constructed.
(see (\ref{eq-stratification-of-iterated-blow-down-t-geq-k})
and (\ref{eq-stratification-of-iterated-blow-down-t-eq-k-minus-1})).
}
\end{rem}

Let us first compute the normal-bundles to the exceptional divisors 
in the iterated blow-ups of $X(r)$. 

\begin{new-lemma}
\label{lemma-normal-bundle-of-divisors-in-full-iterated-blow-up}
Assume that $i \leq t$. 
Let $\beta: B^{[i]}X(r)^t \rightarrow B^{[t]}X(r)^t$ be the iterated blow-up. 
It has been identified as the iterated blow-up
\[
\beta:(B^{[i]}\PP\Omega^1G)(t,W_{B^{[1]}X(r+t)}) \rightarrow
(\PP\Omega^1G)(t,W_{B^{[1]}X(r+t)}).
\]
Moreover, $\beta_t:B^{[t]}X(r)^t \rightarrow B^{[t+1]}X(r)^t$ is the
bundle map $(\PP\Omega^1G)(t,W_{B^{[1]}X(r+t)}) \rightarrow
G(t,W_{B^{[1]}X(r+t)})$ and
$\tilde{f}_{r,t}$ is the bundle map
$\tilde{f}_{r,t}: G(t,W_{B^{[1]}X(r+t)}) \rightarrow B^{[1]}X(r+t)$. 
Let 
\[
\beta':(B^{[i]}\PP\Omega^1G)(t,W_{B^{[1]}Y(r+t)}) \rightarrow
(\PP\Omega^1G)(t,W_{B^{[1]}Y(r+t)})
\]
(resp. $\beta'_t$, $\tilde{f'}_{r,t}$) 
be the analogous morphism with respect to the dual vector 
bundle $W_{B^{[1]}Y(r+t)}$. 
Denote by $h_{B^{[t]}X(r)^t}$ the tautological line sub-bundle analogous to
(\ref{eq-tautological-line-sub-bundle-h}). 
Then, we have the following isomorphisms:
\begin{enumerate}
\item
\label{lemma-item-normal-line-bundle-to-first-exceptional-divisor}
$N_{B^{[i]}X(r)^t/B^{[i]}X(r)}$
is isomorphic to
\[
\beta^*\left\{
h_{B^{[t]}X(r)^t}\otimes(\tilde{f}_{r,t}\circ\beta_t)^*
\left(
\StructureSheaf{B^{[1]}X(r+t)}
\left[-\sum_{k=1}^{\mu(X(r+t))}B^{[1]}X(r+t)^k\right]
\right)
\right\}
\]
\item
\label{lemma-item-two-tautological-h-subbundles-are-isomorphic}
Over the top iterated blow-up, 
$N_{B^{[1]}Y(r)^t/B^{[1]}Y(r)}$ (which is $N_{B^{[1]}X(r)^t/B^{[1]}X(r)}$)
is also isomorphic to
\[
\beta'^*\left\{
h_{B^{[t]}Y(r)^t}
\otimes(\tilde{f'}_{r,t}\circ\beta'_t)^*
\left(
\StructureSheaf{B^{[1]}Y(r+t)}
\left[-\sum_{k=1}^{\mu(X(r+t))}B^{[1]}Y(r+t)^k\right]
\right)\right\}.
\]
Note that $h_{B^{[t]}Y(r)^t}$ stands for 
$h_{(\PP\Omega^1G)(t,W^*_{B^{[1]}X(r+t)})}$ 
(see Remark (\ref{rem-abuse-notation-of-strata-before-defined})).
\end{enumerate}
\end{new-lemma}

\noindent
{\bf Proof:} 
\ref{lemma-item-normal-line-bundle-to-first-exceptional-divisor})
The line bundle $N_{B^{[t]}X(r)^t/B^{[t]}X(r)}$
is the tautological line sub-bundle of the vector bundle 
$\beta^*_{t}\left(N_{B^{[t+1]}X(r)^t/B^{[t+1]}X(r)}\right)$ and Lemma
\ref{lemma-normal-bundle-is-twisted-relative-cotangent-bundle}
identifies the latter as the homomorphism bundle
\[
\beta^*_{t}\left[
\Hom(q_{B^{[t+1]}X(r)^t},\tau_{B^{[t+1]}X(r)^t})\left(
-\sum_{k=t+1}^{\mu(X(r))}B^{[t+1]}X(r)^k
\right)\right].
\]
Hence, $N_{B^{[t]}X(r)^t/B^{[t]}X(r)}$ is isomorphic to 
$\left[h_{B^{[t]}X(r)^t}\otimes \left(
-\sum_{k=t+1}^{\mu(X(r))}B^{[t]}X(r)^k
\right)\right].$
The morphism $\tilde{f}_{r,t}$ is compatible with respect to 
the stratifications up to a shift of indices by $t$
(Condition \ref{cond-compatibility-of-stratifications-wrt-abstract-tyurin}). 
The case $i=t$ of the lemma follows.
The case $1\leq i < t$ is an immediate consequence of the case $i=t$
since the normal line-bundle $N_{B^{[i]}X(r)^t/B^{[i]}X(r)}$ is
simply $\StructureSheaf{B^{[i]}X(r)^t}(B^{[i]}X(r)^t)$. 

\medskip
\noindent
\ref{lemma-item-two-tautological-h-subbundles-are-isomorphic})
Follows immediately from 
part \ref{lemma-item-normal-line-bundle-to-first-exceptional-divisor} 
and the isomorphism 
$\beta^*\left(h_{B^{[t]}X(r)^t}\right)
\cong
(\beta')^*\left(h_{B^{[t]}Y(r)^t}\right)$ which is proven in 
Lemma \ref{lemma-two-tautological-h-subbundles-are-isomorphic}.
%
\EndProof

\bigskip
\noindent
{\bf Proof of Proposition \ref{prop-recursive-construction-of-blow-down}:}
The proof is by induction on $k$.

\noindent
\underline{The case $k=1$}:
The morphisms $\tilde{f}'_{r,0}$ and $\beta'_0$ are the identity morphisms.
Properties
\ref{prop-item-dual-top-blow-up-is-original-top-blow-up} and 
\ref{prop-item-def-of-dual-stratification}
are definitions. 
Properties
\ref{prop-item-smoothness-of-dual-varieties}, 
\ref{prop-item-smooth-exceptional-divisors-in-dual-collection}, and 
\ref{prop-item-canonical-line-bundle-of-dual-collection}
are clear since $B^{[1]}Y(r)^t=B^{[1]}X(r)^t$. 
In Property \ref{prop-item-blow-down-of-dual-strata} the isomorphism 
(\ref{eq-stratification-of-iterated-blow-down-t-geq-k}) 
is the transposition isomorphism when $k=1$ since $\PP{W}_{B^{[1]}Y(r+t)}$
is $\PP{W}^*_{B^{[1]}X(r+t)}$ (see Proposition 
\ref{prop-transposition-for-cotangent-of-grassmannian-bundle})
\[
B^{[1]}Y(r)^t= B^{[1]}X(r)^t \cong
(B^{[1]}\PP\Omega^1G)(t,W_{B^{[1]}X(r+t)})
\IsomRightArrowOf{ref}
(B^{[1]}\PP\Omega^1G)(t,W_{B^{[1]}Y(r+t)}). 
\]
The commutativity of Diagram (\ref{diagram-of-blow-downs}) 
is clear since $\beta'_0$ is the identity morphism. 
Property
\ref{prop-item-normal-bundle-to-exceptional-divisors-in-dual-collection}
is verified in Lemma
\ref{lemma-normal-bundle-of-divisors-in-full-iterated-blow-up}. 

\medskip
\noindent
\underline{Induction Step}: Assume that $B^{[j]}Y(r)^t$, $\tilde{f}'_{r,j-1}$,
and $\beta'_{j-1}$ were defined for $j\leq k$ and that they satisfy the
properties of the proposition. 
We need to construct $B^{[k+1]}Y(r)^t$, 
$\tilde{f}'_{r,k}:B^{[k+1]}Y(r)^k\rightarrow B^{[1]}Y(r+k)$, and
$\beta'_k:B^{[k]}Y(r)^t\rightarrow B^{[k+1]}Y(r)^t$ and prove that the
satisfy the properties of the proposition. 
We first construct the  blowing down of $B^{[k]}Y(r)$ along 
the exceptional divisor $B^{[k]}Y(r)^k$ 
(see (\ref{eq-def-of-blow-down-along-k-th-exceptional-divisor}) below). 
By the induction hypothesis 
(property \ref{prop-item-blow-down-of-dual-strata}), 
$B^{[k]}Y(r)^k$ is isomorphic to
$(\PP\Omega^1G)(k,W_{B^{[1]}Y(r+k)})$. Define
\begin{equation}
\label{eq-def-of-contraction-of-k-th-exceptional-divisor}
\beta'_k:B^{[k]}Y(r)^k\rightarrow G(k,W_{B^{[1]}Y(r+k)}) \ \ \mbox{and}
\end{equation}
\begin{equation}
\tilde{f}'_{r,k}:G(k,W_{B^{[1]}Y(r+k)})\rightarrow B^{[1]}Y(r+k)
\end{equation}
as the natural projections. 
(We will prove below that $B^{[k+1]}Y(r)^k$ is isomorphic to 
$G(k,W_{B^{[1]}Y(r+k)})$).
Property 
\ref{prop-item-normal-bundle-to-exceptional-divisors-in-dual-collection}
and the induction hypothesis imply that $N_{B^{[k]}Y(r)^k}$ is isomorphic to 
\[
h_{B^{[k]}Y(r)^k}\otimes(\tilde{f}'_{r,k}\circ\beta'_k)^*
\left(
\StructureSheaf{B^{[1]}Y(r+k)}
\left[
-\sum_{i=1}^{\mu(X(r+k))}B^{[1]}Y(r+k)^i
\right]
\right).
\]
In particular, $N_{B^{[k]}Y(r)^k}$ restricts 
as $\StructureSheaf{}(-1)$
to each fiber of $\beta'_k$
(which is a projective space $\PP^{[codim(X(r)^k,X(r))-1]}$).
Hence, there exists a smooth projective variety $B^{[k+1]}Y(r)$ and a morphism
\begin{equation}
\label{eq-def-of-blow-down-along-k-th-exceptional-divisor}
\beta'_k:B^{[k]}Y(r) \rightarrow B^{[k+1]}Y(r)
\end{equation}
which is an isomorphism away from $B^{[k]}Y(r)^k$ and contracts the fibers of 
(\ref{eq-def-of-contraction-of-k-th-exceptional-divisor})
(see \cite{artin-algebraization,fujiki-nakano}). We get that the image 
$B^{[k+1]}Y(r)^k$ of $B^{[k]}Y(r)^k$ is isomorphic to 
$G(k,W_{B^{[1]}Y(r+k)})$. Properties
\ref{prop-item-dual-top-blow-up-is-original-top-blow-up},
\ref{prop-item-smoothness-of-dual-varieties}, and
\ref{prop-item-def-of-dual-stratification}
follow. 

\medskip
\noindent
\underline{Verification of Property \ref{prop-item-blow-down-of-dual-strata}}:
The induction hypothesis provides the isomorphisms
\begin{eqnarray*}
B^{[k]}Y(r)^t & \cong & (B^{[k]}\PP\Omega^1G)(t,W_{B^{[1]}Y(r+t)}), \ \ t\geq k
\ \ \mbox{and} 
\\
B^{[k]}Y(r)^{k\cap t} & \cong &
(B^{[k]}\PP\Omega^1G)(t,W_{B^{[1]}Y(r+t)})^k, \ \ \mbox{for} \ \ t>k
\end{eqnarray*}
(see (\ref{eq-exceptional-divisors-on-strata-of-iterated-blow-down-t-geq-k})
and (\ref{eq-stratification-of-iterated-blow-down-t-geq-k})).
Lemma \ref{lemma-inner-square-commutes}
implies that the following diagram is commutative

\begin{equation}
\label{diagram-two-descriptions-of-the-intersection-of-exceptional-divisors}
\begin{array}{ccccc}
B^{[k]}Y(r)^k & \LongRightArrowOf{\cong} &
(\PP\Omega^1G)(k,W_{B^{[1]}Y(r+k)}) & \LongRightArrowOf{\beta'_k} & 
G(k,W_{B^{[1]}Y(r+k)})
\\
\cup \ \uparrow \ \hspace{1ex} & &
\cup \ \uparrow \ \hspace{1ex} & &
\cup \ \uparrow \ \hspace{1ex} 
\\
B^{[k]}Y(r)^{k\cap t} & \LongRightArrowOf{\cong} &
(\PP\Omega^1G)(k,W_{\restricted{B^{[1]}Y(r+k)}{B^{[1]}Y(r+k)^{t-k}}}) &
\LongRightArrowOf{\delta_k} & 
G(k,W_{\restricted{B^{[1]}Y(r+k)}{B^{[1]}Y(r+k)^{t-k}}})
\\
= \ \downarrow \ \hspace{1ex} & &
\cong \ \downarrow \ \hspace{1ex} & &
\cong \ \downarrow \ \hspace{1ex}
\\
B^{[k]}Y(r)^{k\cap t} & \LongRightArrowOf{\cong} &
(B^{[k]}\PP\Omega^1G)(t,W_{B^{[1]}Y(r+t)})^k & \LongRightArrowOf{\gamma_k} &
(B^{[k+1]}\PP\Omega^1G)(t,W_{B^{[1]}Y(r+t)})^k
\\
\cap \ \downarrow \ \hspace{1ex} & &
\cap \ \downarrow \ \hspace{1ex} & &
\cap \ \downarrow \ \hspace{1ex}
\\
B^{[k]}Y(r)^t & \LongRightArrowOf{\cong} & 
(B^{[k]}\PP\Omega^1G)(t,W_{B^{[1]}Y(r+t)}) & \LongRightArrowOf{\alpha_k} & 
(B^{[k+1]}\PP\Omega^1G)(t,W_{B^{[1]}Y(r+t)})
\end{array}
\end{equation}


\noindent
(Diagram (\ref{diagram-inner-square-commutes}) with 
$\PP{H}=\PP{W}_{B^{[1]}Y(r+t)}$ is the middle right-hand square in diagram 
(\ref{diagram-two-descriptions-of-the-intersection-of-exceptional-divisors})).
The morphism $\delta_k$ is precisely the restriction of $\beta'_k$
defined above (\ref{eq-def-of-blow-down-along-k-th-exceptional-divisor}).
We see that both $\beta'_k:B^{[k]}Y(r)^t \rightarrow B^{[k+1]}Y(r)^t$ 
and 
\[
\alpha_k: B^{[k]}Y(r)^t \rightarrow 
(B^{[k+1]}\PP\Omega^1G)(t,W_{B^{[1]}Y(r+t)})
\]
(see diagram 
(\ref{diagram-two-descriptions-of-the-intersection-of-exceptional-divisors}))
are isomorphisms away from the exceptional divisor 
$B^{[k]}Y(r)^{k\cap t}$ and both are the blow-down morphisms determined 
by the same ruling of
$B^{[k]}Y(r)^{k\cap t}$, namely, the contraction $\gamma_k$ in diagram 
(\ref{diagram-two-descriptions-of-the-intersection-of-exceptional-divisors}).
Hence, the scheme theoretic image $B^{[k+1]}Y(r)^t$
of $B^{[k]}Y(r)^t$ is isomorphic to 
$(B^{[k+1]}\PP\Omega^1G)(t,W_{B^{[1]}Y(r+t)})$ and
diagram (\ref{diagram-of-blow-downs}) is commutative also for $\beta'_k$.
Note also that 
the right hand column of diagram 
(\ref{diagram-two-descriptions-of-the-intersection-of-exceptional-divisors})
consists of the spaces $B^{[k+1]}Y(r)^t$ (top node), $B^{[k+1]}Y(r)^{k\cap t}$
(middle two), and $B^{k+1]}Y(r)^t$ (bottom node). The isomorphism
(\ref{eq-exceptional-divisors-on-strata-of-iterated-blow-down-t-geq-k})
follows also for the $k+1$ case. 

\medskip
\noindent
\underline{Verification of Property 
\ref{prop-item-smooth-exceptional-divisors-in-dual-collection}}:
Follows immediately from 
(\ref{eq-stratification-of-iterated-blow-down-t-geq-k}). 

\medskip
\noindent
\underline{Verification of Property 
\ref{prop-item-canonical-line-bundle-of-dual-collection}}: 
The morphism $\beta'_k$ in 
(\ref{eq-def-of-blow-down-along-k-th-exceptional-divisor}) is the
blow-up of $B^{[k+1]}Y(r)$ along $B^{[k+1]}Y(r)^k$. 
Hence, the canonical line bundle $\omega_{B^{[k]}Y(r)}$ of 
$B^{[k]}Y(r)$ is
\[
\left[
(\beta'_k)^*(\omega_{B^{[k+1]}Y(r)})
\right]
\left(
+[codim(X(r)^k,X(r))-1]\cdot B^{[k]}Y(r)^k
\right).
\]
By the induction hypothesis, $\omega_{B^{[k]}Y(r)}$ is identified in
(\ref{eq-canonical-bundle-of-iterated-blow-down}). 
The injectivity of the pullback homomorphism
\[
(\beta'_k)^*: \Pic(B^{[k+1]}Y(r)) \hookrightarrow 
\Pic(B^{[k]}Y(r))
\]
implies that $\omega_{B^{[k+1]}Y(r)}$ is identified by 
(\ref{eq-canonical-bundle-of-iterated-blow-down}) replacing $k$ by $k+1$. 


\medskip
\noindent
\underline{Verification of Property 
\ref{prop-item-normal-bundle-to-exceptional-divisors-in-dual-collection}}:
We need to prove that the normal line-bundle $N_{B^{[k]}Y(r)^t}$, 
i.e., $\StructureSheaf{B^{[k]}Y(r)^t}(B^{[k]}Y(r)^t)$, 
is a pullback $(\beta'^{[k]})^*(A)$ of the 
line bundle $A$ on $B^{[t]}Y(r)^t$ specified in 
(\ref{eq-normal-bundle-to-exceptional-divisors-t-geq-k-in-dual-collection}).
Here, $B^{[t]}Y(r)^t$ denotes the space $(\PP\Omega^1G)(t,W_{B^{[1]}Y(r+t)})$ 
(see Remark \ref{rem-abuse-notation-of-strata-before-defined}). 
We know that $\StructureSheaf{B^{[k]}Y(r)^t}(B^{[k]}Y(r)^t)$ pulls back
to $\StructureSheaf{B^{[1]}Y(r)^t}(B^{[1]}Y(r)^t)$.
In step $k=1$ 
it was shown that $\StructureSheaf{B^{[1]}Y(r)^t}(B^{[1]}Y(r)^t)$
is the pullback $(\beta'^{[1]})^*(A)$. Thus, both
$(\beta'^{[k]})^*(A)$ and 
$\StructureSheaf{B^{[k]}Y(r)^t}(B^{[k]}Y(r)^t)$ pull back to the same 
line bundle on $B^{[1]}Y(r)^t$. Property 
\ref{prop-item-normal-bundle-to-exceptional-divisors-in-dual-collection}
follows from the injectivity of the pullback homomorphism
\[
\Pic(B^{[k]}Y(r)^t)\hookrightarrow \Pic(B^{[1]}Y(r)^t).
\]
This completes the proof of Proposition
\ref{prop-recursive-construction-of-blow-down}. 
\EndProof

\bigskip
Property \ref{prop-item-canonical-line-bundle-of-dual-collection} of
Proposition \ref{prop-recursive-construction-of-blow-down}
implies that $Y(r):=B^{[\mu(X(r))+1]}Y(r)$ has a trivial canonical
line-bundle. Hence, $Y(r)$ is symplectic. 
This completes the proof of Theorem 
\ref{thm-stratified-elementary-transformation}.

\subsection{Isomorphism of cohomology rings}
\label{proof-that-correspondence-induces-isomorphism-of-coho-rings}

We prove Theorem
\ref{thm-correspondence-is-a-linear-combination} in this section. 
In the case of a Mukai elementary transformation, the theorem follows easily
from the work of Huybrechts
(\cite{huybrects-iseparable-points-in-moduli} Theorem 3.4). 
The stratified case is analogous. The main ingredients of the argument are:
\begin{enumerate}
\item 
(A stratified version of Huybrechts' trick) 
The stratified elementary transformation 
\[
M \leftarrow B^{[1]}M \cong B^{[1]}W \rightarrow W
\] 
can be extended to 
a suitably chosen one-parameter deformation $\M \rightarrow T$ 
of the hyperkahler variety $M$. It has the affect of  replacing the special 
fiber $M$ by its dual $W$. 
\item
Proposition \ref{prop-extended-transposition} is used to relate 
the fiber products of dual Grassmannian bundles to
exceptional divisors in the blown-up family $B^{[1]}\M$. 
\end{enumerate}

Choose a smooth family $\M \rightarrow T$ of projective symplectic varieties
over a smooth irreducible one-dimensional base $T$ whose fiber over a point 
$0\in T$ is $M$. The extension class $\epsilon\in H^1(M,TM)$ of 
\[
0 \rightarrow TM \rightarrow \restricted{T\M}{M} \rightarrow 
\StructureSheaf{M} \rightarrow 0
\]
is mapped via the symplectic structure to a $(1,1)$-class $\epsilon'$ 
in $H^1(M,T^*M)$. Assume that the restriction of $\epsilon'$ to
a Grassmannian-fiber in every stratum 
$[M^t\setminus M^{t+1}] \rightarrow [M(t)\setminus M(t)^1]$ does not vanish. 
The assumption on the class $\epsilon$ implies 
that the normal bundle $N_{M^t}$ of $[M^t\setminus M^{t+1}]$ in
$\M$ restricts to every Grassmannian-fiber  $G(t,n\!+\!2t\!-\!1)$
as the non-trivial extension 
(\ref{eq-idempotent-extension-of-cotangent-bundle}). 
(Compare with the deformation considered in section 
\ref{sec-dual-springer-resolutions}). 
In the complex analytic category, any twistor deformation provides such a
family because $\epsilon'$ is a Kahler class. 
Such an algebraic deformation exists by 
the work of Beauville \cite{beauville-varieties-with-zero-c-1} and
the transversality of the following two hyperplanes in the complex moduli
of deformations of $M$: 
1) the kernel of the homomorphism 
$H^1(M,TM) \rightarrow H^{1,1}(M) \rightarrow H^{1,1}(G(t,n\!+\!2t\!-\!1))$
and 
2) the tangent space to the algebraic deformations of the polarized 
projective variety $(M,\StructureSheaf{M}(1))$. 

We regard the stratification of $M$ also as a stratification of $\M$ 
\[
\M \supset M \supset M^1 \supset \cdots \supset M^\mu. 
\]
The analogue of Theorem 
\ref{thm-stratified-elementary-transformation} holds. In other words, 
the top iterated blow-up $B^{[1]}\M$ admits a dual sequence of blow-downs
resulting in a family $\W$ in which the special fiber $W$ is the dual 
of $M$. The proof is identical once we  replace the
transposition isomorphism in Proposition
\ref{prop-transposition-for-cotangent-of-grassmannian-bundle}
by the extended transposition isomorphism in Proposition
\ref{prop-extended-transposition}. 
Denote by 

\smallskip
$\tilde{f}:B^{[1]}\M \rightarrow B^{[1]}\W$
the natural isomorphism, 

$\Gamma(\tilde{f})$ its graph in 
$B^{[1]}\M \times_T B^{[1]}\W$, 

$\Gamma(f)$ the closure of the graph of the 
birational isomorphism in $\M\times_T \W$, 


$D^t$ the image of the exceptional divisor 
$B^{[1]}\M^t$ under the isomorphism  $B^{[1]}\M\rightarrow\Gamma(\tilde{f})$, 

$D^0$ the graph 
of the isomorphism $\tilde{f}_0:B^{[1]}M\IsomRightArrow B^{[1]}W$.  

\smallskip
\noindent
Proposition \ref{prop-extended-transposition} 
identifies 
$D^t$ as the top-iterated blow-up of the fiber-product of dual 
grassmannian bundles over $B^{[1]}M(t)$
\[
D^t \ \cong \ 
B^{[1]}[G(t,W_{B^{[1]}M(t)})\times_{B^{[1]}M(t)}G(t,W^*_{B^{[1]}M(t)})]. 
\]
Above we used the notation of Condition 
\ref{cond-compatibility-of-grassmanian-bundles} 
in Section \ref{sec-stratified-transformations}. 
$D^0$ is embedded in the fiber of $\Gamma(\tilde{f})$ over $0\in T$. 
The fiber of $\Gamma(\tilde{f})$ over $0\in T$ is reduced and its 
irreducible components are $\{D^t\}_{t=0}^\mu$ because 
$\Gamma(\tilde{f})$ is isomorphic to $B^{[1]}\M$. 
Since $T$ is one-dimensional, 
both $\Gamma(\tilde{f})$ and $\Gamma(f)$
are irreducible varieties flat over $T$. 
We have the blow-down morphisms
\begin{eqnarray}
\beta\times\beta \ \ : \ \  B^{[1]}\M \times_T B^{[1]}\W
& \longrightarrow & \M\times_T \W \ \ \mbox{and} 
\label{eq-beta-times-beta-on-the-whole}
\\
\label{eq-beta-times-beta-on-graphs}
\beta\times\beta \ \ : \ \  \hspace{6ex} \Gamma(\tilde{f}) \hspace{6ex} 
& \longrightarrow & 
\Gamma(f). 
\end{eqnarray}
The set-theoretic image of $\Gamma(\tilde{f})$ via 
(\ref{eq-beta-times-beta-on-the-whole}) is $\Gamma(f)$. 
Proposition \ref{prop-extended-transposition} 
implies that, set theoretically, (\ref{eq-beta-times-beta-on-graphs}) 
is bijective away from codimension-two in $\Gamma(f)$. The morphism 
$\Gamma(\tilde{f})\rightarrow T$ factors through 
$\Gamma(f)\rightarrow T$. Hence, in order to prove that the differential of 
$\Gamma(\tilde{f})\rightarrow \M\times_T\W$ is injective away from 
a codimension two locus in $\Gamma(f)$, it suffices to prove that the 
differential of $D^t\rightarrow \M\times_T\W$ is injective away from 
a codimension one locus in $[M\times W]\cap \Gamma(f)$. This follows again from
Proposition \ref{prop-extended-transposition}. 
In the notation of Theorem \ref{thm-correspondence-is-a-linear-combination}, 
we conclude that (\ref{eq-beta-times-beta-on-graphs}) 
is an isomorphism away from codimension-two in $\Gamma(f)$ and
it maps the divisor $D^t$ birationally onto $\Delta_t$.  
It follows that $\sum_{t=0}^\mu\Delta_t$ represents 
the class in the middle dimension Chow group of $M\times W$
of the fiber of $\Gamma(f)$ over $0\in T$.  This completes the proof of 
Theorem \ref{thm-correspondence-is-a-linear-combination}.

\section{
Brill-Noether duality for moduli spaces of sheaves} 
\label{sec-brill-noether-for-k3}
In section \ref{sec-admissible-collections} we restate more precisely Theorem
\ref{thm-lift-of-symmetry-to-elementary-transformations}, 
that $\sigma$ and $\tau$ lift to stratified elementary transformations
(Theorem \ref{thm-mukai-reflection-extends-to-a-stratified-transformation}).
Section \ref{subsec-self-dual-moduli-spaces} describes the automorphisms of 
cohomology rings which arise when Theorems 
\ref{thm-lift-of-symmetry-to-elementary-transformations}
and \ref{thm-correspondence-is-a-linear-combination}
are applied to the self-dual moduli spaces. In sections 
\ref{sec-stability-criteria}, \ref{coherent-systems} and 
\ref{subsec-constructions-of-determinantal-loci}
we collect facts, mostly known, about stability of sheaves,
Le Potier's moduli spaces of coherent systems, and Brill-Noether loci.
In section \ref{subsec-tyurin-extension-morphism}
we prove that two resolutions  of a Brill-Noether stratum $\M_S(v)^t$ 
are isomorphic (Theorem \ref{thm-the-tyurin-extension-isomorphism}).
When applied to the whole moduli, Theorem 
\ref{thm-the-tyurin-extension-isomorphism} implies that two descriptions of 
the closure of the graph of the birational isomorphism corresponding to
$\tau$ are indeed isomorphic. The analogous result for $\sigma$ is
carried out in section \ref{subsec-lazarsfeld-reflection-isomorphism}. 
In section \ref{subsec-admissibility}
we prove that the collection of moduli spaces, endowed with their Brill-Noether
stratifications, is a dualizable collection
in the sense of Definition \ref{def-admisible-stratified-collection}. In
particular, Theorem \ref{thm-correspondence-is-a-linear-combination}
applies once the dual collection is identified. The rigorous 
identification of the dual collection is proven in section 
\ref{sec-dual-collections}. 

\subsection{Dualizable collections} 
\label{sec-admissible-collections}

Let $S$ be a K3 surface, $H$ an ample line bundle on $S$. 
Given a coherent sheaf $F$, its Hilbert polynomial is defined by
\[
P_F(n) := \chi(F\otimes H^n) := h^0(F(n))-h^1(F(n))+h^2(F(n)).
\]
The Hilbert polynomial of a coherent sheaf $F$ on $S$ of rank $r\geq 0$ 
and pure $d$-dimensional support is:
\[
P_F(n) = \frac{l_0}{d!}n^d + \frac{l_1}{(d-1)!}n^{d-1}+ \dots + l_d(F) :=
\left(\frac{r}{2}H^2\right) n^2 + \left( H\cdot c_1(F)\right) n + 
\frac{1}{2}\left( c_1(F)^2-2c_2(F)\right) + 2r.
\]
Note that if $r>0$ then
\[
l_0(F) := rH^2 \ \ \mbox{and} \ \ 
l_1(F) := H\cdot c_1(F),  
\]
while if $r=0$ and $d=1$ then
\[
l_0(F) :=H\cdot c_1(F) \ \ \mbox{and} \ \ 
l_1(F) := \frac{1}{2}\left( c_1(F)^2-2c_2(F)\right) + 2r.
\]

\noindent
If $p$ and $q$ are two polynomials with real coefficients,
we say that $p \succ q$ (resp. $p \succeq q$) if
$p(n) > q(n)$ (resp. $p(n) \geq q(n)$) for all $n$ sufficiently large.

\begin{defi}
\label{def-stability}
\begin{enumerate}
\item
\label{def-item-gieseker-stability}
A coherent sheaf $F$ on $S$ is called {\em $H$-semi-stable} (resp. 
{\em $H$-stable}) if it has support of pure dimension $d$ 
and any non-trivial subsheaf $F' \subset F$, $F' \neq (0)$, $F' \neq F$ 
satisfies 
\[
\frac{P_{F'}}{l_0(F')} \preceq \frac{P_F}{l_0(F)} \ \ \ (\mbox{resp.} \prec).
\]
\item
\label{def-item-slope-stability}
A coherent sheaf $F$ on $S$ is called {\em $H$-slope-semi-stable} 
if it has support of pure dimension $d\geq 1$ 
and for any non-trivial subsheaf $F'$ we have
\[
\frac{l_1(F')}{l_0(F')} 
\leq
\frac{l_1(F)}{l_0(F)}. 
\]
If $d=1$, $H$-slope-stability is defined using above a strict inequality.
If $d=2$, $H$-slope-stability is defined using above a strict inequality
and considering only non-trivial subsheaves $F'$ of lower rank. 
\end{enumerate}
\end{defi}

\noindent
Observe that $H$-slope-stability implies $H$-stability and 
$H$-semi-stability implies $H$-slope-semi-stability. 

\bigskip
Denote by $\M_S(v)$ the moduli space of $H$-semistable
sheaves on $S$ with support of pure dimension $\geq 1$   and 
Mukai vector $v$. (If $r:=\rank(v)$ is larger than $0$,
we consider only torsion free sheaves. If $r=0$ we consider
sheaves with pure $1$-dimensional support). 
$\M_S(v)$ is a projective scheme
\cite{gieseker,simpson}. If furthermore, all sheaves 
parametrized by $\M_S(v)$ are $H$-stable then, if non-empty,
$\M_S(v)$ is smooth of dimension
\[
d(v):= 2 + \langle v,v \rangle
\]
\cite{mukai-symplectic-structure}.

Fix a line bundle $\LB$ on $S$. Assume 
\begin{condition}
\label{cond-linear-system}
\begin{enumerate}
\item 
\label{cond-part-minimality}
$\LB$ is an effective cartier divisor with {\em minimal}
degree
\[
0 < c_1(\LB)\cdot H = 
\min\{H\cdot C' \mid \ \ C' \ \ \mbox{is \ an \ effective \ divisor}\}.
\]
In particular, all curves in the linear system
$\linsys{\LB}$ are reduced and irreducible. 
\item 
The base locus of 
$\linsys{\LB}$ is either empty or zero-dimensional. 
\item \label{cond-part-generic-curve-is-smooth}
The generic curve in $\linsys{\LB}$ is smooth. 
\item
$H^1(S,\LB)=0$.
\end{enumerate}
\end{condition}

\noindent
If, for example, $\Pic(S)$ is $\Integers\cdot \StructureSheaf{S}(1)$
and $\StructureSheaf{S}(1)$ is very-ample, we can take $\LB$
to be $\StructureSheaf{S}(1)$. Note that ampleness of $\LB$ is not assumed. 
If, for example, $\pi:S\rightarrow \PP^1$ is an elliptic K3 
with Picard number 2 and without multiple fibers, 
then we can take $\LB$ to be $\pi^*\StructureSheaf{\PP^1}(1)$.  

The arithmetic genus of a curve in $\linsys{\LB}$ is  
$g := \frac{1}{2}[c_1(\LB)^2+2]$. 
The Hilbert scheme $S^{[d]}$ 
is naturally isomorphic to the moduli space $\M_S(1,\LB,g\!-\!d)$
\[
S^{[d]} \cong 
\M_S(1,\LB,g\!-\!d). 
\]
Simply associate to a subscheme $D$ the $\LB$-twist of its ideal sheaf
$\LB(-D):=\Ideal{S,D}\otimes\LB$. 

\medskip
Denote by $\C \subset S\times \linsys{\LB}$ the universal curve. 
The compactified relative Picard $J^d_{\C}\rightarrow \linsys{\LB}$ 
of degree $d$ is the moduli space 
\[
J^d_{\C} \cong 
\M_S(0,\LB,d\!+\!1\!-\!g). 
\]

\noindent
Observe that the Mukai vector $(1,\LB,0)$ of $S^{[g]}$ is the reflection
of the Mukai vector $(0,\LB,1)$ of $J^g_{\C}$.

Let $X=S^{[g]}$ and $X(r) := 
\M_S(r,\LB,r\!-\!1)$. 
Consider the Brill-Noether stratification
\[
\M_S(r,\LB,r\!-\!1)^t := \{F\mid \ \ h^1(F) \geq t\}.
\]
Our goal is to prove that the collection $\{X(r)^t \}_{r=1}^{\mu}$ 
is dualizable with $n(X)=2$ and 
\[
\mu := \mu(X):= \max\{r\mid \ \ d(r,\LB,r\!-\!1)\geq 0 \} =
\max\{r\mid \ \ r(r\!-\!1) \leq g\}.
\]
Its dual collection is $\{Y(r)^t \}_{r=1}^{\mu}$ where
$Y:= J^g_{\C}$, $Y(r):=\M_S(r\!-\!1,\LB,r)$, and 
\[
\M_S(r\!-\!1,\LB,r)^t := \{F\mid \ \ h^1(F) \geq t\}.
\]

More generally we consider the planar hyperbola $\Hy$ in Figure 
\ref{eq-graph-of-hyperbola}
whose lattice points represent Mukai vectors $v$ with $c_1(v) = \LB$. 


\begin{thm}
\label{thm-mukai-reflection-extends-to-a-stratified-transformation}
Let $v$ be a Mukai vector with $c_1(v)=\LB$. Then the collection
(\ref{eq-diagram-of-M-v}) associated to $\M_S(v)$ is a 
dualizable collection (Definition \ref{def-admisible-stratified-collection})
with $\mu(\M_S(v))$ given by (\ref{eq-length-of-bn-stratification}) and 
$n(\M_S(v))=\Abs{\chi(v)}+1$. Its dual collection is the one 
associated to the Mukai vector $\sigma(v)$ (or, equivalently, $\tau(v)$)
where $\sigma$ and $\tau$ are the reflections defined in 
Figure \ref{eq-graph-of-hyperbola}.
\end{thm}

The theorem is proven in sections \ref{subsec-admissibility} ,
\ref{subsec-lazarsfeld-reflection-isomorphism}, and \ref{sec-dual-collections}.
In section \ref{subsec-admissibility} 
Corollary \ref{cor-the-collection-is-admissible}
we prove that the collection (\ref{eq-diagram-of-M-v}) associated to
$\M_S(v)$ is dualizable. 
In section \ref{sec-dual-collections} we prove that the collections 
associated to $\M_S(v)$ and $\M_S(\sigma(v))$ 
(or equivalently $\M_S(\tau(v))$) are indeed dual 
(Proposition \ref{prop-dualizable-collections-are-dual}). 

\begin{example}
{\rm
Take $v=(1,\LB,b)$ in Theorem 
\ref{thm-mukai-reflection-extends-to-a-stratified-transformation}.
We get that $S^{[g\!-\!b]}$ and $\M_S(b,\LB,1)$ are related by a
stratified Mukai elementary transformation. 
In particular, when $b=0$ we get that $S^{[g]}$ and $J^g_{\C}$ are
related by a stratified Mukai elementary transformation. 
If the genus is in the range $2\leq g \leq 5$, then 
$\mu=1$ and the elementary transformation is Mukai's. 
If $6 \leq g\leq 11$, then $\mu=2$. When the genus is $6$  
we get a Lagrangian $G(2,5)$ in both $S^{[6]}$ and  $J^6_{\C}$. 
The degree of the composition morphism 
$G(2,5)\hookrightarrow J^6_{\C} 
\rightarrow \linsys{\StructureSheaf{S}(1)} \cong \PP^6$ is equal to $5$,
the cardinality of $W^2_6$ on a generic curve. More generally, whenever 
$4g\!+\!1$ is a perfect square, 
$\mu$ is equal to $\frac{-1\!+\!\sqrt{4g\!+\!1}}{2}$ and we have a Lagrangian 
$G(\mu,2\mu\!+\!1)$ in $J^g_{\C}$ mapping to $\PP^g$ via a finite  
morphism whose degree is the cardinality 
${\displaystyle g!\prod_{i=0}^{\mu}\frac{i!}{(\mu+i)!} }$
of $W^\mu_g$ on a generic curve
(use Castelnuovo's formula, Theorem (1.3) of \cite{a-c-g-h}). 
The dual Lagrangian grassmannian $G(\mu,2\mu\!+\!1)$ in $S^{[g]}$ 
parametrizes length $g$ subschemes spanning a $\PP^{g-1-\mu}$ in $\PP^g$. 
}
\end{example}

\begin{example}
\label{example-g-1}
{\rm
Consider the positive rays on the $r$ and $s$ axis in the hyperbola $\Hy$ 
in Figure \ref{eq-graph-of-hyperbola}.  Theorem 
\ref{thm-mukai-reflection-extends-to-a-stratified-transformation} implies 
that $J^{b+g-1}_{\C}$ and $\M_S(b,\LB,0)$, $b\geq 0$, 
are related by a stratified Mukai elementary transformation. 
In particular, when $b=0$ we get that the relative Brill-Noether loci 
$\{\W^r_{g-1}\mid \ -1\leq r \leq \mu-1\}$ constitute a dualizable collection. 
In this case the collection is self-dual (see section
\ref{subsec-self-dual-moduli-spaces}). 
}
\end{example}

In the course of proving Theorem 
\ref{thm-mukai-reflection-extends-to-a-stratified-transformation}
we will consider also the moduli spaces of coherent systems 
$G^0(\chi(v),\M_S(v))$ and their analogue 
$G_1(\chi(v),\M_S(v))$. They are introduced in Section \ref{coherent-systems}. 
It is instructive to compare the moduli  spaces involved to
their analogues in the example of dual Springer resolutions of
section \ref{sec-dual-springer-resolutions}:
\[
\begin{array}{lr}
\begin{array}{ccc}
B^{[1]}\M_S(v) & \cong & B^{[1]}\M_S(\sigma(v))
\\
\downarrow & & \downarrow 
\\
G^0(\chi(v),\M_S(v)) & \cong & G^0(\chi(v),\M_S(\sigma(v)))
\\
\swarrow &  & \searrow 
\\
\M_S(v) \hfill &  & \hfill\M_S(\sigma(v))
\\
\searrow &  & \swarrow
\\
& \bar{\M}
\end{array}
&
\begin{array}{ccc}
B^{[1]}T^*G(t,H) & \cong & B^{[1]}T^*G(t,H^*)
\\
\downarrow & & \downarrow 
\\
\Hom(q_{h-t}(H),\tau_t(H)) & \cong & \Hom(q_{h-t}(H^*),\tau_t(H^*))
\\
\swarrow &  & \searrow 
\\
T^*G(t,H) \hfill & & \hfill T^*G(t,H^*)
\\
\searrow &  & \swarrow
\\
& \overline{{\cal N}^t}
\end{array}
\end{array}
\]

\noindent
If $\tau(v)$ is considered and $\chi(v)\geq 0$, replace 
$G^0(\chi(v),\M_S(\sigma(v)))$ by $G_1(\chi(v),\M_S(\tau(v)))$. 
Above, $\Hom(q_{h-t}(H),\tau_t(H))$ and $\Hom(q_{h-t}(H^*),\tau_t(H^*))$ 
are the natural vector bundles over $Flag(t,h\!-\!t,H)$ and 
$Flag(t,h\!-\!t,H^*)$ respectively. As subvarieties of
the product $T^*G(t,H)\times T^*G(t,H^*)$ both are isomorphic to the closure 
of the graph of the birational isomorphism. 

The existence of the morphism $B^{[1]}\M_S(v)\rightarrow G^0(\chi(v),\M_S(v))$
is proven in Proposition 
\ref{prop-grassmannian-fibrations-of-blown-up-brill-noether-loci}. 
The isomorphism 
$G^0(\chi(v),\M_S(v))  \cong  G^0(\chi(v),\M_S(\sigma(v)))$ is constructed in 
Theorem \ref{thm-lazarsfeld-reflection-isomorphism}. 
The isomorphism $G^0(\chi(v),\M_S(v)) \cong  G_1(\chi(v),\M_S(\tau(v)))$
is constructed in Theorem \ref{thm-the-tyurin-extension-isomorphism}.
Theorem \ref{thm-mukai-reflection-extends-to-a-stratified-transformation}
establishes the isomorphism $B^{[1]}\M_S(v) \cong B^{[1]}\M_S(\sigma(v))$.


Pursuing the analogy with Springer resolutions,
there should exist a singular moduli space $\bar{\M} := \bar{\M}_S([v])$, 
analogous to the closure of a square-zero nilpotent coadjoint orbit.
$\bar{\M}_S([v])$ is 
associated to each $\Integers/2\times \Integers/2$ orbit $[v]$ 
of a Mukai vector $v$ in $\Hy$. 
$\bar{\M}_S([v])$ should parametrize 
equivalence classes of stable sheaves with respect to the following 
equivalence relation: 
If $\chi(F) \geq 0$, given any $t$-dimensional 
subspace $U\subset \Ext^1(F,\StructureSheaf{S})$ corresponding to an extension  
\[
0 \rightarrow U^*\otimes\StructureSheaf{S} \rightarrow E \rightarrow F 
\rightarrow 0
\]
of $F$ by a trivial vector bundle, we identify $F$ with the formal difference 
$E - \StructureSheaf{S}^{\oplus t}$. 
Forgetting the extension data has the affect of contracting the Grassmannian
$G(t,H^0(E))$ in $\M_S(v)$ where $v=v(F)$.  
Similarly, if  $\chi(F) \leq 0$, given any rank $t$ trivial subsheaf 
\[
0 \rightarrow U\otimes\StructureSheaf{S} \rightarrow F \rightarrow Q 
\rightarrow 0,
\]
we identify $F$ with the formal sum $Q + \StructureSheaf{S}^{\oplus t}$, 
forgetting the extension data. For example, 
when the K3 surface $S$ is elliptic and $\LB$ is the class of an elliptic 
fiber, $\bar{J}^0:=\bar{\M}_S(0,\LB,0)$ is a K3 with an ordinary double point. 

It seems plausible that $\bar{\M}_S([v])$ can be constructed in
general by using Geometric Invariant Theory. 
The existence of  $\bar{\M}_S([v])$ would also follow 
from a better understanding of the determinantal line bundles on $\M_S(v)$. 
The restriction homomorphism 
$\Pic_{\M_S(v)} \rightarrow \PP^{\Abs{\chi(v)}+1}$, to a fiber in
the first stratum, has an infinite cyclic kernel. 
We need to know that there exist line bundles in this infinite cyclic kernel 
which are  generated by global sections and give rise to a contraction of the 
Grassmannian-fibrations of the Brill-Noether strata  
\[
\M_S(v) \ \longrightarrow \ \bar{\M}_S([v]). 
\] 
%

Theorem 
\ref{thm-mukai-reflection-extends-to-a-stratified-transformation} and 
the above prediction provide a ``contraction'' of the 
hyperbola $\Hy$ in Figure \ref{eq-graph-of-hyperbola} to its diagonal and the 
two semi-diagonals of moduli spaces with Euler characteristic $\chi(v)$ equal 
to $0$ or $\pm 1$: Let $v=(r,\LB,s)$, assume that 
$r\!+\!s$ is even, and set $v':=v-\frac{1}{2}(\vec{r}+\vec{s})$. Then 
$\bar{\M}_S(v)$ is embedded in $\bar{\M}_S(v')$ as a contracted 
Brill-Noether locus. If $\chi(v)=r\!+\!s$ is odd, set 
$v' = v-\frac{1}{2}(\vec{r}+\vec{s}\pm\vec{1})$ and get an embedding
in a semi-diagonal entry. 

\subsection{Self-dual moduli spaces}
\label{subsec-self-dual-moduli-spaces}

There are two classes of self-dual moduli spaces:
If the Euler characteristic $\chi(v)=r+s$ vanishes, then $\M_S(v)$ is
$\tau$-self-dual. If $r=s$ then $\M_S(v)$ is $\sigma$-self-dual.

If $\M_S(v)$ is $\tau$-self-dual, then the birational isomorphism
$\tau:\M_S(v)\rightarrow \M_S(v)$ is, by definition, the identity
(see Theorem \ref{thm-the-tyurin-extension-isomorphism}). The Grassmannian
fibrations involved all have fibers of type 
$G(a,\ComplexNumbers^{2a})$, $a\geq 1$, which
are canonically isomorphic to their dual $G(a,(\ComplexNumbers^{2a})^*)$. 
Nevertheless, Theorem \ref{thm-correspondence-is-a-linear-combination}
is non-trivial in this case. It states that the cycle
$
\Gamma(id) + \sum_{t=1}^{\mu(v)}\Delta_t
$
in $\M_S(v)\times \M_S(v)$ 
induces an automorphism of order two of the ring $H^*(\M_S(v),\Integers)$.
In particular, we get the identity
\[
\left( \sum_{t=1}^{\mu(v)}\Delta_t\right)^2 \ \ = \ \ 
-2\cdot\left( \sum_{t=1}^{\mu(v)}\Delta_t\right)
\]
analogous to the reflections of a K3 lattice induced 
by a $(-2)$-curve. If, for example, $\mu(v)=1$, then 
the divisor $\Theta := \M_S(v)^1$ is a $\PP^1$-fibration over 
$\M_S(v\!+\!\vec{1})$. The intersection of $\Theta$ with a $\PP^1$-fiber
is $-2$ and  the endomorphism $\Delta_1$ is $-2$ times a projection 
from $H^*(\M_S(v),\Integers)$ onto the image of 
\[
H^*(\M_S(v\!+\!\vec{1}),\Integers) \ \rightarrow \
H^{*+2}(\M_S(v),\Integers). 
\]

\medskip
When $\M_S(v)$ is $\sigma$-self-dual, the birational isomorphism 
$\sigma:\M_S(v)\rightarrow \M_S(v)$ is of order two
(see Theorem \ref{thm-lazarsfeld-reflection-isomorphism}
and Example \ref{example-reflection-among-hilbert-schemes}). 
It is a regular automorphism only in the case of $J^{g-1}=\M_S(0,\LB,0)$
or when $\mu(v)=0$. Regardless of the regularity of
$\sigma$, Theorem \ref{thm-correspondence-is-a-linear-combination}
implies that the cycle 
$
\Gamma(\sigma) + \sum_{t=1}^{\mu(v)}\Delta_t
$
in $\M_S(v)\times \M_S(v)$ 
induces an automorphism of order two of the ring $H^*(\M_S(v),\Integers)$.
In particular, we get the identity
\begin{equation}
\label{eq-relation-involving-graph-of-sigma}
\Gamma(\sigma)^2 + 2\Gamma(\sigma)\circ 
\left(\sum_{t=1}^{\mu(v)}\Delta_t\right)
+ \left(\sum_{t=1}^{\mu(v)}\Delta_t\right)^2 \ \ \equiv \ \ \Gamma(id).
\end{equation}

Consider for example the Hilbert scheme $S^{[3]}$ of a K3 of genus $g=4$.
$S$ is the intersection of a cubic and a quadric $Q$ in $\PP^4$.
Then $v=(1,\LB,1)$, $\mu(v)=1$, and $(S^{[3]})^1$ is a lagrangian 
$\PP^3$ parametrizing collinear subschemes. Clearly, we have a morphism 
from $(S^{[3]})^1$ to the variety $F(Q)$ of lines on $Q$. Since $S$
does not contain a line, $(S^{[3]})^1$ is isomorphic to $F(Q)$.
$(S^{[3]})^1$ is also isomorphic to $\PP{H}^0(E)$ where $E$ is the unique 
stable vector bundle with Mukai vector $(2,\StructureSheaf{S}(1),2)$. 
Theorem \ref{thm-lazarsfeld-reflection-isomorphism}
implies that we have an exact sequence
\[
0 \rightarrow E^* \rightarrow H^0(E)\otimes\StructureSheaf{S}
\rightarrow E \rightarrow 0. 
\]
Hence, $\PP{H}^0(E)$ is canonically isomorphic to its dual $\PP{H}^0(E)^*$. 
The cycle $\Delta_1$ is
the product $\PP{H}^0(E)\times \PP{H}^0(E)^*$ and the endomorphism
\[
\Delta_1 \ : \ H^*(S^{[3]},\Integers) \ \ \rightarrow \ \ 
H^*(S^{[3]},\Integers)
\]
is the projection onto the line $\mbox{span}_{\Integers}\{[\PP^3]\}$ 
sending a class $\alpha$ in $H^6(S^{[3]},\Integers)$ to
$(\alpha\cup [\PP^3])\cdot[\PP^3]$.
Combining (\ref{eq-relation-involving-graph-of-sigma})
with the equality 
\[
\Delta_1([\PP^3]) \ = \ c_3(T^*\PP^3)\cdot [\PP^3] \ = \ -4\cdot [\PP^3],
\]
we get 
that
$\Gamma(\sigma)([\PP^3])$ is $a\cdot[\PP^3]$ where $a$ is $3$ or $5$. 
Indeed, $a=3$. 
The proof of Theorem  \ref{thm-correspondence-is-a-linear-combination}
shows that  $\Gamma(\sigma)$ is isomorphic to the graph of a regular 
automorphism $\tilde{\sigma}$ of $B^{[1]}S^{[3]}$ which restricts to 
$\PP{T}^*\PP{H}^0(E) = Flag(1,3,H^0(E)) 
\subset [\PP{H}^0(E)\times \PP{H}^0(E)^*]$ 
as the involution interchanging the two factors 
(using the above identification of the two factors). 
Using the ``Key formula'' (Proposition 6.7 page 114 in
\cite{fulton}), we get the equality 
\[
\Gamma(\sigma)([\PP^3]) \ = \ 
(\beta_*\tilde{\sigma}_*\beta^*)([\PP^3]) \ = \ 
3[\PP^3]. 
\]

\subsection{Stability criteria} 
\label{sec-stability-criteria}

We will need the following Lemma which, though a slight 
strengthening of the statement 
of \cite{lazarsfeld} Lemma 1.3, is actually proven there. 

\begin{new-lemma}
\label{lazarsfeld-slope-stability-lemma}
\cite{lazarsfeld}
Let $v=(r,\A,s)$ be a Mukai vector with $r\geq 1$. 
Assume that 
$\linsys{\A}$ is not empty and  $\A$ satisfies the 
minimality condition \ref{cond-linear-system} part 
\ref{cond-part-minimality}.
Let $F$ be a torsion-free sheaf with 
$v(F)=v$. Assume further that  $s\geq 1$ or that $h^0(F)\geq r+1$. 
Then
The following are equivalent:
\begin{enumerate}
\item 
\label{lemma-item-stability}
$F$ is $H$-stable.
\item 
\label{lemma-item-slope-stability}
$F$ is $H$-slope-stable.
\item
\label{lemma-item-caracterization-of-slope-stability}
$F$ satisfies all the following conditions
\begin{enumerate}
\item
\label{lemma-item-vanishing-of-h2}
$H^2(F)=(0)$,
\item
\label{lemma-item-generated-by-global-sections}
$F$ is generated by its global sections away from a zero-dimensional 
subscheme.
\end{enumerate}
\item
\label{lemma-item-caracterization-of-slope-stability-weak-version}
Same as \ref{lemma-item-caracterization-of-slope-stability} except that
we replace \ref{lemma-item-generated-by-global-sections} by
\begin{enumerate}
\item[(b')]
$F$ is generated away from a zero-dimensional 
subscheme by any $r+1$ dimensional subspace of $H^0(F)$. 
\end{enumerate}
\end{enumerate}
\end{new-lemma}

\begin{rem}
\label{rem-lazarsfelds-lemma-without-s-geq-1}
{\rm
If we drop both the assumption that $s\geq 1$ and the assumption 
that $h^0(F)\geq r+1$ then 
the following relations hold:
\ref{lemma-item-caracterization-of-slope-stability} $\Rightarrow$ 
\ref{lemma-item-stability}
$\Leftrightarrow$
\ref{lemma-item-slope-stability}.
}
\end{rem}
\noindent
{\bf Proof:} (of Lemma \ref{lazarsfeld-slope-stability-lemma})

\noindent
{\bf \ref{lemma-item-stability}
$\Leftrightarrow$
\ref{lemma-item-slope-stability}:
}
The assumption on the linear system $\linsys{\A}$ implies that
$H$-slope-semi-stability for $F$ is equivalent to 
$H$-slope-stability. Hence $H$-stability is equivalent to 
$H$-slope-stability. 

\noindent
{\bf 
\ref{lemma-item-caracterization-of-slope-stability-weak-version} $\Rightarrow$ 
\ref{lemma-item-caracterization-of-slope-stability}:
} 
If $s\geq 1$ then 
the Euler characteristic of $F$ satisfies 
$\chi_{F}= r+s \geq r+1.$
Hence, $h^0(F) \geq r+1$.

\noindent
{\bf 
\ref{lemma-item-slope-stability} $\Rightarrow$ 
\ref{lemma-item-caracterization-of-slope-stability-weak-version}:
} 
\ref{lemma-item-vanishing-of-h2})
If $F$ is slope-stable then 
$\Hom(F,\StructureSheaf{S})$ vanishes and hence, by Serre's Duality, 
$H^2(F) = 0$. 

\noindent
\ref{lemma-item-generated-by-global-sections})
Let $U \subset H^0(F)$ be a subspace of dimension $\geq r+1$. 
Let $F'$ be the subsheaf generated by the global sections in $U$.
Since $h^0(F')$ is larger than the rank $r'$ of $F'$,
$F'$ can not be the trivial rank $r'$ bundle.
Since $F'$ is generated by its global sections, 
Lemma \ref{lemma-c1-of-sheaf-generated-by-global-sections} implies that 
$c_1(F')$ is represented by an effective (non-zero) divisor $C' \subset S$. 
Condition \ref{cond-linear-system} part 
\ref{cond-part-minimality}
implies that $c_1(\A)\cdot H \leq C'\cdot H$. 
$H$-slope-stability of $F$ implies that $r=r'$. 
If $r=r'$ and $F/F'$ has a one-dimensional support, then
$c_1(F/F')$ is represented by an effective divisor $C''$ and
$c_1(F)= C'+C''$ contradicting the integrality of all curves in 
$\linsys{\A}$. 
We conclude that $F/F'$ has zero-dimensional support and 
\ref{lemma-item-generated-by-global-sections} holds.

\medskip
\noindent
{\bf 
\ref{lemma-item-caracterization-of-slope-stability}
$\Rightarrow$ 
\ref{lemma-item-slope-stability}:
} 
Let $F'\subset F$ be a non-trivial rank $r'$ subsheaf, $0 < r' < r$.
Let $F''_0$ be the quotient $F/F'$ and
$F'':= F''_0/T(F''_0)$ where
$T(F''_0)$ is the torsion subsheaf of $F''_0$.
Since $\Hom(F,\StructureSheaf{S}) \cong H^2(F)^*$ vanishes,
while the projection $F\rightarrow F''$ is surjective, we conclude that
$F''$ is not the trivial vector bundle. Lemma 
\ref{lemma-c1-of-sheaf-generated-by-global-sections}
implies that $c_1(F'')$ is represented by an effective (non-zero)
curve $C''$. $c_1(T(F''_0))$ is represented by an effective (possibly
zero) divisor $D''$, the $1$-dimensional components of the support of
$T(F''_0)$. The minimality condition
\ref{cond-linear-system} part \ref{cond-part-minimality} implies the
inequality
\[
c_1(F''_0)\cdot H = (C''+D'')\cdot H \geq c_1(\A)\cdot H.
\]
We get the inequality
\[
c_1(F')\cdot H = \left( c_1(\A) - c_1(F''_0)\right)\cdot H \leq 0
\]
which implies that $F'$ does {\em not} slope-destabilize $F$
\[
\frac{l_1(F')}{l_0(F')} =
\frac{c_1(F')\cdot H}{r'H^2} \leq 0 <
\frac{c_1(F)\cdot H}{rH^2}=
\frac{l_1(F)}{l_0(F)}.
\]
Hence $F$ is $H$-slope-stable.
\EndProof

\begin{new-lemma}
\label{lemma-tyurins-extension-is-stable}
\begin{enumerate}
\item
\label{lemma-item-equivalence-of-stability-for-non-split-extensions}
Let $U$ be a vector space, 
\begin{equation}
\label{eq-seq-quotient-of-F-by-a-trivial-subbundle}
0 \rightarrow 
U\otimes_{\ComplexNumbers}\StructureSheaf{S} \rightarrow 
F \rightarrow 
Q \rightarrow 0
\end{equation}
a short exact sequence of coherent sheaves on $S$ such that the line-bundle
$\A:=\det(F)=\det(Q)$ satisfies the minimality condition 
\ref{cond-linear-system} part 
\ref{cond-part-minimality}. Then the following are equivalent:
\begin{enumerate}
\item
\label{lemma-item-F-is-stable}
$F$ is $H$-stable.
\item
\label{lemma-item-Q-is-stable}
$Q$ is $H$-stable (with support of pure dimension $1$ or $2$)
and the homomorphism $U^*\rightarrow \Ext^1(Q,\StructureSheaf{S})$
is injective. 
\end{enumerate}
\item
\label{lemma-item-canonical-exact-seq-of-ext}
If stability holds in part 
\ref{lemma-item-equivalence-of-stability-for-non-split-extensions} 
then there are canonical exact sequences
\begin{eqnarray}
0 \rightarrow U \rightarrow H^0(F) \rightarrow H^0(Q) \rightarrow 0, 
\label{eq-zero-cohomology-of-quotient-vs-F}
\\
0 \rightarrow H^1(F) \rightarrow H^1(Q) \rightarrow 
U\otimes H^{0,2}(S) \rightarrow 0
\label{eq-surjective-homomorphism-H1-to-UDual}
\end{eqnarray}
and $h^2(Q)=h^2(F)=0$.
\item
\label{lemma-item-stability-of-quotient-with-support-on-a-curve}
\label{lemma-item-stability-of-quotient-by-trivial-subbubdles}
Let $F$ be an  $H$-stable sheaf of rank $r$, $\det(F)$ as in 
condition \ref{cond-linear-system} part 
\ref{cond-part-minimality}, and $U \subset H^0(F)$ a subspace
of dimension $r'\leq r$. Then the natural sheaf homomorphism 
$i : U\otimes\StructureSheaf{S} \rightarrow F$ is injective 
and the quotient $Q$ is $H$-stable.
\item
\label{lemma-item-characterization-of-trivial-subbubdles}
Let $F'$ be a subsheaf of an $H$-stable sheaf $F$ of rank $0<r'<r$
with $c_1(F)$ as in condition \ref{cond-linear-system} part 
\ref{cond-part-minimality}. Then
$F'$ is generated by its global sections away from a zero-dimensional
subscheme if and only if $F'$ is isomorphic to the trivial rank $r'$ bundle.
\end{enumerate}
\end{new-lemma}

\noindent
{\bf Proof:}

\medskip
\noindent
{\bf \ref{lemma-item-equivalence-of-stability-for-non-split-extensions})}
The proof is by induction on $r=\dim(U)$. 

\noindent
\underline{Case $r=1$:} 
Let $\epsilon\in \Ext^1(Q,\StructureSheaf{S})$ be the extension class of 
\begin{equation}
\label{eq-extension-of-Q-by-O}
0 \rightarrow \StructureSheaf{S} \rightarrow F \rightarrow Q \rightarrow 0. 
\end{equation}

\noindent
\underline{
\ref{lemma-item-Q-is-stable} $\Rightarrow$ \ref{lemma-item-F-is-stable}
}
Assume that $Q$ is $H$-stable and $\epsilon$ does not vanish. 
We need to prove that $F$ is torsion-free and $H$-stable. 
Clearly, the support of every subsheaf of $F$ has dimension at least $1$. 
Assume that $T\subset F$ is a destabilizing subsheaf, 
$H\cdot c_1(T) \geq H\cdot c_1(F)$. 
Denote by $\bar{T}$ its image in $Q$. Clearly  
$T$ can not be a subsheaf of $\StructureSheaf{S}$. 
Hence, $H\cdot c_1(\bar{T}) \geq H\cdot c_1(T) \geq H\cdot c_1(F)$. 
Slope-stability of $Q$ implies that $(Q/\bar{T})$ is supported on a 
zero-dimensional subscheme. Hence $\Ext^1(Q/\bar{T},\StructureSheaf{S})$ 
vanishes and 
$\Ext^1(Q,\StructureSheaf{S}) \hookrightarrow
\Ext^1(\bar{T},\StructureSheaf{S})$ 
is injective.
Let
\[
0\rightarrow \StructureSheaf{S}\rightarrow F' \rightarrow \bar{T} \rightarrow 0
\]
be the pullback extension. Then $F'$ is a subsheaf of $F$ and $F/F'$ 
is isomorphic to $Q/\bar{T}$. 
Since every subsheaf of $F$ has dimension at least $1$, 
$F$ is torsion free (resp. $H$-stable) if and only if 
$F'$ is. Thus, we may assume that $T\rightarrow Q$ is 
surjective and $\bar{T}=Q$.
Since $\epsilon\neq 0$, the kernel $K$ of  $T \rightarrow Q$ 
is a non-zero subsheaf of $\StructureSheaf{S}$. 
The inequality $H\cdot c_1(T) \geq H\cdot c_1(F) = H\cdot c_1(Q)$ implies 
that $\StructureSheaf{S}/K$ is supported on a zero-dimensional subscheme. 
Hence $F$ is torsion free if and only if $T$ is. If $T$ is torsion free
then it is not a destabilizing subsheaf of $E$ since their ranks are equal. 
If the torsion subsheaf $\tau\subset T$ is non-trivial then $\tau$ has pure 
$1$-dimensional support. In that case $Q$ has a pure one-dimensional support. 
Clearly $\tau\cap\StructureSheaf{S}=(0)$ and
$\tau$ embeds as a subsheaf $\bar{\tau}$ of $Q$. Hence the image
of $\epsilon$ in $\Ext^1(\bar{\tau},\StructureSheaf{S})$ vanishes. 
This contradicts the non-vanishing of $\epsilon$ because $Q/\bar{\tau}$
has zero-dimensional support and thus the homomorphism 
$\Ext^1(Q,\StructureSheaf{S})\hookrightarrow
\Ext^1(\bar{\tau},\StructureSheaf{S})$ is injective.

\medskip
\noindent
\underline{
\ref{lemma-item-F-is-stable} $\Rightarrow$ \ref{lemma-item-Q-is-stable}
}
Assume that $F$ is $H$-stable (and hence torsion free of rank $\geq 1$).
Then the extension (\ref{eq-extension-of-Q-by-O}) is non-trivial.
We first prove that $Q$ has support of pure dimension $1$ or $2$.
Suppose $\bar{T}\subset Q$ is a subsheaf with zero-dimensional support.
Let $T\subset F$ be its inverse image in $F$.
Then the extension $0 \rightarrow \StructureSheaf{S} \rightarrow T 
\rightarrow \bar{T} \rightarrow 0$ splits. 
This contradicts the $H$-stability of $F$.  
Suppose now that $\bar{T}\subset Q$ is a subsheaf with one-dimensional support
but $Q$ has a two dimensional support. Then $H\cdot c_1(\bar{T}) > 0$ 
and the inverse image $T$ of
$\bar{T}$ is a destabilizing subsheaf of $F$.
We conclude that $Q$ has support of pure dimension $1$ or $2$.
If $supp(Q)$ is a curve $C$, then $\det(Q) = \StructureSheaf{S}(C)$,
$C$ is an integral curve and $Q$ is a rank $1$ torsion free sheaf on $C$. 
Hence $Q$ is $H$-stable. If $Q$ is torsion free then it must be 
$H$-slope-stable because the inverse image $T$ of any destabilizing subsheaf $\bar{T} \subset Q$ would destabilize $F$. This completes the proof of the
case $r=1$. 

\medskip
\noindent
\underline{Induction step:} Choose a line $U_1\subset U$ and let $\bar{U}$ be 
the quotient. Let $\bar{F}$ be the quotient 
$F/(U_1\otimes\StructureSheaf{S})$. 
We get two exact sequences
\begin{eqnarray}
\label{eq-r-1-extension}
0 \rightarrow \bar{U}\otimes\StructureSheaf{S}
\rightarrow \bar{F} \rightarrow Q\rightarrow 0,
\\
\label{eq-line-extension}
0 \rightarrow U_1\otimes\StructureSheaf{S} \rightarrow F \rightarrow \bar{F}
\rightarrow 0.
\end{eqnarray}
Denote by $\bar{\epsilon}: \bar{U}^*\rightarrow \Ext^1(Q,\StructureSheaf{S})$
the extension class of (\ref{eq-r-1-extension}) and by 
$\epsilon_1: U_1^* \rightarrow \Ext^1(\bar{F},\StructureSheaf{S})$
the class of (\ref{eq-line-extension}). Then 

\smallskip
\noindent
$F$ is $H$-stable $\Longleftrightarrow$ (case $r=1$) 
\\
$\bar{F}$ is $H$-stable and $\epsilon_1\neq 0$ 
$\Longleftrightarrow$ (case $r=\dim(U)-1$) 
\\
$Q$ is $H$-stable, $\bar{\epsilon}$ is injective and $\epsilon_1\neq 0$
$\Longleftrightarrow$  
\\
$Q$ is $H$-stable and $\epsilon : U^* \rightarrow \Ext^1(Q,\StructureSheaf{S})$
is injective. 

\smallskip
\noindent
This completes the proof of part 
\ref{lemma-item-equivalence-of-stability-for-non-split-extensions}.

\medskip
\noindent
{\bf \ref{lemma-item-canonical-exact-seq-of-ext})}
The exact sequences 
(\ref{eq-zero-cohomology-of-quotient-vs-F}) and 
(\ref{eq-surjective-homomorphism-H1-to-UDual}) are
part of the long exact cohomology sequence associated to
(\ref{eq-seq-quotient-of-F-by-a-trivial-subbundle}).

\medskip
\noindent
{\bf \ref{lemma-item-stability-of-quotient-by-trivial-subbubdles})}
Proof by induction on $r'$. The case $r'=1$ is clear (injectivity
is obvious and the stability of $Q$ is
proven in part 
\ref{lemma-item-equivalence-of-stability-for-non-split-extensions}). 

\smallskip
\noindent
\underline{Induction step}: 
Choose a line $U_1\subset U$ and consider the extensions
(\ref{eq-r-1-extension}) and (\ref{eq-line-extension}). 
$\bar{F}$ is stable because $F$ is (use part 
\ref{lemma-item-equivalence-of-stability-for-non-split-extensions}). 
By the induction hypothesis the sheaf homomorphism
$\bar{U}\otimes\StructureSheaf{S} \rightarrow \bar{F}$ is injective. 
Hence $U\otimes\StructureSheaf{S} \rightarrow F$ is also injective. 
The stability of $Q$ follows from part 
\ref{lemma-item-equivalence-of-stability-for-non-split-extensions}. 

\noindent
{\bf \ref{lemma-item-characterization-of-trivial-subbubdles})}
By Lemma \ref{lemma-c1-of-sheaf-generated-by-global-sections}
either $F'$ is the trivial rank $r'$ bundle, or
$c_1(F')$ is represented by an effective divisor $C'$.
$H$-slope-stability of $F$ and the minimality condition
\ref{cond-linear-system} part 
\ref{cond-part-minimality} rules out the latter.
\EndProof

\begin{new-lemma}
\label{lemma-c1-of-sheaf-generated-by-global-sections}
\cite{lazarsfeld}
Let $U$ be a torsion-free sheaf on a smooth projective surface.
If $U$ is generated by its global sections away from a zero-dimensional
subscheme, then
$c_1(U)$ is represented by an effective (or zero) divisor. 
In that case, $c_1(U)=0$ if and only if $U$ is a trivial vector bundle. 
\end{new-lemma}

\subsection{Coherent systems}
\label{coherent-systems}

Consider the coarse moduli space $G^0(t,\M_S(v))$ parametrizing pairs 
$(F,U)$ consisting of an $H$-stable sheaf $F$ with $v(F)=v$ and a 
$t$-dimensional subspace $U$ of $H^0(S,F)$. 
Le Potier constructed this moduli space as a projective scheme coarsely 
representing a functor in \cite{le-potier-coherent} Theorem 4.12. 
Le Potier's  semi-stability condition is more relaxed
and does not imply the stability of $F$. Nevertheless, $G^0(t,\M_S(v))$ 
embeds as a Zariski open subset in the stable locus of Le Potier's moduli. 

\begin{rem}
\label{rem-GM-v-is-a-hilbert-scheme}
Let $t$ be a positive integer and set $v':=v+(t,0,t)$. 
It should be easy to check 
that $G^0(\chi(v'),\M_S(v'))$ is also a union of components 
of the Hilbert scheme parametrizing subvarieties  in $\M_S(v)$ 
isomorphic to the grassmannian $G(\chi(v)\!+\!t,\chi(v)\!+\!2t)$.
\end{rem}

A family of such pairs over a Noetherian scheme $G$ is a pair ($\F$, $\tau$).
$\F$ is a sheaf over $G\times S$, flat over $S$, such that its 
restriction $\restricted{\F}{S_g}$ is stable with Mukai vector $v$ 
over every closed point $g$ in $G$. 
The family of vector spaces of sections is encoded by a rank $t$ 
locally free $\StructureSheaf{G}$-module $\tau$ in the following way:
the dual $\tau^*$ is a quotient of the relative Ext sheaf 
$\RelExt_p^2(\F,\omega_p)$ where 
$p:G\times S\rightarrow G$ is the projection and $\omega_p$ is the relative 
dualizing sheaf. Here we use 1) Serre's Duality 
$\Ext^2_S(F_g,\omega_S) \cong H^0(S,F_g)^*$, and 2) the base change theorem for
relative Ext sheaves (\cite{lange} Theorem 1.4) which implies that the natural 
homomorphism
$\RelExt_p^2(\F,\omega_p\restricted{)}{g} \rightarrow \Ext^2_S(F_g,\omega_S)$ 
is surjective for all closed points $g$ in $G$ (use the vanishing of
$Ext^3_S(F_g,\omega_S)$). 

The constructions in section \ref{subsec-tyurin-extension-morphism}
require the following descriptions of $\tau$ as a subsheaf of $p_*\F$.

\begin{new-lemma}
\label{lemma-families-of-subspaces-of-sections}
Let $\F$ be a flat family over $G$ of stable sheaves with Mukai vector $v$ 
and $\tau$ a locally free $\StructureSheaf{G}$-module. 
The following data are equivalent: 
\begin{enumerate}
\item
\label{lemma-item-quotient-of-ext-2}
A surjective homomorphism $\RelExt_p^2(\F,\omega_p)\rightarrow \tau^*.$
\item
\label{lemma-item-subsheaf-of-push-forward-sheaf}
An injective homomorphism $\iota : \tau \hookrightarrow p_*\F$ satisfying the 
property that its restriction
$\iota_g : \restricted{\tau}{g} \rightarrow (p_*\F\restricted{)}{g}$
is injective for all closed points $g$ in $G$. 
\end{enumerate}
\end{new-lemma}

\noindent
{\bf Proof:}
\underline{
\ref{lemma-item-quotient-of-ext-2} $\Rightarrow$ 
\ref{lemma-item-subsheaf-of-push-forward-sheaf})
}
Choose  a very ample line bundle $\StructureSheaf{S}(n)$ on $S$
with the property that the sheaf $\RelExt_p^i(\F(n),\omega_p)$ vanishes
for $i=0,1$, is locally free for $i=2$, and the sheaf $p_*\F(n)$ is locally 
free. A choice of a section of $\StructureSheaf{S}(n)$ yields 
the long exact sequence of relative Ext sheaves associated to the short exact 
sequence
\[
0 \rightarrow \F \hookrightarrow \F(n) \rightarrow Q \rightarrow 0. 
\]
Since $\RelExt^3_p(Q,\omega_p)$ vanishes, the homomorphism
$\RelExt_p^2(\F(n),\omega_p)\rightarrow \RelExt_p^2(\F,\omega_p)$ 
is surjective. 
Hence, $\tau$ embeds in  $\RelExt_p^2(\F(n),\omega_p)^*$ as a
subsheaf of the image of $\RelExt_p^2(\F,\omega_p)^*$. 
Similarly, the homomorphism $p_*\F \hookrightarrow p_*\F(n)$
is injective. Serre's Duality identifies 
$\RelExt_p^2(\F(n),\omega_p)^*$ with $p_*\F(n)$ and exhibits $\tau$
as a subsheaf of $p_*\F(n)$ which is contained in $p_*\F$. 
The surjectivity of the composition
\[
\RelExt_p^2(\F(n),\omega_p\restricted{)}{g} \rightarrow 
\RelExt_p^2(\F,\omega_p\restricted{)}{g} \rightarrow \restricted{\tau}{g}^*
\]
implies that the composition 
\[
\restricted{\tau}{g}\rightarrow 
(p_*\F\restricted{)}{g} \rightarrow 
(p_*\F(n)\restricted{)}{g}
\] 
is injective for all closed points $g$ in $G$. It follows that 
$\restricted{\tau}{g}\rightarrow (p_*\F\restricted{)}{g}$ is injective.

\noindent
\underline{
\ref{lemma-item-subsheaf-of-push-forward-sheaf} $\Rightarrow$ 
\ref{lemma-item-quotient-of-ext-2})
}
The proof is similar. 
\EndProof

\bigskip
We will study also the  coarse moduli space $G_1(t,\M_S(v))$  of pairs
$(F,U)$ consisting of an $H$-stable sheaf $F$ with $v(F)=v$ and a 
$t$-dimensional subspace $U$ of $\Ext^1(F,\StructureSheaf{S})$. 
The corresponding functor (from Noetherian schemes  to sets) associates to a
Noetherian scheme $G$ the set of equivalence classes of 
pairs $(\F,\tau)$ as in the case of coherent 
systems except that $\tau^*$ is a quotient of the higher direct image sheaf
$R^1_{p_*}(\F\otimes \omega_p)$. Here we use the vanishing of 
$R^2_{p_*}(\F\otimes \omega_p)$ and the
base change theorem for cohomology (\cite{hartshorne} Theorem III.12.11) 
to conclude that the natural homomorphism 
$R^1_{p_*}(\F\otimes \omega_p\restricted{)}{g} \rightarrow 
H^1(S,F_g\otimes\omega_S)$ is surjective for all
closed points $g$ in $G$. An analogue of Lemma 
\ref{lemma-families-of-subspaces-of-sections} translates the data of a 
surjective homomorphism $R^1_{p_*}(\F\otimes \omega_p)\rightarrow \tau^*$ 
to the data of a homomorphism 
$\tau \hookrightarrow \RelExt^1_p(\F,\StructureSheaf{G\times S})$
injective in each fiber. 
The existence of $G_1(t,\M_S(v))$ was proven in \cite{le-potier-coherent}
Theorem 5.6 in the case where $\mbox{rank}(v)=0$, i.e., when the sheaves have 
pure one-dimensional support. For a general Mukai vector with $c_1(v)$ 
satisfying Condition \ref{cond-linear-system} part \ref{cond-part-minimality},
the existence of $G_1(t,\M_S(v))$ as a projective scheme follows from
the proof of Theorem \ref{thm-the-tyurin-extension-isomorphism} 
where the functor is shown to be equivalent to the functor  
represented by $G^0(t,\M_S(v+\vec{t}))$. 

The following Lemma will be needed in the proof of 
Theorem \ref{thm-the-tyurin-extension-isomorphism}.  
\begin{new-lemma}
\label{lemma-flatness}
(\cite{matsumura} Application 2 page 150) \\
Let $G$ be a Noetherian scheme and 
$\tau \RightArrowOf{u} E \rightarrow Q \rightarrow 0$ 
a right exact sequence of coherent sheaves over $G\times S$. 
Assume that $E$ is flat over $G$. Then the following are equivalent:

\begin{enumerate}
\item
$u$ is injective and $Q$ is flat over $G$,
\item
The restriction 
$\restricted{u}{g} : \restricted{\tau}{g} \rightarrow \restricted{E}{g}$ 
of $u$ to the fiber over every closed point $g$ in $G$ is injective.
\end{enumerate}
\end{new-lemma}

\subsubsection{Universal properties}

Fix a Mukai vector $v$ with $c_1(v)=\LB$. Choose a very ample line bundle 
$\StructureSheaf{S}(n)$ on $S$ satisfying the property: 
The higher cohomologies $H^i(S,F)$, $i>0$,  vanish for
every stable sheaf $F$ on $S$ with Mukai vector $v$. 
Fix also a non-zero section $\gamma$ of $\StructureSheaf{S}(n)$. 
We get a functorial injective homomorphism 
$\otimes\gamma : \F \hookrightarrow \F(n)$ for every family of sheaves on $S$. 
If a universal sheaf $\F_v$ exists over $S\times \M_S(v)$, then 
the pushforward $p_*\F_v(n)$ is a vector bundle on $\M_S(v)$ and 
we get a projective bundle $\PP(p_*\F_v(n))$ over $\M_S(v)$. 
The projective bundle exists even if $\F_v$ does not. We abuse notation and 
denote this  universal projective bundle by
$\PP(p_*\F_v(n))$ (it is locally trivial in the \'{e}tale topology). 
Over $G^0(t,\M_S(v))$ we have a tautological 
$\PP^{t\!-\!1}$ bundle which is a subbundle
\[
\PP\tau := \PP\tau_{G^0(t,\M_S(v))} \ \subset \ \pi^*\PP(p_*\F_v(n))
\]
of the pullback of $\PP(p_*\F_v(n))$ to $G^0(t,\M_S(v))$.

\begin{prop}
\label{prop-universal-propery-of-G0}
The moduli space $G^0(k,\M_S(v))$ satisfies the following universal property.
Assume we are given 
\begin{enumerate}
\item
a scheme $T$ of finite type, a Zariski open covering $\{T_\alpha\}$, 
finite surjective 
\'{e}tale morphisms $\tilde{T}_\alpha \rightarrow T_\alpha$, 
a family $\E_\alpha$, flat over $\tilde{T}_\alpha$,  
of stable sheaves on $\tilde{T}_\alpha \times S$ 
with Mukai vector $v$,
\item
a rank $k$ locally free sheaf $W_\alpha$ on $\tilde{T}_\alpha$, and 
\item
a homomorphism of $\StructureSheaf{\tilde{T}_\alpha}$-modules
\[
i_\alpha : W_\alpha \hookrightarrow p_*\E_\alpha
\]
which is injective on each fiber,
\end{enumerate}
satisfying the compatibility condition: the pull-back of 
$(\E_\alpha, i_\alpha)$ to $\tilde{T}_\alpha\times_{T}\tilde{T}_\beta$ 
is equivalent, as a coherent system, to the pull-back of 
$(\E_\beta, i_\beta)$.

\noindent
Then there exists 
a unique pair $(\kappa,\{\delta_\alpha\})$ consisting of a morphism
\[ 
\kappa:T \rightarrow G^0(k,\M_S(v))
\]
and a  collection of isomorphisms
\begin{equation}
\label{eq-sub-bundles-of-global-sections-are-isomorphic}
\delta_\alpha : \PP{W}_{\tilde{T}_\alpha} \IsomRightArrow 
(\kappa^*\PP\tau\restricted{)}{\tilde{T}_\alpha}
\end{equation}
such that the composition of $\kappa$ with the forgetful morphism to 
$\M_S(v)$ is the classifying morphism of $\E$ and 
the following diagrams commute:
\begin{equation}
\label{eq-diagram-isomorphism-of-global-sections-bundles-induced-by-sheaf-iso}
{
\divide\dgARROWLENGTH by 2
\begin{diagram}
\node{\PP(p_{*}\E_\alpha(n))}
\arrow{e,t}{\cong}
\node{\kappa^*\PP(p_*\F_v(n)\restricted{)}{\tilde{T}_\alpha}}
\\
\node{\PP{W}_{\tilde{T}_\alpha}}
\arrow{n,l}{i\otimes\gamma}
\arrow{e,tb}{\delta_\alpha}{\cong}
\node{(\kappa^*\PP\tau\restricted{)}{\tilde{T}_\alpha}}
\arrow{n}
\end{diagram}
}
\end{equation}
\end{prop}

Note that 
the fact that $G^0(k,\M_S(v))$ is a coarse moduli scheme implies additional 
universal properties (see \cite{git} Definition 7.4). 

\medskip
\noindent
{\bf Proof:}
Since $G^0(k,\M_S(v))$ is a coarse moduli space, we get a unique collection 
$\{(\kappa_\alpha,\delta_\alpha)\}$ where 
$\kappa_\alpha:\tilde{T}_\alpha\rightarrow G^0(k,\M_S(v))$ is the
classifying morphism. 
The compatibility condition implies that 
the two compositions 
\begin{eqnarray*}
\tilde{T}_\alpha\times_{T}\tilde{T}_\beta & \rightarrow & 
\tilde{T}_\alpha \RightArrowOf{\kappa_\alpha}  G^0(k,\M_S(v)) \ \ \mbox{and} 
\\
\tilde{T}_\alpha\times_{T}\tilde{T}_\beta & \rightarrow & 
\tilde{T}_\beta \RightArrowOf{\kappa_\beta}  G^0(k,\M_S(v))
\end{eqnarray*}
are equal. This implies that the collection $\{\kappa_\alpha\}$ descends and
patches to a global morphism $\kappa$ (see \cite{milne} Chapter I Theorem 
2.17).
\EndProof


\begin{prop}
\label{prop-universal-propery-of-G1}
The universal property of $G_1(k,\M_S(v))$, analogous to that in Proposition 
\ref{prop-universal-propery-of-G0}, 
holds as well.
\end{prop}

\subsection{Construction of the Brill-Noether loci}
\label{subsec-constructions-of-determinantal-loci}

We review the construction of the Brill-Noether loci 
in order to set up the notation used throughout the rest of the paper. 
Assume for simplicity of notation that there exists a universal sheaf
$\F$ over $\M_S(v)\times S$. We indicate in Remark 
\ref{rem-non-existence-of-universal-sheaf} the changes necessary 
if a universal sheaf exists only locally. We carry out the construction in 
the case $\chi(v) \geq 0$. (The case $\chi(v) < 0$ is identical, but for a 
shift by $\chi(v)$ of the index $t$ in $\M_S(v)^t$). 

Choose a section $\gamma$ of $H^{\otimes n}$ 
for $n$ sufficiently large so that
$H^i(F(n))$ vanishes for $i>0$ and for all $F$ in $\M_S(v)$. 
Let $\Gamma\in \linsys{H^{n}}$ be the zero divisor of $\gamma$. 
Consider the short exact sequence of sheaves over $\M_S(v)\times S$

\begin{equation}
\label{eq-short-exact-seq-of-sheaves-before-taking-direct-image}
0 \rightarrow \F \RightArrowOf{\gamma} \F(n) 
\rightarrow \F(n\restricted{)}{\Gamma\times \M_S(v)}
\rightarrow 0.
\end{equation}

\noindent
Denote by $p: \M_S(v)\times S\rightarrow \M_S(v)$ the projection.
We get the long exact sequence of the right-derived functor
\begin{equation}
\label{eq-exact-seq-of-higher-direct-image}
0 \rightarrow 
p_*\F 
\rightarrow 
p_*\F(n)
\rightarrow 
p_*\left(\F(n\restricted{)}{\Gamma\times M_S(v)}\right)
\rightarrow 
R^1_{p_*}\F 
\rightarrow 0.
\end{equation}
Note that $H$-slope-stability of all the sheaves parametrized by
$\M_S(v)$ implies that $R^2_{p_*}\F$ vanishes. Hence, so does
$R^1_{p_*}\left(\F(n\restricted{)}{\Gamma\times \M_S(v)}\right)$. 
Clearly, 
$R^2_{p_*}\left(\F(n\restricted{)}{\Gamma\times \M_S(v)}\right)$
vanishes as well. Thus 
$V_1:=p_*\left(\F(n\restricted{)}{\Gamma\times \M_S(v)}\right)$
is a locally free sheaf of rank $(rH^2)n^2+(H\cdot\LB)^n+r(1-g_\Gamma)$.
As the genus  $g_\Gamma$ equals $\frac{2+n^2H^2}{2}$ the rank is 
$
\left(\frac{rH^2}{2}\right)n^2+(H\cdot\LB)^n.
$
Let $V_0$ be $p_*\F(n)$ and
\begin{equation}
\label{eq-locally-free-presentation-of-R1}
\rho : V_0 \rightarrow V_1
\end{equation}
the restriction map. $V_0$ is locally free of rank
\[
\rank(V_0)= \left(\frac{rH^2}{2}\right)n^2+(H\cdot\LB)^n + \chi(v)
=\rank(V_1)+ \chi(v).
\]
In our new notation, (\ref{eq-exact-seq-of-higher-direct-image}) becomes
\begin{equation}
\label{eq-exact-seq-is-locally-free-presentation-of-R1}
0 \rightarrow 
p_*\F 
\rightarrow 
V_0 
\RightArrowOf{\rho} 
V_1 
\rightarrow 
R^1_{p_*}\F 
\rightarrow 0.
\end{equation}
Truncating the leftmost sheaf in 
(\ref{eq-exact-seq-is-locally-free-presentation-of-R1})
we get a locally free presentation of $R^1_{p_*}\F$. 
Observe that Corollary \ref{cor-verification-of-codimension-condition}
implies that $\rho$ is generically surjective if $\chi(v)\geq 0$. 

\bigskip
Define the determinantal loci $\M_S(v)^t$ for $t\geq 0$ as the
subscheme of $\M_S(v)$ which is the zero locus of 
$\Wedge{\rank(V_1)+ 1- t}\rho$. 
Note that (\ref{eq-locally-free-presentation-of-R1}) is a
locally-free presentation of $R^1_{p_*}\F$. 
A standard argument shows that the subscheme
structure of $\M_S(v)^t$ is independent of the choices of $\F$, $n$ and
$\gamma$ (being independent of the locally-free presentation, 
it is independent of $n$ and $\gamma$. Any other choice $\F'$  of
a universal sheaf would result in  $R^1_{p_*}\F'$ which is a twist 
of $R^1_{p_*}\F$ by a line bundle on $\M_S(v)$. Hence $\M_S(v)^t$
is independent of $\F$). 
Note that $\M_S(v)^t$ is supported by the set of sheaves $F$
with $h^1(F) \geq t$
\[
\M_S(v)^t = \{ F \mid \ h^1(F) \geq t \}.
\]
Using the language of Fitting ideals (and their rank) one can show that 
the Brill-Noether loci $\M_S(v)^t$ represent a functor 
(see \cite{a-c-g-h} Remark 3.2 page 179). 

For future reference we note that if 
we take the relative $\RelExt_{p}(\bullet,\omega_S)$ of 
(\ref{eq-short-exact-seq-of-sheaves-before-taking-direct-image})
we get the exact sequence of sheaves on $\M_S(v)$
\[
0
\rightarrow
\RelExt_{p}^1(\F,\omega_S)
\rightarrow
\RelExt_{p}^2(\F(n\restricted{)}{\Gamma\times\M_S(v)},\omega_S)
\rightarrow
\RelExt_{p}^2(\F(n),\omega_S)
\rightarrow
\RelExt_{p}^2(\F,\omega_S)
\rightarrow
0.
\]

\noindent
Serre's Duality identifies it as the locally free presentation
of $\RelExt_{p}^2(\F,\omega_S)$
dual to (\ref{eq-exact-seq-is-locally-free-presentation-of-R1})

\begin{equation}
\label{eq-exact-seq-is-locally-free-presentation-of-Ext1}
0
\rightarrow
\RelExt_{p}^1(\F,\omega_S)
\rightarrow
V_1^*
\RightArrowOf{\rho^*}
V_0^*
\rightarrow
\RelExt_{p}^2(\F,\omega_S)
\rightarrow
0.
\end{equation}

\begin{rem}
\label{rem-non-existence-of-universal-sheaf}
{\rm
We do have a universal sheaf if $c_2(v)=g-2$ or $g$ and
$c_1(v)=\LB$ because in this case $s-r=\pm 1$ 
and hence $gcd(r,c_1(v)^2,s)=1$. If $m:=gcd(r,c_1(v)^2,s)>1$ 
then we only have a global quasi-universal sheaf of similitude $m$
(see \cite{mukai-hodge} Appendix 2). Nevertheless, 
in general, we do have a universal sheaf locally 
(in the complex or \'{e}tale topology) over $\M_S(v)$. 
It is possible to carry out the construction of the Brill-Noether loci 
using only the local existence of a universal sheaf. 
We carry out the constructions locally, 
show independence of the choice of the universal sheaf, and
conclude that the constructions glue as a global algebraic object.
In fact, the stability (and hence simplicity) of the sheaves parametrized by 
$\M_S(v)$ imply (as in \cite{mukai-hodge} Appendix 2)
that the vector bundles $\SheafHom(\Wedge{k}V_0,\Wedge{k}V_1)$ 
and their sections $\Wedge{k}\rho$, $k\geq 1$, exist globally and depend 
canonically on $n$ and $\gamma$ even though the universal sheaf $\F$ 
and the vector bundles $V_0$, $V_1$ exist only locally. 
}
\end{rem}

\subsection{Tyurin's extension morphism}
\label{subsec-tyurin-extension-morphism}


The determinantal loci $\M_S(v)^t$, the moduli spaces 
$G^0(t,\M_S(v))$, $G_1(t,\M_S(v))$ and the forgetful morphisms 
\begin{eqnarray*}
G_1(t,\M_S(v)) \rightarrow \M_S(v)^t, \\
G^0(\chi(v)+t,\M_S(v)) \rightarrow \M_S(v)^t
\end{eqnarray*}
were constructed for the relevant Mukai vectors. 
We extend results of Tyurin 
\cite{tyurin-cycles-curves-surfaces} 
and construct an isomorphism
\[
\tilde{f} : G_1(t,\M_S(v)) \IsomRightArrow
G^0(t,\M_S(v'))
\]
where $v'=v+(t,0,t)$.

Let $F$ be a sheaf of rank $r\geq 0$ and
$V \subset \Ext^1(F,\StructureSheaf{S})$ a $t$-dimensional subspace. 
There exists a canonical rank $r+t$ sheaf $f(F,V)$ and 
an extension 
\begin{equation}
\label{eq-extension-defining-the-tyurin-morphism}
0 \rightarrow V^*\otimes \StructureSheaf{S}
\rightarrow f(F,V) \rightarrow
F \rightarrow 0.
\end{equation}
Simply define $f(F,V)$ via the canonical class in
$\Ext_S^1(F,V^*\otimes \StructureSheaf{S})$ using the isomorphism 
\[
\Ext_S^1(F,V^*\otimes \StructureSheaf{S}) \cong 
\Ext_S^1(F,\StructureSheaf{S})\otimes V^* \cong 
\Hom[V,\Ext_S^1(F,\StructureSheaf{S})].
\] 
It is easy to see that the Mukai vector of $f(F,V)$ is 
\[
v(f(F,V)) = v(F)+(t,0,t).
\]

The sheaf $f(F,V)$ is torsion-free and $H$-slope-stable 
if and only if $F$ is $H$-slope-stable 
(Lemma \ref{lemma-tyurins-extension-is-stable}). 
Moreover, we have the exact sequences
\begin{eqnarray}
0 \rightarrow V^* \rightarrow H^0(f(F,V)) 
\rightarrow H^0(F) \rightarrow 0 \ \ \mbox{and} 
\\
0 \rightarrow H^1(f(F,V)) \rightarrow H^1(F) 
\rightarrow V^*\otimes H^{0,2}(S) 
\rightarrow 0.
\label{eq-connecting-hommomorphism-from-h1-to-h2-for-tyurin-extension}
\end{eqnarray}

\begin{thm}
\label{thm-the-tyurin-extension-isomorphism}
Let $v=(r,\LB,s)$ be a Mukai vector with $r\geq 0$ and $t$ a positive integer. 
There is a natural isomorphism 
\begin{equation}
\label{eq-tyurin-extension-isomorphism}
\tilde{f} : G_1(t,\M_S(v)) \IsomRightArrow
G^0(t,\M_S(v'))
\end{equation}
where $v'=v+(t,0,t)$. The isomorphism $\tilde{f}$ is compatible with
the Brill-Noether stratification with a shift of indices by $t$: 
For $k\geq t$, $\tilde{f}$ realizes 
$G_1(t,\M_S(v)^k)\setminus G_1(t,\M_S(v)^{k+1})$ as a 
$G(t,\chi(v')+k-t)$-bundle over 
$\M_S(v')^{k-t}\setminus \M_S(v')^{k+1-t}$.
\end{thm}

\noindent
{\bf Proof:}
We prove the equivalence of the two functors coarsely represented by
$G_1(t,\M_S(v))$ and $G^0(t,\M_S(v'))$. 
Let $G$ be a Noetherian scheme, $\F_{v'}$ a sheaf on $G\times S$, 
flat over $S$, whose restriction to $S_g$ is a stable sheaf with Mukai vector 
$v'$ for all closed points $g$ in $G$. 
Let $\tau$ be a rank $t$ locally free $\StructureSheaf{G}$-subsheaf of 
$p_*\F_{v'}$ satisfying the property that 
$\restricted{\tau}{g}\rightarrow (p_*\F_{v'}\restricted{)}{g}$ 
is injective for all closed points $g$ in $G$. We get a short exact sequence of sheaves on $G\times S$:
\begin{equation}
\label{eq-functorial-quotient-by-trivial-subbundle}
0 
\rightarrow
p^*\tau
\rightarrow 
\F_{v'}
\rightarrow
Q
\rightarrow
0.
\end{equation}
By Lemma \ref{lemma-tyurins-extension-is-stable}
part \ref{lemma-item-stability-of-quotient-by-trivial-subbubdles}, 
the restriction of (\ref{eq-functorial-quotient-by-trivial-subbundle})
to $S_g$ is exact for all closed points $g$ in $G$. 
Lemma \ref{lemma-flatness} 
implies that $Q$ is flat over $G$. 
Slope stability of the bundles $(\F_{v'}\restricted{)}{g}$, $g\in G$,
implies that the $\StructureSheaf{G}$-module 
$p_*\SheafHom(\F_{v'},\StructureSheaf{G\times S})$ vanishes and the following
short exact sequence of $\StructureSheaf{G}$-modules 
\[
0 \rightarrow \tau^* 
\hookrightarrow \RelExt^1_{p_*}(Q,\StructureSheaf{G\times S})
\rightarrow \RelExt^1_{p_*}(\F_{v'},\StructureSheaf{G\times S})
\rightarrow 0
\]
is part of the long exact sequence of relative extension sheaves. 

Conversely, let $\F_v$ be a sheaf on $G\times S$, 
flat over $S$, whose restriction to $S_g$ is a stable sheaf with Mukai vector 
$v$ for all closed points $g$ in $G$. 
Denote by $\omega_p$ the relative dualizing sheaf $p^*\omega_S$. 
Let $\tau$ be a rank $t$ locally free $\StructureSheaf{G}$-module, 
$\epsilon : \tau  \hookrightarrow 
\SheafExt^1_{p_*}(\F_{v},\omega_p)$ 
an injective homomorphism into the
relative extension sheaf  satisfying the property 
that $\restricted{\tau}{g} \rightarrow 
\RelExt^1_{p_*}(\F_{v},\omega_p\restricted{)}{g}$ 
is injective for all closed points $g$ in $G$. 
We get a section $\epsilon$ of 
$H^0(G,\RelExt^1_{p_*}(\F_v,p^*\tau^*\otimes \omega_p))$. 
The Grothendieck spectral sequence
\[
H^p(G,\RelExt^q(\F_v,p^*\tau^*\otimes\omega_p)) 
\Rightarrow \Ext^{p+q}_{G\times S}(\F_v,p^*\tau^*\otimes\omega_p)
\]
gives the exact sequence:
\begin{eqnarray*}
0 & \rightarrow &
H^1(G,p_*\SheafHom(\F_v,p^*\tau^*\otimes\omega_p)) 
\rightarrow
\Ext^{1}_{G\times S}(\F_v,p^*\tau^*\otimes\omega_p)
\rightarrow 
H^0(G,\RelExt^1_p(\F_v,p^*\tau^*\otimes\omega_p))
\\
& \rightarrow &
H^2(G,p_*\SheafHom(\F_v,p^*\tau^*\otimes\omega_p)).
\end{eqnarray*}

\noindent
But the sheaf $p_*\SheafHom(\F_v,p^*\tau^*\otimes\omega_p)$
vanishes because $\F_v$ is a family of torsion free sheaves $F$ 
with vanishing $H^2(F)$. Hence we get an isomorphism 
\[
\Ext^{1}_{G\times S}(\F_v,p^*\tau^*\otimes\omega_p)
\IsomRightArrow 
H^0(G,\RelExt^1_p(\F_v,p^*\tau^*\otimes\omega_p))
\]
and the section $\epsilon$ determines an extension 
of sheaves on $G\times S$:

\begin{equation}
\label{eq-functorial-tautological-extension}
0 \rightarrow
p^*\tau^*\otimes\omega_p \hookrightarrow \E
\rightarrow \F_v
\rightarrow 0.
\end{equation}

\noindent
Since $p^*\tau^*$ and $\F_v$ are flat over $G$, so is $\E$. Lemma
\ref{lemma-tyurins-extension-is-stable} part 
\ref{lemma-item-equivalence-of-stability-for-non-split-extensions}
implies that $\E$ restricts to $S_g$ as a stable sheaf with Mukai vector 
$v'$ for all closed points $g$ in $G$. 

The two constructions are inverse of each other modulo the equivalence 
of objects parametrized. In other words, 
given the data $\tau_0 \hookrightarrow p_*\F_{v'}$
in the first construction and the data 
$\tau_1 \hookrightarrow \SheafExt^1_{p_*}(\F_{v},\omega_p)$ 
in the second construction, there exists a line bundle $K$ on $G$ 
such that tensoring the short exact sequence 
(\ref{eq-functorial-tautological-extension}) by $K$ we get the sequence 
(\ref{eq-functorial-quotient-by-trivial-subbundle}). In particular, 
$\tau_0 \cong (\tau_1^*\otimes p_*\omega_p)\otimes K$
which is a relative version of
Lemma \ref{lemma-tyurins-extension-is-stable} part 
\ref{lemma-item-canonical-exact-seq-of-ext} combined with 
the relative version of Serre's Duality.

\bigskip
Finally we prove that the isomorphism $\tilde{f}$ in
(\ref{eq-tyurin-extension-isomorphism}) is compatible with the Brill-Noether 
stratifications. By definition, these stratifications are the 
pull-back of the Brill-Noether stratifications on $\M_S(v)$ and $\M_S(v')$ 
via the forgetful morphism. 
Given a Noetherian scheme $G$ and data $\F_{v'}$, $\tau_0$ as in 
(\ref{eq-functorial-quotient-by-trivial-subbundle})
and $\F_{v}$, $\tau_1$ as in 
(\ref{eq-functorial-tautological-extension}) we get 
on $G \times S$ 
the short exact sequence
\begin{equation}
\label{eq-universal-tautological-bundle-is-a-subsheaf}
0 \rightarrow
p^*\tau
\rightarrow 
\F_{v'}
\rightarrow
\F_v\otimes p^*K
\rightarrow 
0
\end{equation}

\noindent
for some line bundle $K$ on $G$ 
(where $\tau\cong \tau_0\cong(\tau_1^*\otimes p_*\omega_p)\otimes K$). 
The long exact sequence of higher direct image brakes into two short 
exact sequences on $G$ one of which is
\begin{equation}
\label{eq-quotient-is-tau}
0 \rightarrow
R^1_{p_*}\F_{v'} 
\rightarrow 
(R^1_{p_*}\F_{v})\otimes K
\RightArrowOf{j}
\tau\otimes_{\ComplexNumbers}H^{0,2}(S)
\rightarrow
0.
\end{equation}

\noindent
Since the quotient $\tau\otimes_{\ComplexNumbers}H^{0,2}(S)$
is locally free of rank $t$, the determinantal stratification 
defined by the sheaf $R^1_{p_*}\F_{v'}$ coincides with that 
of $(R^1_{p_*}\F_{v})\otimes K$ with indices shifted by $t$. 
To see that, choose a locally free presentation
\[
A \RightArrowOf{e}
B \RightArrowOf{\eta}
(R^1_{p_*}\F_{v})\otimes K
\rightarrow 0
\]
with $\rank(A) \geq \rank(B)$. 
The kernel $B_1:=\ker(j\circ \eta)$ is locally free because 
$B$ and $\tau$ are
and $(j\circ \eta)$ is surjective. 
Denote by $a$, $b$ and $b_1$ the corresponding ranks. 
The image $e(A)$ is a subsheaf
of $B_1$ and we get a commutative diagram of locally free presentations
(and short exact columns):

\[
{
\divide\dgARROWLENGTH by 2
\begin{diagram}
\node{A}
\arrow{e,t}{e_1}
\arrow{s,l}{=}
\node{B_1}
\arrow{e}
\arrow{s,l}{\cap}
\node{R^1_{p_*}\F_{v'}}
\arrow{e}
\arrow{s,l}{\cap}
\node{0}
\\
\node{A}
\arrow{e,t}{e}
\node{B}
\arrow{e,t}{\eta}
\arrow{s,l}{j\circ\eta}
\node{(R^1_{p_*}\F_{v})\otimes K}
\arrow{e}
\arrow{s,l}{j}
\node{0}
\\
\node[2]{\tau}
\arrow{e,t}{=}
\node{\tau}
\end{diagram}
}
\]

\noindent
The stratifications of $G$ associated to $R^1_{p_*}\F_{v}$ 
and $R^1_{p_*}\F_{v'}$ are defined by 
\begin{eqnarray*}
G_1(v)^r & := & \left(
\Wedge{[(a-b)+r]}e = 0
\right), \ \ \mbox{and} 
\\
G^0(v')^r & := & \left(
\Wedge{[(a-b_1)+r]}e_1 = 0
\right).
\end{eqnarray*}

\noindent
Clearly $(\Wedge{m}e=0)$ and $(\Wedge{m}e_1=0)$ define
the same subscheme for all $m\geq 0$. The equality 
\[
[(a-b_1)+r] = [(a-b)+r] + t
\] 
verifies that the indices are shifted by $t$. 
\EndProof


\begin{cor}
\label{cor-verification-of-codimension-condition}
Let $v=(r,\LB,s)$ be a Mukai vector, 
$t\geq 1$ an integer and $\mu(v)$ as in (\ref{eq-length-of-bn-stratification}).
If $t > \mu(v)$ then $\M_S(v)^t$ is empty. 
If $0 \leq t \leq \mu(v)$ then the 
codimension of $\M_S(v)^t$ in $\M_S(v)$ is $(\Abs{\chi(v)}+t)t$. 
\end{cor}

\noindent
{\bf Proof:}
Assume first that $\chi(v)$ is non-negative. 
The proof is by induction on $d(v)$ where $v$ ranges in the set
$\{v_0\!+\!\vec{t} \ \mid \ t \in \Integers \ \mbox{and} \  
\chi(v_0\!+\!\vec{t})\geq 0\}$
for some fixed vector $v_0$. Notice that $d(v_0\!+\!\vec{t})$ is a 
decreasing function of $t$ on this set. 
If $d(v)$ is negative, $\M_S(v)^t$ is empty for any integer $t\geq 0$.
If $d(v)$ is non-negative and $t=0$ then 
$\M_S(v)^0=\M_S(v)$ is non-empty by Lemma \ref{non-emptiness-lemma}.
If $d(v+(1,0,1))$ is negative but $d(v)$ is non-negative then 
$\mu(v)=0$. Hence, the first step of the induction involves $v$
with $\mu(v)=0$ for which Corollary 
\ref{cor-verification-of-codimension-condition} holds. 
Assume that the statement holds for all $v''$ with $d(v'')< d(v)$. 
Note that $G_1(t,\M_S(v)^t)\setminus G_1(t,\M_S(v)^{t+1})$
is isomorphic to $\M_S(v)^t\setminus \M_S(v)^{t+1}$. 
Set $v' := v\!+\!\vec{t}$ with $t\geq 1$. 
The induction hypothesis implies that either 
$\M_S(v')$ is empty, or $\M_S(v')^{1}$ is a proper 
(possibly empty) subscheme of $\M_S(v')$. 
By definition, if $t\leq \mu(v)$,  then $\mu(v')=\mu(v)-t$.
If $t>\mu(v)$ then the induction hypothesis implies that 
$\M_S(v')$ is empty.
Hence, Theorem \ref{thm-the-tyurin-extension-isomorphism} implies
that $\M_S(v)^t$ is also empty.
If $t\leq\mu(v)$ then the induction hypothesis implies that 
$\M_S(v')$ is not empty and
Theorem  \ref{thm-the-tyurin-extension-isomorphism} implies
that the dimension 
$\dim[G_1(t,\M_S(v)^t)\setminus G_1(t,\M_S(v)^{t+1})]$
is equal to
\[
t(r+s+t)+\dim[\M_S(v')^{0}\setminus \M_S(v')^{1}]=t(r+s+t)+d(v').
\]
Thus, the codimension of $\M_S(v)^t$ in $\M_S(v)$
is 
\[
d(v)-\left(t(r+s+t)+d(v')\right)= t(r+s+t). 
\]

The case of negative Euler characteristic $\chi(v)$ can be reduced to the
non-negative  Euler characteristic case via the isomorphism

\noindent
$[\M_S(v)^t\setminus \M_S(v)^{t\!+\!1}] \ = \  
[G_1(t\!-\!\chi(v),\M_S(v)) \setminus 
G_1(t\!-\!\chi(v),\M_S(v)^{t\!+\!1})] \ \cong$\\
$[G^0(t\!-\!\chi(v),\M_S(v\!+\!\vec{t}\!-\!\vec{\chi(v)}))\setminus 
G^0(t\!-\!\chi(v),\M_S(v\!+\!\vec{t}\!-\!\vec{\chi(v)})^1)].$ 

\noindent
The first equality is a definition. The second isomorphism 
follows from Theorem \ref{thm-the-tyurin-extension-isomorphism}. 
The Mukai vector $v\!+\!\vec{t}\!-\!\vec{\chi(v)}$ has Euler characteristic
$2t\!-\!\chi(v)$ which is non-negative. 
\EndProof

\begin{new-lemma}
\label{non-emptiness-lemma}
Let $v=(r,\LB,s)$ be a Mukai vector.
If $d(v)$ is non-negative then $\M_S(v)$ is non-empty.
\end{new-lemma}

\noindent
{\bf Proof:}
We may assume that $r$ is non-negative (otherwise, replace $v$ by
$(\sigma\circ\tau)(v):=(-r,\LB,-s)$). 
Let $v':=(0,\LB,s-r)$. Note that $\M_S(v')$ is 
the compactified relative Jacobian $J_{\C}^{g-1+s-r}$. 
Condition \ref{cond-linear-system} 
Part \ref{cond-part-generic-curve-is-smooth} assures us the existence
of a smooth curve in $\linsys{\LB}$. 
Use the non-emptiness results
from classical Brill-Noether theory for a smooth curve in
$\linsys{\LB}$ to conclude that 
$G_1(r,\M_S(v'))$ is  not empty if
$d(v)=2g-2rs$ is non-negative.
Then use Theorem
\ref{thm-the-tyurin-extension-isomorphism}
to conclude that $G^0(r,\M_S(v))$ is non-empty.
\EndProof

\subsection{Dualizability} 
\label{subsec-admissibility}

Let $v=(a,\LB,b)$ be a Mukai vector as in Theorem 
\ref{thm-mukai-reflection-extends-to-a-stratified-transformation}.
We prove in this section that, if $\chi(v)\geq 0$, the collection
\begin{equation}
\label{eq-admissible-collection}
\{ \M_S(v+\vec{r})^t\mid \ 0 \leq r \leq \mu(v) \}
\end{equation}
is a dualizable collection 
(Definition \ref{def-admisible-stratified-collection}). 
Similarly, for $\chi(v)\leq 0$ the collection with $\M_S(v\!+\!\vec{r})^t$
replaced by $\M_S(v\!-\!\vec{r})^t$ is dualizable. 

Assume first that $\chi(v)$ and $\rank(v)$ are non-negative. 
See Remark \ref{rem-proof-of-dualizability-for-negative-euler-characteristic}
for the general case. 
We also assume, for simplicity of notation, that we can 
fix universal sheaves $\F_{v'}$ over the moduli spaces $\M_S(v')$ 
in our collection. In Remark 
\ref{rem-proof-of-dualizability-without-universal-sheaf}
we indicate the changes necessary when the universal sheaves exist only 
locally. 
We obtain the homomorphisms $\rho_{v'}$ 
(\ref{eq-exact-seq-is-locally-free-presentation-of-R1}),
the iterated blow-ups $\beta:B^{[k]}\M_S(v')\rightarrow \M_S(v')$,
and the elementary transforms 

\begin{equation}
\label{eq-elementary-transform-of-rho}
B^{[k]}(\rho_{v'}) :  B^{[k]}V_{v',0}\rightarrow B^{[k]}V_{v',1}. 
\end{equation}

\noindent
Recall that $B^{[k]}V_{v',0}$ is simply the pullback 
of $V_{v',0}$ by the iterated blow-up, while 
$B^{[k]}V_{v',1}$ is the result of an iteration 
of two operations: pullback via a single blow-up morphism followed by 
a Hecke-transformation along the exceptional divisor
(see (\ref{eq-blown-up-homomorphism})).


\begin{prop}
\label{prop-grassmannian-fibrations-of-blown-up-brill-noether-loci}
\begin{enumerate}
\item
\label{prop-item-constant-rank}
$B^{[k+1]}(\rho_{v'}):B^{[k+1]}V_{v',0}\rightarrow B^{[k+1]}V_{v',1}$ 
has constant rank 
over $B^{[k+1]}\M_S(v')^{k}$, $0\leq k \leq \mu(v')$. 
The kernel of its restriction 
\[
W_{v'}^k:= \ker\left(
(B^{[k+1]}\rho_{v'}\restricted{)}{B^{[k+1]}\M_S(v')^{k}}
\right)
\] 
is a  rank $k+\chi(v')$-vector bundle
over $B^{[k+1]}\M_S(v')^{k}$ which is a sub-sheaf of 
the pullback $\beta^*\left(p_*\left(\F_{v'}
\restricted{\right)}{\M_S(v')^{k}}\right)$. 
\item
\label{prop-item-dual-of-cokernel-is-subsheaf-of-ext}
The cokernel $U_{v'}^k$
of the restriction of $B^{[k+1]}(\rho_{v'})$ to $B^{[k+1]}\M_S(v')^{k}$
is a rank $k$ vector bundle
\[
U_{v'}^k:= \coker\left(
(B^{[k+1]}\rho_{v'}\restricted{)}{B^{[k+1]}\M_S(v')^{k}}
\right)
\] 
over $B^{[k+1]}\M_S(v')^{k}$. Its dual $(U_{v'}^k)^*$ is a subsheaf
of the twisted pullback 
\begin{equation}
\label{eq-twisted-pullback-ext-sheaf}
\beta^*\left(
\RelExt^1_{p_*}\left( (\F_{v'}\restricted{)}{\M_S(v')^{k}},\omega_{S}\right)
\right)
\left(+\sum_{i=k+1}^{\mu(v')}B^{[k+1]}\M_S(v')^i\right).
\end{equation}
\item
\label{prop-item-constant-rank-on-full-iterated-blow-up}
$B^{[1]}(\rho_{v'}):B^{[1]}V_{v',0}\rightarrow B^{[1]}V_{v',1}$ 
is surjective. 
The kernel $W_{v'}:=W_{v'}^0$ 
of its restriction is a rank $\chi(v')$-vector bundle
over $B^{[1]}\M_S(v')$. 
\item
\label{prop-item-universal-property-of-G-k-W}
The Grassmannian bundle  $G(k,W_{v'})$ over $B^{[1]}\M_S(v')$
satisfies the universal property analogous to that of
$G^0(k,\M_S(v'))$  for the restricted class of those 
families $(\E,W_T)$ over $T$ as in Proposition 
\ref{prop-universal-propery-of-G0} which satisfy the additional property:

\medskip
\noindent
The inverse image $(\kappa^{-1}\Ideal{\M_S(v')^t})\cdot\StructureSheaf{T}$ 
of the ideal sheaves $\Ideal{\M_S(v')^t}$ of all the Brill-Noether loci
with $1 \leq t \leq \mu(v')$ 
via the classifying morphism $\kappa:T\rightarrow \M_S(v')$ 
are all invertible sheaves of ideals on $T$.
\item
\label{prop-item-universal-property-of-B-M-k}
$B^{[k+1]}\M_S(v')^k$ satisfies the universal property 
analogous to that of $G_1(k,\M_S(v'))$ for the restricted class of those 
families $(\E,U_T)$ over $T$ as in Proposition 
\ref{prop-universal-propery-of-G1} which satisfy the additional property:

\medskip
\noindent
The inverse image $(\kappa^{-1}\Ideal{\M_S(v')^t})\cdot\StructureSheaf{T}$ 
of the ideal sheaves $\Ideal{\M_S(v')^t}$ of all the Brill-Noether loci
with $k+1 \leq t \leq \mu(v')$ 
via the classifying morphism $\kappa:T\rightarrow \M_S(v')$ 
are all all invertible sheaves of ideals on $T$.
\item 
\label{prop-item-blown-up-tyurin-isomorphism}
There is a canonical isomorphism 
\begin{equation}
\label{eq-blown-up-tyurin-isomorphism}
\DoubleTilde{f}_{v',k}:
B^{[k+1]}\M_S(v')^{k}
\IsomRightArrow
G(k,W_{v'+\vec{k}}).
\end{equation}
Denote by 
\[
\tilde{f}_{v',k}:B^{[k+1]}\M_S(v')^{k}
\rightarrow
B^{[1]}\M_S(v'+\vec{k})
\]
the composition of $\DoubleTilde{f}_{v',k}$ with the
projection to $B^{[1]}\M_S(v'+\vec{k})$. 
\item
\label{prop-item-smoothness}
$B^{[k+1]}\M_S(v')^k$ is smooth, non-empty and of codimension
$k\cdot(\Abs{\chi(v')}\!+k)$, for $0\leq k \leq \mu(v')$.
\item
\label{prop-item-tau-is-U-and-q-is-W}
Let $\tau^k_{v'}$ and $q^k_{v'}$ be the universal
sub- and quotient bundles of
$\tilde{f}^*_{v',k}(W_{v'+\vec{k}})$
(which are vector bundles over
$B^{[k+1]}\M_S(v')^k$).
There exists a line bundle $K$ on $B^{[k+1]}\M_S(v')^k$ 
and isomorphisms
\begin{eqnarray}
\label{eq-tau-is-twisted-coker}
\tau^{k}_{v'} \left(
-\sum_{j=k+1}^{\mu(v')}B^{[k+1]}\M_S(v')^j
\right)
& \cong & U^{k}_{v'}\otimes K
\\
q^{k}_{v'} & \cong &  W^{k}_{v'}\otimes K. 
\label{eq-q-is-twisted-ker}
\end{eqnarray}
\item
\label{prop-item-cotangent-bundle-is-same-hom-bundle-as-petri-bundle}
The isomorphism 
\begin{equation}
\label{eq-contraction-with-pulled-back-two-form}
N_{B^{[k+1]}\M_S(v')^k/B^{[k+1]}\M_S(v')} 
\IsomRightArrowOf{\Contract{\sigma_{\M_S(v')}}}
\Hom(q^{k}_{v'},\tau^{k}_{v'})\otimes
\StructureSheaf{B^{[k+1]}\M_S(v')^k}\left(
-\sum_{j=k+1}^{\mu(v')}B^{[k+1]}\M_S(v')^j
\right)
\end{equation}
induced by the symplectic structure of 
$\M_S(v')$ (see Lemma 
\ref{lemma-normal-bundle-is-twisted-relative-cotangent-bundle})
is conjugated via 
(\ref{eq-tau-is-twisted-coker}) and 
(\ref{eq-q-is-twisted-ker})
to the dual of the Petri map of $B^{[k+1]}(\rho_{v'})$
\begin{equation}
\label{eq-petri-homomorphism-on-iterated-blow-up-of-moduli}
\phi_{v',k}:\Hom(U^k_{v'},W^k_{v'})
\IsomRightArrow
N^*_{B^{[k+1]}\M_S(v')^k/B^{[k+1]}\M_S(v')}. 
\end{equation}
\end{enumerate}
\end{prop}

\noindent
{\bf Proof:}
The proof is by ascending induction on the pair $(\mu(v'),\mu(v')-k)$
with respect to the usual order:
$(a,b)>(a',b')$ if $a>a'$ or $a=a'$ and $b>b'$. 
Note that for a fixed Mukai vector $v$, the map 
\[(r,k) \mapsto (\mu(v+\vec{r}),\mu(v+\vec{r})-k)\]
is strictly decreasing in both $r$ and $k$. For each fixed $v$ 
the subset 
\[
\{(r,k) \ \mid \ v\!+\vec{r} \in \Hy \ \mbox{and} \ k\leq \mu(v\!+\!\vec{r})\}
\]
of the hyperbola in Figure \ref{eq-graph-of-hyperbola} is a finite set. 

\medskip
\noindent
The case $(\mu(v),\mu(v)-k) = (0,0)$. 
In this case $B^{[1]}\M_S(v)^0 = \M_S(v)$ and the Brill-Noether stratification
is trivial: $h^1(F)=0$ for every sheaf parametrized by $\M_S(v)$. 
$B^{[1]}\rho_v=\rho_v$ is surjective, $U^0_v$ is the zero bundle.
Both $\DoubleTilde{f}_{v,0}$ and 
$\tilde{f}_{v,0}$ are the identity.

\bigskip
\noindent
The induction step:
Assume that the proposition (and the dualizability conditions)
hold for all $(v',k')$ with $(\mu(v'),\mu(v')-k') < (\mu(v),\mu(v)-k)$.

\bigskip
\noindent
\ref{prop-item-constant-rank})
By the induction hypothesis $B^{[k+2]}(\rho_v)$ has constant rank
over $B^{[k+2]}(\M_S(v)^{k+1})$ and $B^{[k+2]}(\M_S(v)^{k+1})$ 
is smooth. 
Corollary \ref{cor-verification-of-codimension-condition} implies that 
the generic rank on every component of $B^{[k+1]}(\M_S(v)^{k})$ is
$k\!+\!\chi(v)$. Lemma 
\ref{lemma-compatibility-of-two-stratifications-on-exceptional-divisor}
part \ref{lemma-item-two-stratifications-coincide} implies that
$B^{[k+1]}(\rho_v)$ has constant rank over $B^{[k+1]}(\M_S(v)^{k})$.

\medskip
\noindent
\ref{prop-item-dual-of-cokernel-is-subsheaf-of-ext})
The only non-obvious statement is that $(U^k_v)^*$
is a subsheaf of the $\RelExt^1_p$-sheaf
(\ref{eq-twisted-pullback-ext-sheaf}).
Using Corollary
\ref{cor-ker-of-blown-up-e-vs-coker-of-blown-up-dual-e}
inductively, we get the isomorphisms
\begin{eqnarray*}
U^k_v := \coker(B^{[k+1]}\rho_v) & \cong &
\coker\left[
(B^{k+1}(B^{[k]}\rho_v)^*)^*
\right](-B^{[k+1]}\M_S(v)^{k+1}) \cong \cdots
\\
& \cong &
\coker\left[
(B^{[k+1]}\rho_v^*)^*
\right]\left(-\sum_{i=k+1}^{\mu(v)}B^{[k+1]}\M_S(v)^{i}\right). 
\end{eqnarray*}

\noindent
It remain to show that 
\[
\left(
\coker\left[
(B^{[k+1]}\rho_v^*)^*
\restricted{\right]}{B^{[k+1]}\M_S(v)^k}
\right)^*
\cong
\ker\left[
(B^{[k+1]}\rho^*_v\restricted{)}{B^{[k+1]}\M_S(v)^k}
\right]
\] 
is a subsheaf of 
$\beta^*\left(
\RelExt^1_{p_*}(\F_{v},\omega_{S}\restricted{)}{\M_S(v)^{k}}
\right)$. 
Part \ref{prop-item-constant-rank} implies
that $(B^{[k+1]}\rho^*_{v})$ has constant rank over
$B^{[k+1]}\M_S(v)^k$ and the kernel of its
restriction is a subsheaf of
$\beta^*\left(\ker(\rho^*_v\restricted{)}{\M_S(v)^k}\right)$.
Serre's Duality (\ref{eq-exact-seq-is-locally-free-presentation-of-Ext1})
identifies $\ker(\rho^*_v\restricted{)}{\M_S(v)^k}$ 
with 
$\RelExt^1_{p_*}(\F_{v},\omega_{S}\restricted{)}{\M_S(v)^{k}}$.

\medskip
\noindent
\ref{prop-item-constant-rank-on-full-iterated-blow-up})
This is a special case of part
\ref{prop-item-constant-rank}.

\medskip
\noindent
\ref{prop-item-universal-property-of-G-k-W})
%
We prove the analogue of Proposition 
\ref{prop-universal-propery-of-G0} in the case of a single family 
$\E$ of stable sheaves with Mukai
vector $v$ and $W_T$ a rank $k$ locally free subsheaf of $\pi_{T_*}\E$. 
The patching argument needed in the general case is identical to that in the  
proof of  Proposition \ref{prop-universal-propery-of-G0}. 
Assume that the classifying morphism 
$\kappa:T\rightarrow \M_S(v)$ of $\E$ pulls back all Brill-Noether loci
to Cartier divisors on $T$. 
Let $\tilde{\kappa}:T\rightarrow G^0(k,\M_S(v))$ be the classifying morphism
of $(\E,W_T)$. We claim that there exists a natural lift
$\kappa^{[1]}:T\rightarrow B^{[1]}\M_S(v)$ of $\kappa$. $\kappa^{[1]}$ 
exists by the universal property of blowing-up and induction. 
There is a subtlety in the $i$-th step of the induction: 
we need to check that if the inverse image in $T$ of 
$\Ideal{B^{[i\!+\!1]}\M_S(v)^t}$, $t\leq i$, 
via $\kappa^{[i\!+\!1]}$ is invertible on $T$, 
then the inverse image of its {\em strict} transform 
$\Ideal{B^{[i]}\M_S(v)^t}$ via $\kappa^{[i]}$ is also invertible on $T$. 
This follows from the comparison between the strict transform and the inverse 
image via the blow-up map in Theorem \ref{thm-vainsencher} part 
\ref{thm-item-vain1}. 
The existence of $\kappa^{[1]}$ follows. 
We have a canonical homomorphism (up to a scalar factor)
\[
\pi_{T_*}\E\otimes K \IsomRightArrow \kappa^*(p_*\F_v)
\]
for some line bundle
$K$ on $T$. 
Hence we have a canonical injective bundle homomorphism 
$W_T\otimes K\hookrightarrow \kappa^*V_0$ obtained as a composition
\begin{equation}
\label{eq-W-T-injects-to-V0}
W_T\otimes K \hookrightarrow (\pi_{T_*}\E)\otimes K \cong 
\kappa^*(p_*\F_v)\hookrightarrow \kappa^*V_0
\end{equation}
(see (\ref{eq-exact-seq-is-locally-free-presentation-of-R1})).
In particular, we have a canonical
morphism to the Grassmannian bundle $T\rightarrow G(k,V_0)$
over $\M_S(v)$, namely the composition
\[
T \RightArrowOf{\tilde{\kappa}} G(k,\M_S(v)) \hookrightarrow G(k,V_0).
\] 
But $\kappa^*V_0$ is isomorphic to 
$(\kappa^{[1]})^*(B^{[1]}V_0)$. Hence we get a canonical
morphism 
\[
\alpha:T \rightarrow G(k,B^{[1]}V_0)
\]
to the Grassmannian bundle over $B^{[1]}\M_S(v)$. 
We need to show that $\alpha$ 
factors through $G(k,W_v)$. 
By definition of $\tilde{\kappa}$, the composition
(\ref{eq-W-T-injects-to-V0})
factors through $(\tilde{\kappa})^*\tau_{G^0(k,\M_S(v))}$. 
We need to show that 
$(\tilde{\kappa})^*\tau_{G^0(k,\M_S(v))}$ is also in the kernel of

\begin{equation}
\label{eq-pullbak-of-blown-up-homo}
B^{[1]}(\kappa^*(\rho_v)):B^{[1]}\kappa^*(V_0)\rightarrow 
B^{[1]}\kappa^*(V_1)
\end{equation}

\noindent
(the kernel of (\ref{eq-pullbak-of-blown-up-homo}) 
is canonically isomorphic to
$(\kappa^{[1]})^*W_v$). 
$B^{[1]}\kappa^*(V_0)$ is equal to $\kappa^*(V_0)$
and $B^{[1]}\kappa^*(V_1)$ is, by definition, a sub-sheaf of
$\kappa^*(V_1)$ which contains the image of $\kappa^*(\rho_v)$. 
Hence, the sheaf theoretic kernel of (\ref{eq-pullbak-of-blown-up-homo})
is the same as that of $\kappa^*(\rho_v)$. We conclude that indeed 
$\alpha$ factors through $G(k,W^k_v)$. 


\medskip
\noindent
\ref{prop-item-universal-property-of-B-M-k})
Part \ref{prop-item-dual-of-cokernel-is-subsheaf-of-ext} implies the equality  
$B^{[k\!+\!1]}\M_s(v)^k = G(k,U^k_v)$. The rest of the proof is  
analogous to that of Part \ref{prop-item-universal-property-of-G-k-W}.

\medskip
\noindent
\ref{prop-item-blown-up-tyurin-isomorphism})
The statement is a tautology for $k=0$. Assume that $k \geq 1$. 
The induction hypothesis implies that $G(k,W_{v+\vec{k}})$ is smooth. 
By part \ref{prop-item-dual-of-cokernel-is-subsheaf-of-ext} for 
$(v,k)$ 
and the universal property of $G_1(k,M_S(v))$ (Proposition
\ref{prop-universal-propery-of-G1}) there is a canonical classifying morphism 
\[
\alpha: B^{[k+1]}\M_S(v)^k \rightarrow G_1(k,\M_S(v))
\]
compatible with the Brill-Noether stratifications. 
Theorem \ref{thm-the-tyurin-extension-isomorphism} implies that we have a 
canonical isomorphism
\begin{equation}
\label{eq-this-isomorphism-f-will-be-replaced-in-the-negative-euler-char-case}
f: G_1(k,\M_S(v)) \IsomRightArrow G^0(k,\M_S(v+\vec{k}))
\end{equation}
compatible (up to a shift by $k$) with the Brill-Noether stratification.
We claim that the composition $f\circ \alpha$ lifts to morphism
\[
\DoubleTilde{f}_{v,k}:  B^{[k+1]}\M_S(v)^k \rightarrow
G^0(k,W_{v+\vec{k}}).
\]
Indeed, this would follow from 
the universal property of $G^0(k,W_{v+\vec{k}})$ 
(part \ref{prop-item-universal-property-of-G-k-W} of the Proposition)
provided that we check that $f\circ \alpha$ pulls back all the 
Brill-Noether strata to Cartier divisors on 
$B^{[k+1]}\M_S(v)^k$. 
By the compatibility of $f$ and $\alpha$ with respect to the 
(total-transform of the) Brill-Noether stratifications, 
all we need to check is that the subschemes 
$B^{[k+1]}\M_S(v)^t$, $t\geq k+1$ are all Cartier divisors
in $B^{[k+1]}\M_S(v)^k$. This is indeed the case. 

Conversely, $G(k,W_{v+\vec{k}})$ admits a morphism  to 
$G(k,\M_S(v+\vec{k}))$  
\[
g: G(k,W_{v+\vec{k}}) \rightarrow G(k,\M_S(v+\vec{k})).
\] 
The existence of $g$ follows from the universal property of 
$G(k,\M_S(v+\vec{k}))$ (Proposition
\ref{prop-universal-propery-of-G0}).
The morphism $g$ is compatible with the Brill-Noether stratification.
Hence, the composition 
\[
f^{-1}\circ g: G(k,W_{v+\vec{k}}) \rightarrow 
G_1(k,\M_S(v))
\]
is compatible with the 
Brill-Noether stratification (up to a shift by $k$).
We claim that $f^{-1}\circ g$ lifts to a morphism
\[
\DoubleTilde{g}_{v,k}: G(k,W_{v+\vec{k}}) \rightarrow
B^{[k+1]}\M_S(v)^k.
\] 
Indeed, that  would follow from the universal property of 
$B^{[k+1]}\M_S(v)^k$
(part \ref{prop-item-universal-property-of-B-M-k} of the Proposition)
provided that we check that $f^{-1}\circ g$ pulls back all the 
Brill-Noether strata to Cartier divisor on $G(k,W_{v+\vec{k}})$. 
By the compatibility of $f^{-1}\circ g$ with respect to the 
(total-transform of the) Brill-Noether stratifications, 
all we need to check is that the subschemes 
$G^0(k,W_{v+\vec{k}})^t$, $1\leq t \leq \mu(v)-k$
are all Cartier divisors
in $G(k,W_{v+\vec{k}})$. This is indeed the case.

Both $\DoubleTilde{f}_{v,k}$ and $\DoubleTilde{g}_{v,k}$ 
are birational isomorphisms and 
$\DoubleTilde{f}_{v,k}\circ \DoubleTilde{g}_{v,k}$,
$\DoubleTilde{g}_{v,k}\circ \DoubleTilde{f}_{v,k}$ are regular morphisms
which are generically the identity. $G(k,W_{v+\vec{k}})$ 
is smooth by the induction hypothesis. Hence 
$\DoubleTilde{f}_{v,k}\circ \DoubleTilde{g}_{v,k}$ is the identity. 
By properness, both $\DoubleTilde{f}_{v,k}$ and $\DoubleTilde{g}_{v,k}$ 
are surjective. Hence, $\DoubleTilde{f}_{v,k}$ is also injective. 
Since $G(k,W_{v+\vec{k}})$ 
is smooth, $\DoubleTilde{f}_{v,k}$ must be an isomorphism 
(by Zariski's Main Theorem).

\medskip
\noindent
\ref{prop-item-smoothness}) The smoothness and codimension 
formula follow immediately from
the induction hypothesis
and part \ref{prop-item-blown-up-tyurin-isomorphism}. 
The non-emptiness follows from Corollary 
\ref{cor-verification-of-codimension-condition}. 

\medskip
\noindent
\ref{prop-item-tau-is-U-and-q-is-W})
We identify $B^{[k+1]}\M_S(v)^k$ with $G^0(k,W_{v+\vec{k}})$. 
We denote the pullback of the universal sheaves to
$B^{[k+1]}\M_S(v)^k \times S$ by the same notation $\F_v$ and 
$\F_{v+\vec{k}}$. 
Consider the two exact sequences 
(\ref{eq-functorial-quotient-by-trivial-subbundle}) and
(\ref{eq-functorial-tautological-extension}) over $G$ where $G$ denotes both
$G^0(k,\M_S(v\!+\!\vec{k}))$ and $G_1(k,\M_S(v))$. 
Pulling back (\ref{eq-functorial-quotient-by-trivial-subbundle}) and
(\ref{eq-functorial-tautological-extension})
to $B^{[k+1]}\M_S(v)^k \times S$ via the classifying morphisms of the families 
$[\F_{v+\vec{k}},\tau^k_v]$ and
$[\F_v,(U^k_v)^*\left(
-\sum_{i=k+1}^{\mu(v)}B^{[k+1]}\M_S(v)^i
\right)]$
we get the two short exact sequences:
\begin{equation}
\label{eq-blown-up-universal-quotient-by-trivial-subbundle}
0 \rightarrow
p^*\tau^k_v 
\rightarrow
\F_{v+\vec{k}}
\rightarrow
Q
\rightarrow
0, \ \ \mbox{and} 
\end{equation}

\begin{equation}
\label{eq-blown-up-universal-tautological-extension}
0 \rightarrow
p^*U^k_v \left(
\sum_{i=k+1}^{\mu(v)}B^{[k+1]}\M_S(v)^i
\right)
\rightarrow
\E
\rightarrow
\F_v
\rightarrow
0. 
\end{equation}

We already know that there exists a line bundle $K$ on $G_1(k,\M_S(v))$
such that (\ref{eq-functorial-tautological-extension})$\otimes K$ is
isomorphic to (\ref{eq-functorial-quotient-by-trivial-subbundle}). 
Hence the pullback of this line bundle to $B^{[k+1]}\M_S(v)^k$ (denoted 
also by $K$) gives the isomorphism between the short exact sequences 
(\ref{eq-blown-up-universal-tautological-extension})$\otimes K$ 
and (\ref{eq-blown-up-universal-quotient-by-trivial-subbundle}). 
In particular, we get the isomorphism 
(\ref{eq-tau-is-twisted-coker}). 

\medskip
Next we construct the isomorphism (\ref{eq-q-is-twisted-ker}). 
Pushing forward 
(\ref{eq-blown-up-universal-tautological-extension})$\otimes K$ 
and (\ref{eq-blown-up-universal-quotient-by-trivial-subbundle})
via the projection 
$p:B^{[k+1]}\M_S(v)^k \times S \rightarrow B^{[k+1]}\M_S(v)^k$ 
we get the short exact sequence
\begin{equation}
\label{eq-push-forward-of-universal-tautological-extension}
0
\rightarrow
\tau^k_v
\rightarrow
p_*\F_{v+\vec{k}}
\rightarrow
(p_*\F_v)\otimes K
\rightarrow 0.
\end{equation}

\noindent
Both $q_v^k$ and $W^k_v\otimes K$ are subsheaves of 
$(p_*\F_v)\otimes K$ of maximal rank. 
At a generic point (in the dense open subset
$B^{[k+1]}\M_S(v)^k \setminus \cup_{i=k+1}^{\mu(v)}B^{[k+1]}\M_S(v)^i$) 
all three sheaves are equal. 
Moreover, $(p_*\F_v)\otimes K$ is a
subsheaf of the vector bundle $p_*\F_v(n)\otimes K=B^{[k+1]}V_{v,0}\otimes K$
(see (\ref{eq-exact-seq-is-locally-free-presentation-of-R1})).
By part \ref{prop-item-constant-rank}, 
$W^k_v\otimes K$ is a subbundle of $B^{[k+1]}V_{v,0}\otimes K$. 
The short exact sequence 
(\ref{eq-push-forward-of-universal-tautological-extension})
implies that $q_v^k$ is also a subbundle of $p_*\F_v(n)\otimes K$ 
(the injective sheaf homomorphism is also injective on each fiber). 
These two subbundles of $B^{[k+1]}V_{v,0}\otimes K$ are equal on a dense open
subset. We get two sections of the Grassmannian bundle 
$G(\chi(v)+k,B^{[k+1]}V_{v,0}\otimes K)\rightarrow B^{[k+1]}\M_S(v)^k$ 
which are equal on a dense open subset. 
Smoothness of $B^{[k+1]}\M_S(v)^k$ implies that they are 
globally equal.

\medskip
\noindent
\ref{prop-item-cotangent-bundle-is-same-hom-bundle-as-petri-bundle}) 
It suffices to prove the statement over $\M_S(v)^k\setminus \M_S(v)^{k+1}$.
Let $F$ be such a sheaf and
\[
0 \rightarrow H^1(F)\otimes\omega_S \HookRightArrowOf{\iota} E \rightarrow F
\rightarrow 0
\]
the natural extension. It is known that the Petri map 
$\phi:H^1(F)^*\otimes H^0(F) \rightarrow T^*_{[F]}\M_S(v)$
becomes the Yoneda product under the natural identifications
$H^1(F)^*\cong \Ext^1(F,\omega_S)$, 
$H^0(F)\cong \Hom(\omega_S,F\otimes\omega_S)$, and 
$T^*_{[F]}\M_S(v)\cong\Ext^1(F,F\otimes\omega_S)$. 
Let $W\in G(k,H^0(E))$ represent the point $\iota(H^1(F)\otimes H^0(\omega_S))$
and 
\[
\psi : T_{W}G(k,H^0(E)) \rightarrow T_{[F]}\M_S(v)
\]
the differential of the natural embedding. It is easy to verify that 
$\psi$ is also induced by the Yoneda product under the natural identification 
of $T_{W}G(k,H^0(E))$ with $\Ext^1(F,\StructureSheaf{S})\otimes H^0(F)$ 
and $T_{[F]}\M_S(v)$ with $\Ext^1(F,F)$. It follows that the following 
diagram commutes
\[
\begin{array}{ccc}
\Ext^1(F,\omega_S)\otimes H^0(F) & \LongRightArrowOf{\phi} &
\Ext^1(F,F\otimes\omega_S)
\\
\uparrow & & \uparrow
\\
\Ext^1(F,\StructureSheaf{S})\otimes H^0(F) & \LongRightArrowOf{\psi} &
\Ext^1(F,F)
\end{array}
\]
where the vertical isomorphisms are induces by cup-product with the symplectic
structure of $S$. Part
\ref{prop-item-cotangent-bundle-is-same-hom-bundle-as-petri-bundle}
follows since the right vertical isomorphism is the one induced by the 
symplectic structure on $\M_S(v)$.
\EndProof

\begin{cor}
\label{cor-the-collection-is-admissible}
The collection
$
\{\M_S(v+\epsilon\cdot\vec{r})^t
\mid \ 0 \leq r \leq \mu(v) \}
$
is dualizable. Above $\epsilon$ is $1$ if $\chi(v)$ is positive,
$-1$ if  $\chi(v)$ is negative and either one if $\chi(v)=0$.
\end{cor}

\noindent
{\bf Proof:}
Dualizability condition \ref{cond-codimension} follows from 
Mukai's dimension formula $\dim \M_S(v) = 2+\langle v,v\rangle$.
Conditions \ref{cond-abstract-tyurin-morphism} and 
\ref{cond-compatibility-of-grassmanian-bundles}
are verified in part 
\ref{prop-item-blown-up-tyurin-isomorphism} of Proposition 
\ref{prop-grassmannian-fibrations-of-blown-up-brill-noether-loci}
(where $\PP{W}_{B^{[1]}\M_S(v+\vec{t})}$ is denoted by
$\PP{W}_{v+\vec{t}}$).
The compatibility of the two stratifications in 
Condition \ref{cond-compatibility-of-stratifications-wrt-abstract-tyurin}
is proven exactly as in the last part of the proof of
Theorem \ref{thm-the-tyurin-extension-isomorphism}. 
Condition \ref{cond-compatibility-of-stratifications-on-proper-transforms}
follows from Lemma 
\ref{lemma-compatibility-of-two-stratifications-on-exceptional-divisor}
part \ref{lemma-item-two-stratifications-coincide} and 
part \ref{prop-item-cotangent-bundle-is-same-hom-bundle-as-petri-bundle}
of Proposition 
\ref{prop-grassmannian-fibrations-of-blown-up-brill-noether-loci}.
\EndProof

\begin{rem}
\label{rem-proof-of-dualizability-without-universal-sheaf}
{\rm
In the absence of  global universal sheaves, the statement of Proposition
\ref{prop-grassmannian-fibrations-of-blown-up-brill-noether-loci}
should be modified as follows:
The vector bundles $B^{[k]}V_{v',0}$, $B^{[k]}V_{v',1}$, $W^k_{v'}$,
$U^k_{v'}$ exist only locally, but the corresponding projective bundles 
exist globally. The {\em vector} bundles 
$\SheafHom(B^{[k]}V_{v',0},B^{[k]}V_{v',1})$ and their sections 
$B^{[k]}(\rho_{v'})$ in (\ref{eq-elementary-transform-of-rho}) 
exist globally (see Remark \ref{rem-non-existence-of-universal-sheaf}). 
The statement involving equation 
(\ref{eq-twisted-pullback-ext-sheaf}) is valid only locally.
The {\em vector} bundles $\SheafHom(q^k_{v'},\tau^k_{v'})$ in
(\ref{eq-contraction-with-pulled-back-two-form}) and 
$\SheafHom(W^k_{v'},U^k_{v'})$ in 
(\ref{eq-petri-homomorphism-on-iterated-blow-up-of-moduli}) exist globally.
Equations (\ref{eq-tau-is-twisted-coker}) and 
(\ref{eq-q-is-twisted-ker}) should be replaced by their projectivization 
but the simultaneous conjugation by one and the inverse of the other in part 
\ref{prop-item-cotangent-bundle-is-same-hom-bundle-as-petri-bundle} 
of the proposition is a global isomorphism of vector bundles. 

The proof goes through if we replace the global family $\F_v$ 
by an open covering of $\M_S(v)$ (in the complex or \'{e}tale topology) 
and a local family over each open set. The main point is that 
any choice of gluing transformations (a cochain over $\M_S(v)\times S$)
fails only slightly to be a cocycle. 
The coboundary of any such choice of gluing transformations 
is necessarily the pullback of a \v{C}eck $2$-coboundary for the sheaf 
$\StructureSheaf{\M_S(v)}^\times$ of invertible functions over $\M_S(v)$. 
This follows from the simplicity of the sheaves involved (see 
\cite{mukai-hodge} Appendix 2). 
}
\end{rem}

\subsection{Lazarsfeld's reflection isomorphism}
\label{subsec-lazarsfeld-reflection-isomorphism}

\begin{thm}
\label{thm-lazarsfeld-reflection-isomorphism}
\begin{enumerate}
\item
There exists a natural isomorphism
\begin{equation}
\label{eq-lazarsfeld-reflection-on-t-stratum}
\tilde{q}_t: G^0(a+b+t,\M_S(a,\LB,b)) \ \LongIsomRightArrow \ 
G^0(a+b+t,\M_S(b+t,\LB,a+t)), 
\end{equation}
for all integers $a$, $b$, $t$  satisfying $a\geq 0$ and $b+t\geq 0$. 
Equivalently, 
there exist  natural isomorphisms 
\begin{eqnarray*}
G^0(k,\M_S(v)) & \LongIsomRightArrow & G^0(k,\M_S(\sigma\circ\tau(v-\vec{k}))),
\ \ \ \mbox{for} \ \ \ k\geq \rank(v)\geq 0, \ \ \mbox{and}
\\
G^0(\chi(v),\M_S(v)) & \LongIsomRightArrow & G^0(\chi(v),\M_S(\sigma(v)))
\hspace{8ex} \mbox{for}  \ \ \chi(v)\geq \rank(v)\geq 0.
\end{eqnarray*}
\item
The compositions $\tilde{q}_t \circ \tilde{q}_{-t}$ and 
$\tilde{q}_{-t} \circ \tilde{q}_{t}$ are both the identity morphism.
\end{enumerate}
\end{thm}

\begin{example}
\label{example-reflection-among-hilbert-schemes}
{\rm
Consider the case where $v=(1,\LB,g\!-\!d)$, 
$1\leq d \leq 2g\!-\!2$,
and $k=2$. Then $\M_S(v)=S^{[d]}$, $G^0(2,S^{[d]})$ parametrizes pairs
consisting of a $\PP^{g-2}$ in $\PP^g$ containing a length $d$ subscheme $D$, 
while $\M_S(\sigma\circ\tau(v\!-\!\vec{2}))=S^{[2g\!-\!2\!-\!d]}$. 
The $\PP^{g-2}$ intersects $S$ in a length $2g\!-\!2$ subscheme $\tilde{D}$ 
and the isomorphism 
\[
G^0(2,S^{[d]}) \ \cong \ G^0(2,S^{[2g\!-\!2\!-\!d]})
\]
maps a pair $D\subset\PP^{g-2}$ to the complementary pair 
$D^\perp\subset\PP^{g-2}$. If $\tilde{D}$ is reduced, 
then $D^\perp$ is the set theoretic complement  $\tilde{D}\setminus D$.
For a general complete intersection $\tilde{D}$, the dualizing sheaf 
$\omega_{\tilde{D}}$ is a free   
$\StructureSheaf{\tilde{D}}$-module of rank $1$ and the ideal
$\Ideal{\tilde{D},D^\perp}$ of $D^\perp$ as a subscheme of $\tilde{D}$ 
is the annihilator of 
$\Ideal{\tilde{D},D}\otimes_{\StructureSheaf{\tilde{D}}} \omega_{\tilde{D}}$
under the perfect local duality pairing
$
H^0(\StructureSheaf{\tilde{D}})  \otimes H^0(\omega_{\tilde{D}})
\rightarrow \ComplexNumbers. 
$
}
\end{example}

Theorem \ref{thm-lazarsfeld-reflection-isomorphism} is geometrically 
intuitive also in the following  special case:

\medskip
\noindent
{\bf The case $(a,b)=(1,0)$ and $t=0$}: 
Let $\C \subset S \times \linsys{\LB}$ be the universal curve. 
$\M_S(1,\LB,0)$ is $S^{[g]}$ and 
$\M_S(0,\LB,1)$ is the relative Picard $\Pic^g(\C/\linsys{\LB})$. 
Both $G(1,\M_S(1,\LB,0))$ 
and $G(1,\M_S(0,\LB,1))$ are isomorphic to the universal relative 
Hilbert scheme $\Hilb_g(\C)$ 
of length $g$ subschemes of curves in the linear system 
$\linsys{\LB}$. 

\begin{equation}\label{diagram-abel-jacobi}
{\divide\dgARROWLENGTH by 2
\begin{diagram}
\node[4]{\Hilb_g(\C)}
\arrow{sw,l}{a}
\arrow{se,t}{e} \\
\node{\Pic^g_{\C}}
\node[2]{\Pic^{-g}_{\C}}
\arrow[2]{w,t}{\SheafHom_{\C}(\bullet,\StructureSheaf{\C})}
\node[2]{S^{[g]}}
\end{diagram}
}
\end{equation}

\noindent
Note that both $e$ and $a$ are surjective regular morphisms.
 The morphism 
$e:\Hilb_g(\C) \rightarrow S^{[g]}$ 
is a birational isomorphism since $g$ points in 
general position on $S$ determine a unique curve in the linear system 
$\linsys{\LB}$ passing through them. 
The morphism  
$a:\Hilb_g(\C/\linsys{\LB}) \rightarrow J^{-g}$ 
is the Abel-Jacobi map (see \cite{altman-kleiman} for its construction). 
It is a birational isomorphism by the 
Abel-Jacobi theorem. Given an integral curve $C$ with planar singularities, 
the morphism $\SheafHom_{C}(\bullet,\StructureSheaf{C})$ is an involution of 
the compactified Picard.  Set theoretically, 
we need only the fact that a rank $1$ torsion free sheaf $F$ on $C$ is a 
reflecsive $\StructureSheaf{C}$-module and
$\chi\left(\SheafHom_{C}(F,\omega_{C})\right) = -\chi(F)$. This is a special 
case of a theorem of Serre which states that the functor 
$\SheafExt^c(\bullet,\StructureSheaf{S})$ is exact and involutive on the 
category of Cohen-Macaulay modules of codimension $c$. Serre's Theorem 
enters the picture via:

\begin{new-lemma}
\label{lemma-hom-over-C-is-ext1-over-S}
Let $F$ be an $\StructureSheaf{S}$-module  with pure one-dimensional support 
$C\subset S$. Then $F$ is a Cohen-Macaulay module and 
we have a natural isomorphism 
\[
\SheafHom_S(F,\omega_C) \cong \SheafExt^1_S(F,\StructureSheaf{S}).
\]
\end{new-lemma}

\noindent
{\bf Proof:}
Consider the long exact sequence of extension sheaves obtained by taking 
$\SheafHom_S(F,\bullet)$ of the short exact sequence
\[
0\rightarrow \StructureSheaf{S}\rightarrow \StructureSheaf{S}(C)
\rightarrow \omega_C \rightarrow 0. 
\]
The homomorphism 
$\SheafExt^1_S(F,\StructureSheaf{S})\rightarrow 
\SheafExt^1_S(F,\StructureSheaf{S}(C))$ 
vanishes since $C$ is also the support of 
$\SheafExt^1_S(F,\StructureSheaf{S})$. 
Hence the connecting homomorphism is an isomorphism. 
\EndProof

\bigskip
\noindent
{\bf Proof:} (of Theorem \ref{thm-lazarsfeld-reflection-isomorphism})
The proof consists of two steps. In the first step we describe the bijection 
(\ref{eq-lazarsfeld-reflection-on-t-stratum}) set theoretically 
taking care of stability issues. In the second step we work in families and 
prove that the two functors coarsely represented by the two moduli spaces
are equivalent.

\smallskip
\noindent
{\bf Step I} (stability) 
The set theoretic description of the bijection 
(\ref{eq-lazarsfeld-reflection-on-t-stratum}) is given below in 
(\ref{eq-lazarsfelds-reflection-of-a-pair-rk-U-larger-rk-F})
and
(\ref{eq-lazarsfelds-reflection-of-a-pair-on-a-curve}).
A stable sheaf $F\in \M_S(a,\LB,b)$ and a subspace 
$U\in G(a+b+t,H^0(F))$ determine the exact sequence:

\begin{equation}
\label{eq-E-U-F}
0 \rightarrow
E(\tilde{F},U) \HookRightArrowOf{i}
U\otimes_{\ComplexNumbers}\StructureSheaf{S} \RightArrowOf{ev}
\tilde{F} \rightarrow 0
\end{equation}

\noindent
where $\tilde{F}$ is the subsheaf of $F$ generated by $U$. 
Our assumptions imply that $a+b+t \geq a$. 
We have two cases which are slightly different:
$a+b+t=a$ and $a+b+t > a$. 

\noindent
{\bf The case $a+b+t > a$}:
Lemma \ref{lazarsfeld-slope-stability-lemma} 
implies that $U$ generates $F$ away from a zero-dimensional subscheme. 
Thus $F/\tilde{F}$ has a zero-dimensional support
\begin{equation}
\label{eq-F-mod-tilde-F}
0 \rightarrow \tilde{F} \RightArrowOf{e} F \rightarrow F/\tilde{F} \rightarrow 0.
\end{equation}

\noindent
If $F$ is generated by $U$, we define $q(F,U)$ to be the vector bundle 
$E(F,U)^*\in \M_S(b+t,\LB,a+t)$.  Otherwise, 
we proceed to define the stable torsion free sheaf $q(F,U)$ as a sub-sheaf of
$E(\tilde{F},U)^*$ (see (\ref{eq-def-of-q-F-U})). 
We have the locally free presentations 
\begin{eqnarray}
0 & \rightarrow &
F  \hookrightarrow  
F^{**} \rightarrow F^{**}/F \rightarrow 0 \ \ \mbox{and} 
\label{eq-presentation-of-F}
\\
0 & \rightarrow &
\tilde{F} \hookrightarrow F^{**} \rightarrow F^{**}/\tilde{F} \rightarrow 0.
\end{eqnarray}
Note that $E(\tilde{F},U)$ in (\ref{eq-E-U-F})
is  locally free and 

\begin{equation}
\label{eq-locally-free-resolution-via-double-dual}
0 \rightarrow 
E(\tilde{F},U) \HookRightArrowOf{i}
U\otimes_{\ComplexNumbers}\StructureSheaf{S} \RightArrowOf{ev}
F^{**} \rightarrow
F^{**}/\tilde{F} \rightarrow 0
\end{equation}

\noindent
is a locally free resolution of $F^{**}/\tilde{F}$.
The Local Duality Theorem (\cite{griffiths-harris} page 693)
implies that 

\begin{eqnarray}
\SheafExt^1_S(F^{**}/F,\StructureSheaf{S}) & = &
\SheafExt^1_S(F/\tilde{F},\StructureSheaf{S}) = 0,
\label{eq-sheaf-ext1-of-a-sheaf-with-zero-dimensional-support-with-O-S-vanish}
\\
\SheafExt^2_S(F,\StructureSheaf{S}) & = & 0, \ \ \mbox{and}, 
\ \ \mbox{using} \ \ \mbox{(\ref{eq-presentation-of-F})}, 
\label{eq-sheaf-ext2-of-a-torsion-free-sheaf-with-O-S-vanish}
\\
\SheafExt^1_S(F,\StructureSheaf{S}) & \cong & 
\SheafExt^2_S(F^{**}/F,\StructureSheaf{S}).
\label{eq-eq-sheaf-ext1-of-a-torsion-free-sheaf-is-sheaf-ext2}
\end{eqnarray}

\noindent
Moreover, we have a perfect pairing
\begin{equation}
\label{eq-local-duality-pairing}
H^0(\SheafExt^2(F/\tilde{F},\omega_S))
\otimes_{\ComplexNumbers}H^0(F/\tilde{F}) 
\rightarrow \ComplexNumbers. 
\end{equation}

\noindent
Taking $\SheafHom(\bullet,\StructureSheaf{S})$ of (\ref{eq-F-mod-tilde-F})
and using
(\ref{eq-sheaf-ext1-of-a-sheaf-with-zero-dimensional-support-with-O-S-vanish})
and 
(\ref{eq-sheaf-ext2-of-a-torsion-free-sheaf-with-O-S-vanish})
we get the short exact piece of the long exact sequence of 
local Exts:

\begin{equation}
\label{eq-short-exact-sequence-of-local-ext}
0 \rightarrow
\SheafExt^1_S(F,\StructureSheaf{S}) \HookRightArrowOf{e^*}
\SheafExt^1_S(\tilde{F},\StructureSheaf{S}) \rightarrow 
\SheafExt^2_S(F/\tilde{F},\StructureSheaf{S}) \rightarrow 0.
\end{equation}

\noindent
Taking $\SheafHom(\bullet,\StructureSheaf{S})$ of 
(\ref{eq-E-U-F}) we get the long exact sequence

\begin{equation}
\label{eq-Hom-of-E-U-F}
0 \rightarrow F^* \RightArrowOf{ev^*}
U^*\otimes_{\ComplexNumbers}\StructureSheaf{S}
\RightArrowOf{i^*}
E(\tilde{F},U)^* \RightArrowOf{\eta}
\SheafExt^1_S(\tilde{F},\StructureSheaf{S}) \rightarrow 0.
\end{equation}

\noindent
Moding out (\ref{eq-Hom-of-E-U-F}) by $F^*$ we get an extension class
\[
\epsilon(\tilde{F},U) \in \Ext^1_S\left[
\SheafExt^1_S(\tilde{F},\StructureSheaf{S}),\coker(ev^*)
\right].
\]

\noindent
Pulling back (\ref{eq-Hom-of-E-U-F})
via the injective homomorphism $e^*: \SheafExt^1_S(F,\StructureSheaf{S})
\hookrightarrow \SheafExt^1_S(\tilde{F},\StructureSheaf{S})$ 
we get the class 
\[
\epsilon(F,U) \in \Ext^1_S\left[
\SheafExt^1_S(F,\StructureSheaf{S}),\coker(ev^*)
\right]
\]
representing an exact sequence defining $q(F,U)$

\begin{equation}
\label{eq-def-of-q-F-U}
0 \rightarrow F^* \RightArrowOf{ev^*} 
U^*\otimes_{\ComplexNumbers}\StructureSheaf{S} \RightArrowOf{i^*}
q(F,U) \RightArrowOf{\eta\circ e^*} 
\SheafExt^1_S(F,\StructureSheaf{S}) \rightarrow 0.
\end{equation}

\noindent
Moding out (\ref{eq-Hom-of-E-U-F}) by (\ref{eq-def-of-q-F-U}) 
and using the exactness of (\ref{eq-short-exact-sequence-of-local-ext})
we see that the sheaf $q(F,U)$ fits in the short exact sequence

\begin{equation}
\label{eq-q-F-U-is-a-subsheaf-of-E}
0 \rightarrow q(F,U) \hookrightarrow 
E(\tilde{F},U)^* \rightarrow \SheafExt^2_S(F/\tilde{F},\StructureSheaf{S})
\rightarrow 0.
\end{equation}

\noindent
By definition (\ref{eq-E-U-F})
\[
c_2(E(\tilde{F},U)) = 2g-2-c_2(\tilde{F}) = 2g-2 - c_2(F) + 
length(F/\tilde{F}).
\]

\noindent
Local Duality (\ref{eq-local-duality-pairing}) 
and (\ref{eq-q-F-U-is-a-subsheaf-of-E})
imply that $c_2(q(F,U)) = 2g-2-c_2(F).$ 
Clearly, $rank(q(F,U)) = b+t$ and $\det(q(F,U)) = \LB$. Hence,
the Mukai vector of $q(F,U)$ is 
\[
v(q(F,U))=(b+t,\LB,a+t). 
\]
It is easy to see that the dual sheaf does not have any global sections 
$H^0(q(F,U)^*) = H^0(E(\tilde{F},U))= 0$. 
Moreover, $q(F,U)$ is generated by its global sections away from a 
zero-dimensional sub-scheme (see 
(\ref{eq-def-of-q-F-U})
and
(\ref{eq-eq-sheaf-ext1-of-a-torsion-free-sheaf-is-sheaf-ext2})). 
Lemma \ref{lazarsfeld-slope-stability-lemma} implies that 
$q(F,U)$ is stable. As $H^0(F^*)$ vanishes, the homomorphism $i^*$
in (\ref{eq-def-of-q-F-U}) embeds $U^*$ in $H^0(S,q(F,U))$. We get a pair
$(q(F,U),W) \in G^0(a+b+t,\M(b+t,\LB,a+t))$. 
We define (\ref{eq-lazarsfeld-reflection-on-t-stratum}) 
by
\begin{equation}
\label{eq-lazarsfelds-reflection-of-a-pair-rk-U-larger-rk-F}
\tilde{q}_t(F,U) := (q(F,U),W).
\end{equation}

\medskip
\noindent
{\bf The case $(a,b)=(1,0)$ and $t=0$ revisited}: 
This case is part  of the more general case $a+b+t=a$. 
We treat it separately as a warm-up. 
Let $(\LB(-D),s)$ be a pair consisting of the $\LB$-twisted ideal sheaf 
of a length $g$ subscheme
$D \subset S$ and a section $s\in \PP H^0(S,\LB(-D))$. The zero-locus of
$s$ determines, as a section of $\LB$, a curve $C$ containing $D$ as a 
subscheme.
We get the triple $(C,\Ideal{C,D}\otimes\LB,s_1)$ where 
$s_1\in \PP\Ext_S^1(\Ideal{C,D}\otimes\LB,\StructureSheaf{S})$ is 
the extension class
of 
\[
0 \rightarrow \StructureSheaf{S} \RightArrowOf{\cdot s}
\LB(-D) \rightarrow \Ideal{C,D}\otimes\LB  \rightarrow 0.
\]
Serre's Duality on $S$ identifies 
$\PP\Ext_S^1(\Ideal{C,D}\otimes\LB,\StructureSheaf{S})$
with $\PP H^1(S,\Ideal{C,D}\otimes\LB)^*$. Serre's Duality on $C$ 
identifies the latter with $\PP H^0(C,\Hom(\Ideal{C,D},\StructureSheaf{C}))$. 
Denote the corresponding section by $s_0$. 
The triple $(C,\Hom(\Ideal{C,D},\StructureSheaf{C}),s_0)$ is 
a point in $G^0(1,\Pic^g(\C/\linsys{\LB}))$.

Conversely, given a triple $(C,L,s_0)$ representing a point in
$G^0(1,\Pic^g(\C/\linsys{\LB}))$ we use Serre's Duality on $S$ and on $C$ to 
interpret $s_0$ as an extension class $s_1$ of $\StructureSheaf{S}$-modules 
\[
0 \rightarrow \StructureSheaf{S}  \hookrightarrow
F \RightArrowOf{j} \SheafHom(L,\omega_C) \rightarrow 0.
\]
Lemma \ref{lemma-tyurins-extension-is-stable} implies that
$F$, which is $f(\SheafHom(L,\omega_C),span(s_1))$, is a stable rank $1$ 
torsion free sheaf. 
Composing $j$ with the evaluation at $s_0$ 
$\Hom(L,\omega_C) \RightArrowOf{s_0} \omega_C$ we
get a non-zero homomorphism  
\[
\bar{s}_0:F \rightarrow \omega_C \cong \restricted{\LB}{C}.
\]
It is easy to check that $\bar{s}_0$ maps to zero by the 
connecting homomorphism $\delta$ of the long exact sequence 
\[
0 = \Hom(F,\StructureSheaf{S}) \rightarrow \Hom(F,\LB) \rightarrow
\Hom(F,\restricted{\LB}{C}) \RightArrowOf{\delta} \Ext^1(F,\StructureSheaf{S})
\]
(replace $(F,V)$ in 
(\ref{eq-connecting-hommomorphism-from-h1-to-h2-for-tyurin-extension})
by $(\SheafHom(L,\omega_C),span(s_1))$ and apply Serre's Duality).
Hence $\bar{s}_0$ determines a non-trivial homomorphism 
\[
i : F \hookrightarrow \LB.
\]
We recover an ideal sheaf 
of a length $g$ subscheme $D \subset S$ by setting $\Ideal{S,D}$ 
to be the image of
$F\otimes \LB^{-1}\HookRightArrowOf{i} \StructureSheaf{S}$.

\bigskip
\noindent
{\bf The case $a+b+t = a$}:
In that case, $t=-b$ and $\M_S(b+t,\LB,a+t)$ is
$\M_S(0,\LB,a-b)$ which is also $\Pic_{\C}^{a-b+g-1}$. 
The evaluation homomorphism in (\ref{eq-E-U-F}) is injective and 
$E(\tilde{F},U)$ is the zero sheaf. 
We have a short exact sequence
\begin{equation}
\label{eq-U-F-Q}
0 \rightarrow U\otimes\StructureSheaf{S} 
\rightarrow F 
\rightarrow Q 
\rightarrow 0.
\end{equation}
The quotient sheaf $Q$ is a torsion $\StructureSheaf{S}$-module 
which is supported, as a rank $1$ torsion free sheaf, on a curve $C$
in the linear system $\linsys{\LB}$ (see Lemma
\ref{lazarsfeld-slope-stability-lemma}
Part \ref{lemma-item-stability-of-quotient-with-support-on-a-curve}). 
In particular, $Q$ represents a sheaf in
$\M_S(0,\LB,b-a)$ which is $\Pic_{\C}^{b-a+g-1}$. 
As in the case $(a,b)=(1,0)$ and $t=0$, we carry out the reflection on
the curve $C$ supporting $Q$. We define

\begin{equation}
\label{eq-reflection-on-a-curve}
q(F,U) := \SheafHom_C(Q,\omega_C). 
\end{equation}

\noindent
The exact sequence (\ref{eq-U-F-Q}) corresponds to an embedding 
$i:U^* \hookrightarrow \Ext^1_S(Q,\StructureSheaf{S})$. 
Serre's Duality on $S$ represents $U$ as a quotient of 
$H^1(S,Q\otimes \omega_S)$. Serre's Duality on $C$ represents $U^*$ as a
subset $W$ of $H^0(q(F,U))$. We define

\begin{equation}
\label{eq-lazarsfelds-reflection-of-a-pair-on-a-curve}
\tilde{q}_t(F,U) := (q(F,U),W). 
\end{equation}

\medskip 
Lemma \ref{lemma-hom-over-C-is-ext1-over-S} 
identifies the reflection (\ref{eq-reflection-on-a-curve}) as an operation on 
$\StructureSheaf{S}$-modules: we have an equality
$q(F,U) = \SheafExt^1_S(Q,\StructureSheaf{S})$.

\bigskip
\noindent
{\bf The case $a+b+t > a$ revisited:} 
The reflection (\ref{eq-lazarsfelds-reflection-of-a-pair-rk-U-larger-rk-F}),
while canonical, is inconvenient to work with in families. 
Upon a choice of an $a$-dimensional  subspace $W\subset U$, we can express
(\ref{eq-lazarsfelds-reflection-of-a-pair-rk-U-larger-rk-F}) 
as the conjugation of the reflection  along a curve
(\ref{eq-lazarsfelds-reflection-of-a-pair-on-a-curve}) by the extension
isomorphism in Theorem \ref{thm-the-tyurin-extension-isomorphism}. 
We need  (\ref{eq-lazarsfelds-reflection-of-a-pair-rk-U-larger-rk-F}) 
in order to prove that the conjugation is independent off the choice of $W$. 

Let $W$ be an $a$-dimensional  subspace of $U$ and assume that 
$a+b+t > a$. Then $Q:=F/(W\otimes\StructureSheaf{S})$ is a stable sheaf with 
pure one-dimensional support (depending on $W$). 
$U/W$ is a subspace of $H^0(Q)$ and
$W^*$ is a subspace of 
$\Ext^1_S(Q,\StructureSheaf{S})$. 
Consequently, $W^*$ is a subspace of 
$H^0(\SheafExt^1_S(Q,\StructureSheaf{S}))$
and $(U/W)$ is a subspace of 
$\Ext^1_S[\SheafExt^1_S(Q,\StructureSheaf{S}),
\StructureSheaf{S}]$. Theorem \ref{thm-the-tyurin-extension-isomorphism}
implies that $(U/W)^*$ is a subspace of global sections of the stable sheaf 
$f[\SheafExt^1_S(Q,\StructureSheaf{S}),U/W]$
extending $\SheafExt^1_S(Q,\StructureSheaf{S})$
by $(U/W)^*\otimes \StructureSheaf{S}$. We get that $U^*$ is isomorphic to 
the inverse image in 
$H^0(f[\SheafExt^1_S(Q,\StructureSheaf{S}),U/W])$ of $W^*$ 
under the quotient sheaf homomorphism 
$f[\SheafExt^1_S(Q,\StructureSheaf{S}),U/W] \rightarrow 
\SheafExt^1_S(Q,\StructureSheaf{S})$.

\begin{new-lemma}
\label{lemma-lazarsfelds-reflection-is-a-conjugation-of-reflection-along-curve}
The sheaves $q(F,U)$ and $f(q(f^{-1}(F,W),U/W),(U/W)^*)$ are isomorphic and
the isomorphism conjugates the embeddings of 
$U^*$ as a subspace of global sections. 
\end{new-lemma}

\noindent
{\bf Proof:}
In the special case where $U$ generates $F$ and $F$ is locally free, 
the Lemma follows from the commutative diagram with short exact rows and 
columns
\[
\begin{array}{ccccc}
0 & \hookrightarrow & (U/W)^*\otimes\StructureSheaf{S} & \RightArrowOf{=} &
(U/W)^*\otimes\StructureSheaf{S} 
\\
\downarrow & & \downarrow & & \downarrow 
\\
F^* &  \hookrightarrow & U^*\otimes\StructureSheaf{S} & \rightarrow & q(F,U)
\\
=\downarrow\hspace{1em} & & \downarrow & & \downarrow 
\\
F^* & \hookrightarrow & W^*\otimes\StructureSheaf{S} & \rightarrow & 
\SheafExt^1_S(Q,\StructureSheaf{S}) 
\end{array}
\]
The proof of the general case is a laborious unwinding of the definition
of the reflection 
(\ref{eq-lazarsfelds-reflection-of-a-pair-rk-U-larger-rk-F}). 
\EndProof

\bigskip
\noindent
{\bf Step II}: (of the proof Theorem 
\ref{thm-lazarsfeld-reflection-isomorphism})
We work out the relative version of the maps $\tilde{q}_t$ in 
(\ref{eq-lazarsfelds-reflection-of-a-pair-rk-U-larger-rk-F})
and
(\ref{eq-lazarsfelds-reflection-of-a-pair-on-a-curve}). 
Consequently, the functors, which are coarsely
represented by $G^0(k,\M_S(v))$ and $G^0(k,\M_S(\sigma\circ\tau(v-\vec{k})))$
are equivalent. 

Set $r:=\rank(v)$. 
The case $r=k$ was proven by Le Potier 
(Theorem 5.7 in \cite{le-potier-coherent} 
when $k=r=0$ and Theorem 5.12 when $k=r>0$). 
The general case $k\geq r$ is a conjugate of the case $r=k=0$ 
as in Lemma 
\ref{lemma-lazarsfelds-reflection-is-a-conjugation-of-reflection-along-curve}. 
Assume given 
1) a family $\F_v$ over $S\times T$ flat over a scheme $T$, 
2) a locally free $\StructureSheaf{T}$-module $\U$ of rank $k$ and 
3) a homomorphism $i: \U \hookrightarrow p_*(\F_v)$ 
injective on each fiber. 
Choose an open covering $\{T_j\}$ of $T$ and local sections 
of $G(r,\U)$ corresponding to subbundles $\W_j$ of 
$\restricted{\U}{T_j}$. By Theorem 
\ref{thm-the-tyurin-extension-isomorphism}, 
we have an equivalence of functors 
$G^0(r,\M_S(v)) \cong G_1(r,\M_S(v-\vec{r}))$. 
We get a short exact sequence flat over each $T_j$
\[
0\rightarrow \W_j \hookrightarrow \F_v \rightarrow 
\F_{v-\vec{r},j}
\rightarrow 0
\]
and homomorphisms, injective on each fiber,
\begin{eqnarray}
\label{eq-W-dual-is-a-local-subsheaf}
\W^* & \hookrightarrow & 
\RelExt^1_p(\F_{v-\vec{r},j},\StructureSheaf{T_j\times S})
\\
\label{eq-U-mod-W-is-a-local-subsheaf}
\U/\W_j & \hookrightarrow & p_*\F_{v-\vec{r},j}.
\end{eqnarray}

\noindent
Note that if $k>r$ the families $\F_{v-\vec{r},j}$ need not patch. 
(Even their support, which is of relative dimension 1, need not patch!) 
Applying Le Potier's Theorem  5.7 we get families 
\[
\F_{\sigma\circ\tau(v-\vec{r}),j}:= 
\SheafExt^1_{T_j\times S}(\F_{v-\vec{r},j},\StructureSheaf{T_j\times S})
\]
flat over $T_j$ and the analogues of (\ref{eq-W-dual-is-a-local-subsheaf}) 
and (\ref{eq-U-mod-W-is-a-local-subsheaf}).
Applying Theorem 
\ref{thm-the-tyurin-extension-isomorphism} once more we get 
families $\F_{\sigma\circ\tau(v-\vec{k}),j}$ of rank $k\!-\!r$ stable sheaves
and natural homomorphisms, injective on each fiber,
\[
\iota_j \ : \ 
(\U^*\restricted{)}{T_j} \ \hookrightarrow \ 
p_*\F_{\sigma\circ\tau(v-\vec{k}),j}.
\]
We claim that the the families $\F_{\sigma\circ\tau(v-\vec{k}),j}$ 
and the homomorphisms $\iota_j$ 
patch naturally to a global family $\F_{\sigma\circ\tau(v-\vec{k})}$ 
and a global homomorphism $\iota$. Indeed, Lemma 
\ref{lemma-lazarsfelds-reflection-is-a-conjugation-of-reflection-along-curve}
and the simplicity of the sheaves parametrized imply that 
the the sheaf 
$\ p_*\SheafHom(\F_{\sigma\circ\tau(v-\vec{k}),i},
\F_{\sigma\circ\tau(v-\vec{k}),j})$ is a 
line bundle on $T_i\cap T_j$ and it has a canonical invertible section 
$\phi_{i,j}$ satisfying $\phi_{i,j}\circ \iota_j=\iota_i$. 
The collection $\{\phi_{i,j}\}$ is a $1$-cocycle which
glues $\{\F_{\sigma\circ\tau(v-\vec{k}),j}\}$ because its restriction 
to the invariant subsheaves $\restricted{\U^*}{T_j}$ is.  
This completes the proof of Theorem 
\ref{thm-lazarsfeld-reflection-isomorphism}.
\EndProof

\begin{rem}
\label{rem-proof-of-dualizability-for-negative-euler-characteristic}
{\rm
We can now indicate the modifications necessary in the statement and 
proof of the analogue of Proposition 
\ref{prop-grassmannian-fibrations-of-blown-up-brill-noether-loci}
in case $\rank(v)\geq 0$ and $\chi\leq 0$. Recall that the Brill-Noether
stratification is indexed now by $h^0$. 
In the statement of parts 
\ref{prop-item-constant-rank} and 
\ref{prop-item-dual-of-cokernel-is-subsheaf-of-ext} 
the rank of $W^k_{v'}$ is $k$ and the rank of $U^k_{v'}$ is $k-\chi(v')$. 
In the statement of part \ref{prop-item-blown-up-tyurin-isomorphism}
we need to replace $G(k,W_{v'+\vec{k}})$ by $G(k,U^*_{v'\!-\!\vec{k}})$. 
In this case, even though $\rank(v)$ is non-negative, we may end-up with a 
non-empty Brill-Noether locus $\M_S(v)^k$ such that $v\!-\!\vec{k}$ has 
negative rank. If $k > \rank(v)$, we define $U_{v\!-\!\vec{k}}$ to be 
$W^*_{\sigma\circ\tau(v\!-\!\vec{k})}$ and $W_{v\!-\!\vec{k}}$ to be 
$U^*_{\sigma\circ\tau(v\!-\!\vec{k})}$. 
We also define $G_1(k,\M_S(v\!-\!\vec{k}))$ to be 
$G^0(k,\M_S(\sigma\circ\tau(v\!-\!\vec{k})))$ if $k > \rank(v)$. 

In the proof of part \ref{prop-item-blown-up-tyurin-isomorphism}
we need to replace the isomorphism 
(\ref{eq-this-isomorphism-f-will-be-replaced-in-the-negative-euler-char-case})
by
\[
f \ : \ G^0(k,\M_S(v)) \ \LongIsomRightArrow \ G_1(k,\M_S(v\!-\!\vec{k})).
\]
If $k\leq \rank(v)$, this isomorphism follows from Theorem 
\ref{thm-the-tyurin-extension-isomorphism}. 
If $k> \rank(v)$, the above definitions translate it to the isomorphism
\[
G^0(k,\M_S(v)) \ \cong \ G^0(k,\M_S(\sigma\circ\tau(v\!-\!\vec{k}))).
\]
The latter isomorphism is precisely Theorem 
\ref{thm-lazarsfeld-reflection-isomorphism}.
}
\end{rem}

\subsection{The collections are dual} 
\label{sec-dual-collections}

We prove in this section that the two dualizable collections
in Theorem \ref{thm-mukai-reflection-extends-to-a-stratified-transformation}
are dual to each other. This completes the proof of the Theorem. 

\begin{prop}
\label{prop-dualizable-collections-are-dual}
Let $v=(r,\LB,s)$ be a Mukai vector and 
$\LB$ a line bundle satisfying 
Condition \ref{cond-linear-system}. 
Then there exist  natural isomorphisms 
\begin{eqnarray}
\label{eq-lazarsfelds-reflection-of-complete-colineations}
\tilde{q}: B^{[1]}\M_S(v)  & \IsomRightArrow & B^{[1]}\M_S(\sigma(v)), 
\ \ \ \mbox{and} 
\\
\label{eq-tyurins-reflection-of-complete-colineations}
\tilde{q}: B^{[1]}\M_S(v)  & \IsomRightArrow & B^{[1]}\M_S(\tau(v))
\end{eqnarray}
compatible with respect to the Brill-Nother loci.
\end{prop}

\noindent
{\bf Proof:} 
Proposition \ref{prop-dualizable-collections-are-dual} follows 
from the universal properties of the coarse moduli spaces involved. 
It suffices to construct the isomorphism 
(\ref{eq-lazarsfelds-reflection-of-complete-colineations}) 
in case $v=(a,\LB,b)$ and both $a$ and $b$ are non-negative integers. 
Proposition \ref{prop-grassmannian-fibrations-of-blown-up-brill-noether-loci}
part \ref{prop-item-constant-rank} produces the pair 
$(B^{[1]}\M_S(v), {\PP}W_v)$. 
By the universal property of $G^0(\chi(v),\M_S(v))$ we get a morphism 
\[
B^{[1]}\M_S(v) \ \rightarrow \ G^0(\chi(v),\M_S(v))
\]
(see Proposition \ref{prop-universal-propery-of-G0}). Similarly, we get
a morphism 
\[
B^{[1]}\M_S(\sigma(v))  \ \rightarrow \ G^0(\chi(\sigma(v)),\M_S(\sigma(v))). 
\]
By Theorem \ref{thm-lazarsfeld-reflection-isomorphism}, 
the two moduli spaces  
$G^0(\chi(v),\M_S(v))$ and $G^0(\chi(\sigma(v)),\M_S(\sigma(v)))$ 
are isomorphic. 
Part \ref{prop-item-universal-property-of-G-k-W}
of Proposition 
\ref{prop-grassmannian-fibrations-of-blown-up-brill-noether-loci} 
implies that there exist  morphisms 
\begin{eqnarray*}
B^{[1]}\M_S(v) & \rightarrow & G(0,W_{\sigma(v)}) \ = \ B^{[1]}\M_S(\sigma(v)),
\ \ \ \mbox{and} 
\\
B^{[1]}\M_S(\sigma(v))  & \rightarrow & G(0,W_{v}) \ = \ B^{[1]}\M_S(v).
\end{eqnarray*}
The composition of the two morphisms is generically the identity
morphism (by Corollary \ref{cor-verification-of-codimension-condition}). 
It must be globally the identity since 
$B^{[1]}\M_S(v)$ and $B^{[1]}\M_S(\sigma(v))$ are smooth
(part \ref{prop-item-smoothness} of Proposition 
\ref{prop-grassmannian-fibrations-of-blown-up-brill-noether-loci}). 

\bigskip
\noindent
It suffices to construct the isomorphism 
(\ref{eq-tyurins-reflection-of-complete-colineations})
in case $a\geq 0$ and $b\leq 0$. 
The construction is similar to that of 
(\ref{eq-lazarsfelds-reflection-of-complete-colineations}).
Simply use Theorem \ref{thm-the-tyurin-extension-isomorphism} instead of
Theorem \ref{thm-lazarsfeld-reflection-isomorphism}. 
\EndProof

\end{document}